\documentclass[10pt]{amsart}
\usepackage{amssymb}
\usepackage{amscd}

\theoremstyle{definition}
\newtheorem{thm}{Theorem}[section]
\newtheorem{lem}[thm]{Lemma}
\newtheorem{prp}[thm]{Proposition}
\newtheorem{dfn}[thm]{Definition}
\newtheorem{cor}[thm]{Corollary}

\newtheorem{cnv}[thm]{Convention}
\newtheorem{rmk}[thm]{Remark}
\newtheorem{ntn}[thm]{Notation}
\newtheorem{exa}[thm]{Example}

\newtheorem{qst}[thm]{Question}

\newcommand{\beq}{\begin{equation}}
\newcommand{\eeq}{\end{equation}}
\newcommand{\beqr}{\begin{eqnarray*}}
\newcommand{\eeqr}{\end{eqnarray*}}
\newcommand{\bal}{\begin{align*}}
\newcommand{\eal}{\end{align*}}
\newcommand{\bei}{\begin{itemize}}
\newcommand{\eei}{\end{itemize}}
\newcommand{\limi}[1]{\lim_{{#1} \to \infty}}

\newcommand{\af}{\alpha}
\newcommand{\bt}{\beta}
\newcommand{\gm}{\gamma}
\newcommand{\dt}{\delta}
\newcommand{\ep}{\varepsilon}
\newcommand{\zt}{\zeta}
\newcommand{\et}{\eta}
\newcommand{\ch}{\chi}
\newcommand{\io}{\iota}
\newcommand{\te}{\theta}
\newcommand{\ld}{\lambda}
\newcommand{\sm}{\sigma}

\newcommand{\ph}{\varphi}
\newcommand{\ps}{\psi}
\newcommand{\rh}{\rho}
\newcommand{\om}{\omega}
\newcommand{\ta}{\tau}

\newcommand{\Gm}{\Gamma}

\newcommand{\Sm}{\Sigma}

\newcommand{\Q}{{\mathbf{Q}}}
\newcommand{\Z}{{\mathbf{Z}}}
\newcommand{\R}{{\mathbf{R}}}
\newcommand{\C}{{\mathbf{C}}}
\newcommand{\N}{{\mathbf{N}}}

\newcommand{\id}{{\mathrm{id}}}
\newcommand{\ev}{{\mathrm{ev}}}
\newcommand{\sint}{{\mathrm{int}}}

\newcommand{\dist}{{\mathrm{dist}}}
\newcommand{\sa}{{\mathrm{sa}}}
\newcommand{\spec}{{\mathrm{sp}}}

\newcommand{\diag}{{\mathrm{diag}}}
\newcommand{\supp}{{\mathrm{supp}}}
\newcommand{\rank}{{\mathrm{rank}}}

\newcommand{\spn}{{\mathrm{span}}}
\newcommand{\card}{{\mathrm{card}}}
\newcommand{\Aut}{{\mathrm{Aut}}}
\newcommand{\Ad}{{\mathrm{Ad}}}

\newcommand{\Hom}{{\mathrm{Hom}}}
\newcommand{\Ext}{{\mathrm{Ext}}}
\newcommand{\Ker}{{\mathrm{Ker}}}

\newcommand{\Zq}[1]{{\Z_{#1}}}
\newcommand{\Zqt}{\Zq{2}}
\newcommand{\Zqh}{\Zq{3}}
\newcommand{\Zqf}{\Zq{4}}
\newcommand{\Zqs}{\Zq{6}}
\newcommand{\Zqn}{\Zq{n}}
\newcommand{\Cs}[3]{C^* (\Zq{#1}, #2, #3)}
\newcommand{\Csw}[3]{C^* (\Zq{#1}, \, #2, \, #3)}
\newcommand{\CZnAa}{\Cs{n}{A}{\af}}

\newcommand{\dirlim}{\displaystyle \lim_{\longrightarrow}}

\newcommand{\andeqn}{\,\,\,\,\,\, {\mbox{and}} \,\,\,\,\,\,}
\newcommand{\QED}{\rule{0.4em}{2ex}}
\newcommand{\ts}[1]{{\textstyle{#1}}}
\newcommand{\ds}[1]{{\displaystyle{#1}}}
\newcommand{\ssum}[2]{{\ts{ {\ds{\sum}}_{#1}^{#2} }}}

\newcommand{\ca}{C*-algebra}

\newcommand{\suca}{simple unital C*-algebra}
\newcommand{\pisca}{purely infinite simple C*-algebra}
\newcommand{\ct}{continuous}
\newcommand{\pj}{projection}

\newcommand{\nbhd}{neighborhood}

\newcommand{\hm}{homomorphism}
\newcommand{\wolog}{without loss of generality}
\newcommand{\Wolog}{Without loss of generality}

\newcommand{\ifo}{if and only if}

\newcommand{\mops}{mutually orthogonal \pj s}

\newcommand{\hme}{homeomorphism}

\newcommand{\cfn}{continuous function}
\newcommand{\hsa}{hereditary subalgebra}

\newcommand{\mvnt}{Murray-von Neumann equivalent}
\newcommand{\mvnc}{Murray-von Neumann equivalence}
\newcommand{\aRp}{tracial Rokhlin property}
\newcommand{\sRp}{strict Rokhlin property}
\newcommand{\fd}{finite dimensional}
\newcommand{\uct}{Universal Coefficient Theorem}
\newcommand{\suct}{satisfies the \uct}

\newcommand{\rsz}[1]{\raisebox{0ex}[0.8ex][0.8ex]{$#1$}}

\renewcommand{\S}{\subset}
\newcommand{\ov}{\overline}
\newcommand{\SM}{\setminus}
\newcommand{\I}{\infty}

\title[Crossed products and tracial Rokhlin]{Crossed
     products by finite cyclic group actions with the tracial
     Rokhlin property}

\author{N.\  Christopher Phillips}

\date{26 June 2003}

\address{Department of Mathematics, University  of Oregon,
       Eugene OR 97403-1222, USA.}

\email[]{ncp@darkwing.uoregon.edu}

\subjclass{Primary 46L35, 46L55;
 Secondary 22D25, 46L40.}
\thanks{Research partially supported by NSF grant DMS 0070776.}

\begin{document}

\setcounter{section}{-1}

\begin{abstract}
We define the \aRp\  for actions of finite cyclic groups
on stably finite simple unital \ca s.
We prove that the crossed product
of a stably finite simple unital \ca\  with tracial rank zero
by an action with this property again has tracial rank zero.
Under a kind of weak approximate innerness assumption and
one other technical condition, we prove that if the
action has the \aRp\  and the crossed product has tracial rank zero,
then the original algebra has tracial rank zero.
We give examples showing how the tracial Rokhlin property differs
from the Rokhlin property of Izumi.

We use these results,
together with work of Elliott-Evans and Kishimoto,
to give an inductive proof that
every simple higher dimensional noncommutative torus is an AT~algebra.
We further prove that the crossed product of
every simple higher dimensional noncommutative torus by the flip
is an AF algebra,
and that the crossed products of irrational rotation algebras
by the standard actions of $\Zqh$, $\Zqf$, and $\Zqs$ are
simple AH algebras with real rank zero.
In the case of $\Zqf$, we recover Walters' result that the
crossed product is AF for a dense $G_{\dt}$-set of rotation numbers.
\end{abstract}

\maketitle

\section{Introduction}\label{Sec:Intro}

\indent
A higher dimensional noncommutative torus is a generalization
of the rotation algebra $A_{\te}$ to more generators.
See Section~\ref{Sec:NCT} for details.
The most important result of this paper is that
every simple higher dimensional noncommutative torus
is an AT~algebra, that is, a direct limit of finite direct sums
of \ca s of the form $C (S^1, M_n)$ for varying values of $n$.
Elliott and Evans proved~\cite{EE} that the irrational rotation
algebras are AT~algebras,and Boca~\cite{Bc}
showed that ``most'' simple higher dimensional noncommutative toruses
are AT~algebras.
See the introduction to Section~\ref{Sec:NCT} for more of the history.

Our proof is by induction on the number of generators.
Every higher dimensional noncommutative torus can be written
as a successive crossed product by $\Z$,
and the proof of Corollary~6.6 of~\cite{Ks4} uses an inductive
argument which works whenever the intermediate crossed
products are all simple.
Our contribution is a method,
involving crossed products by actions of finite cyclic groups
which have what we call the \aRp,
for replacing any simple higher dimensional noncommutative torus
by one to which something close to this argument applies.
Our proof depends heavily on H.\  Lin's classification theorem
for simple \ca s with tracial rank zero~\cite{Ln15}.

Our original motivation for the \aRp\  is that it enables us to
prove that every simple higher dimensional noncommutative torus
is an AT~algebra.
In retrospect, however, the following motivation is perhaps better.
In~\cite{Iz} and~\cite{Iz2},
Izumi has started an intensive study of finite group actions with
the Rokhlin property, which here, to minimize confusion,
we call the strict Rokhlin property.
The strict Rokhlin property imposes severe restrictions on the
relation between the K-theory of the original algebra,
the action of the group on this K-theory,
and the K-theory of the crossed product.
See especially Section~3 of~\cite{Iz2}.
In particular, the results there show that none of the actions
appearing in any of our main applications
has the strict Rokhlin property.
Accordingly, a less restrictive version of the strict Rokhlin property
is needed.
The definition of the \aRp\  is obtained from the definition of
the strict Rokhlin property in roughly the same way
that H.\  Lin's definition of a tracially AF \ca~\cite{LnTAF}
is obtained from the definition of an AF algebra.
We note that a similar concept, for integer actions on AF~algebras,
has appeared in~\cite{BEK}.

There are standard actions of $\Zqt$, $\Zqh$, $\Zqf$, and $\Zqs$
on the rotation algebras,
and the action of $\Zqt$ is defined on every
higher dimensional noncommutative torus as well.
See the introductions to Sections~\ref{Sec:NCFT}
and~\ref{Sec:CrPrdFlip} for details and some history.
When the algebra is simple, all these actions have the \aRp.
As a consequence, we show that the crossed products are always
simple AH~algebras with no dimension growth and real rank zero.
For the action of $\Zqt$ on an irrational rotation algebra,
it was already known~\cite{BK} that the crossed is~AF.
For a dense $G_{\dt}$-set of rotation numbers $\te$,
the computation of K-theory of the crossed product
$C^* (\Zqf, A_{\te})$ implies that this algebra is~AF.
We thus recover the main result of~\cite{Wl2}.
That paper, however, gives no information on the structure
of $C^* (\Zqf, A_{\te})$ for other irrational values of $\te$,
while our results imply in particular that $C^* (\Zqf, A_{\te})$
has real rank zero and stable rank one for all irrational $\te$.
Our results are completely new for the
actions of $\Zqh$ and $\Zqs$ on the irrational rotation algebras.
In the case of the flip on a
simple higher dimensional noncommutative torus,
the crossed product was already known to be~AF in
``most'' cases~\cite{Bc};
we are able to compute the K-theory and show that the
crossed product is always~AF.

Roughly the first half of this paper is a direct line to
the proof that every simple noncommutative torus is an AT~algebra.
We develop the theory of crossed products by actions with
the \aRp\  just far enough for this application.
We start the second half of the paper with a theorem,
Theorem~\ref{ARPFromPosElts2}, which we find very useful
for proving that an action has the \aRp.
Then we consider the standard actions of $\Zqh$, $\Zqf$, and $\Zqs$
on the irrational rotation algebras,
and the flip action of $\Zqt$
on a simple higher dimensional noncommutative torus.
Finally, we return to the general theory,
and in particular present several examples involving UHF algebras
which illustrate the limits of the theorems.

We now give a brief outline of the sections.
Further introductory material can be found in
Sections~\ref{Sec:ARP},
\ref{Sec:NCT}, \ref{Sec:NCFT}, and~\ref{Sec:Exs}.
We start in Section~\ref{Sec:ARP}
by introducing the \aRp\  for actions of finite cyclic groups.
In Section~\ref{Sec:CrPrd}, we prove that crossed products of \ca s with
tracial rank zero by such actions again have tracial rank zero.
In Section~\ref{Sec:Dual}, we give some conditions under which the
dual action of an action with the \aRp\  again has the \aRp, and
we examine the fixed point algebra of an action with the \aRp.
Section~\ref{Sec:TAI} introduces the notion of a
tracially approximately inner automorphism,
which is needed in Section~\ref{Sec:Dual}.

Then we turn to higher dimensional noncommutative toruses.
Section~\ref{Sec:NCT} contains various preliminaries.
In Section~\ref{Sec:ARokhOnNCT}, we prove that if $A$ is a simple
noncommutative torus, then the automorphism which multiplies one
of the standard unitary generators by $\exp (2 \pi i / n)$
generates an action  of $\Zqn$ with the \aRp.
In Section~\ref{Sec:NCTInd}, we use this result, work of Kishimoto
on crossed products by $\Z$, and H.\  Lin's classification theorem,
to construct an inductive proof that every simple noncommutative torus
is an AT~algebra.

In the next three sections, we consider crossed products of
irrational rotation algebra
and higher dimensional noncommutative toruses
by actions of finite cyclic groups which have already been studied.
In Section~\ref{Sec:NCFT},
we show that the noncommutative Fourier transform
on an irrational rotation algebra generates an action of
$\Zqf$ with the \aRp,
and use this to show that the crossed product
is always a simple AH algebra with real rank zero,
and is AF for ``many'' rotation numbers.
This section contains a very useful criterion for an action
on a simple \ca\  with tracial rank zero with unique tracial state
to have the \aRp,
a criterion which does not mention any \pj s.
In Section~\ref{Sec:OtherActions}, we show that the crossed
products of an irrational rotation algebra by the ``standard''
actions of $\Zqh$ and $\Zqs$ are simple AH algebras with real rank zero.
Then in Section~\ref{Sec:CrPrdFlip} we show that the crossed product
of a simple higher dimensional noncommutative torus by the flip
is always AF.

In the last three sections, we return to the general theory.
Section~\ref{Sec:More} contains several results which were not
proved earlier because they are not needed for our main applications.
In particular, the crossed product an AF algebra
by an action with the (strict) Rokhlin property is again AF.
In Section~\ref{Sec:Exs},
we give examples of actions on the $2^{\infty}$ UHF algebra
which show that many of our results can't be improved.
Section~\ref{Sec:Qst} contains some interesting questions
to which we do not know the answers.

We use the notation $\Zq{n}$ for $\Z / n \Z$;
the $p$-adic integers will not appear in this paper.
If $A$ is a \ca\  and $\af \colon A \to A$ is an automorphism
such that $\af^n = \id_A$, then we write
$\CZnAa$ for the crossed product of
$A$ by the action of $\Zqn$ generated by $\af$.

We write $A_{\sa}$ for the set of selfadjoint elements of a \ca\  $A$.
We write $p \precsim q$ to mean that the \pj\  $p$ is \mvnt\  to a
sub\pj\  of $q$, and $p \sim q$ to mean that
$p$ is \mvnt\  to $q$.
Also, $[a, b]$ denotes the additive commutator $a b - b a$.

I am grateful to Masaki Izumi for discussions concerning the
(strict) Rokhlin property.
In particular, he suggested Example~\ref{CAR3}.
I am grateful to Marc Rieffel for discussions concerning
higher dimensional noncommutative toruses,
and in particular for pointing out that it was not known
whether a higher dimensional noncommutative torus is isomorphic to
its opposite algebra.
See Corollary~\ref{IsomOpp};
this remains open in the nonsimple case.
I am also grateful to to Sam Walters for discussions related to the
noncommutative Fourier transform,
especially his work on the K-theory of the crossed product by this
action.
Finally, I would like to thank Hiroyuki Osaka for carefully
reading earlier versions of this paper, and catching a number
of misprints and suggesting many improvements,
and Hanfeng Li for useful comments.

\section{The tracial Rokhlin property}\label{Sec:ARP}

\indent
In this section we introduce the \aRp\  and several related
properties.
We observe several elementary relations and consequences,
and we prove several useful equivalent formulations.

\begin{dfn}\label{ARPDfn}
Let $A$ be a stably finite \suca\  %
and let $\af \in \Aut (A)$ satisfy $\af^n = \id_A$.
We say that the action of $\Zqn$ generated by $\af$ has the
{\emph{tracial Rokhlin property}} if for every finite set
$F \S A$, every $\ep > 0$, every $N \in \N$,
and every nonzero positive element $x \in A$,
there are \mops\  $e_0, e_1, \dots, e_{n - 1} \in A$ such that:
\begin{itemize}
\item[(1)]
$\| \af (e_j) - e_{j + 1} \| < \ep$ for $0 \leq j \leq n - 2$.
\item[(2)]
$\| e_j a - a e_j \| < \ep$ for $0 \leq j \leq n - 1$ and all $a \in F$.
\item[(3)]
With $e = \sum_{j = 0}^{n - 1} e_j$, the \pj\  $1 - e$ is \mvnt\  to a
\pj\  in the \hsa\  of $A$ generated by $x$.
\item[(4)]
For every $j$ with $0 \leq j \leq n - 1$,
there are $N$ \mops\  $f_1, f_2, \dots, f_N \leq e_j$,
each of which is \mvnt\  to the \pj\  $1 - e$ of~(3).
\end{itemize}
\end{dfn}

In this definition, $e = 1$ is allowed, in which case conditions~(3)
and~(4) are vacuous.

\begin{cnv}\label{ARPConv}
With the notation as in Definition~\ref{ARPDfn}, we always take
$e_n = e_0$, and we will usually take $e = \sum_{j = 0}^{n - 1} e_j$.
\end{cnv}

We do not require that $\sum_{j = 0}^{n - 1} e_j = 1$,
as Izumi does for the Rokhlin property in Definition~3.1 of~\cite{Iz}.
The terminology is motivated by H.\  Lin's
``tracially AF''~\cite{LnTAF} and
``tracial topological rank''~\cite{LnTTR},
in whose definitions there is also a ``small'' leftover projection.

We recall Izumi's definition,
specialized to the case of a finite cyclic group, but to emphasize
the difference, we call it the strict Rokhlin property here.

\begin{dfn}\label{ERPDfn}
Let $A$ be a unital \ca\  and
let $\af \in \Aut (A)$ satisfy $\af^n = \id_A$.
We say that the action of $\Zqn$ generated by $\af$ has the
{\emph{strict Rokhlin property}} if for every finite set
$F \S A$, and every $\ep > 0$,
there are \mops\  $e_0, e_1, \dots, e_{n - 1} \in A$ such that:
\begin{itemize}
\item[(1)]
$\| \af (e_j) - e_{j + 1} \| < \ep$ for $0 \leq j \leq n - 2$.
\item[(2)]
$\| e_j a - a e_j \| < \ep$ for $0 \leq j \leq n - 1$ and all $a \in F$.
\item[(3)]
$\sum_{j = 0}^{n - 1} e_j = 1$.
\end{itemize}
\end{dfn}

\begin{rmk}\label{SRPImpTRP}
If an action of $\Zqn$ on a unital \ca\  $A$ has the
strict Rokhlin property, then it has the \aRp.
\end{rmk}

If $\af$ is approximately inner,
requiring $\sum_{j = 0}^{n - 1} e_j = 1$ forces $[1_A] \in K_0 (A)$
to be divisible by $n$,
and therefore rules out many \ca s of interest.
In fact, the strict Rokhlin property imposes much more stringent
conditions on the K-theory.
Theorem~3.3 and Lemma~3.2(1) of~\cite{Iz2} show that if
a nontrivial finite group $G$ acts on a simple unital \ca\  $A$
in such a way that the induced action on $K_* (A)$ is trivial,
and if one of $K_0 (A)$ and $K_1 (A)$ is a nonzero free abelian
group, then $\af$ does not have the strict Rokhlin property.
Theorem~3.3 and the discussion preceding Theorem~3.4 of~\cite{Iz2}
show that if in addition $G$ is cyclic of order $n$, then
the strict Rokhlin property implies that
$K_* (A)$ is uniquely $n$-divisible.
It follows that the actions considered in the applications in
this paper never have the strict Rokhlin property.

For the \aRp\  to be likely to hold, the \ca\  must have a
reasonable number of \pj s.
For reference, we recall here the definition of the property that
seems most relevant.

\begin{dfn}\label{SPD}
Let $A$ be a \ca.
We say that $A$ has {\emph{Property~(SP)}} if every nonzero
\hsa\  in $A$ contains a nonzero \pj.
\end{dfn}

\begin{lem}\label{RImpSP}
Let $A$ be a stably finite \suca,
and let $\af \in \Aut (A)$ be an automorphism which
satisfies $\af^n = \id_A$ and such that the action of $\Zqn$
generated by $\af$ has the tracial Rokhlin property.
Then $A$ has Property~(SP) or the action
generated by $\af$ has the strict Rokhlin property.
\end{lem}

\begin{proof}
If $A$ does not have Property~(SP),
then there is a nonzero positive element $x \in A$
which generates a \hsa\  which contains no nonzero \pj.
\end{proof}

Our definition also seems too weak for use with \pisca s.
For example,
Condition~(3) of Definition~\ref{ARPDfn} would then be vacuous,
and Condition~(4) would be automatically
satisfied whenever $e \neq 0$.
On the other hand, the strict Rokhlin property still seems too strong,
since even in a unital \pisca,
$[1]$ need not be divisible by $n$ in $K_0 (A)$.

For the \aRp\  as we have defined it to be useful,
not only stable finiteness but also
some condition on comparison of \pj s seems to be necessary.
Although we will not make immediate use of them,
we recall here two well known conditions of this type.
Others will appear later.

\begin{dfn}\label{CancD}
Let $A$ be a unital \ca.
We say that $A$ has {\emph{cancellation of \pj s}} if
whenever $n \in \N$ and $p, q, e \in M_n (A)$ are \pj s
such that $e$ is orthogonal to both $p$ and $q$,
and $p + e \sim q + e$, then $p \sim q$.
\end{dfn}

\begin{dfn}\label{OrdDetD}
Let $A$ be a unital \ca.
We say that the {\emph{order on \pj s over $A$ is determined by traces}}
if whenever $n \in \N$ and $p, q \in M_n (A)$ are \pj s such that
$\ta (p) < \ta (q)$ for all tracial states $\ta$ on $A$, then
$p \precsim q$.
\end{dfn}

This is just Blackadar's Second Fundamental Comparability Question.
See 1.3.1 in~\cite{Bl3}.

Finally, we do not attempt to formulate the appropriate definition for
actions on \ca s which are not simple.
At the very least, one should ask that, given an arbitrary
nonzero element $a \in A$, one can in addition require
$\| e a e \| \geq \| a \| - \ep$.
See Definition~2.1 of~\cite{LnTAF}.
Quite possibly one should impose other conditions as well.

It is convenient to have a formally stronger version of the \aRp,
in which $\| \af (e_{n - 1}) - e_0 \| < \ep$ as well, and in
which the defect \pj\  is $\af$-invariant.
We will use a simple case of the following lemma,
which is stated more generally for later use.
We note for emphasis: in it, $\dt$ depends on $E$, $N$, and $\ep$,
but not on $\ph$ or the particular \pj s $p, e_0, e_1, \dots, e_N$.

\begin{lem}\label{FixRel}
Let $E$ be a \fd\  \ca,
let $S$ be a complete system of matrix units for $E$,
let $N \in \N$, and let $\ep > 0$.
Then there exists $\dt > 0$ such that whenever $A$ is a unital \ca,
$p \in A$ is a \pj, $\ph \colon E \to p A p$ is a unital \hm, and
$e_0, e_1, \dots, e_N \in A$ are \pj s with
$\sum_{j = 0}^N e_j = 1$ and such that
$\| e_j \ph (s) - \ph (s) e_j \| \leq \dt$
for all $j$ and all $s \in S$,
then there is a unitary $u \in A$ with $\| u - 1 \| < \ep$ such that
$(u e_j u^*) \ph (a) = \ph (a) (u e_j u^*)$
for all $j$ and all $a \in E$.
\end{lem}

\begin{proof}
Replacing $E$ by $E \oplus \C$
and $\ph$ by $(a, \ld) \mapsto \ph (a) + \ld (1 - p)$,
we see that it suffices to prove the result under the additional
restriction $p = 1$ (with a different value of $\ep$).

Now apply Lemma~2.5.10 of~\cite{LnBook} to all dimensions up to
$\dim (E)$,
finding $\dt > 0$ such that whenever $B$ and $C$ are subalgebras
of a unital \ca\  $A$ such that $\dim (B) \leq \dim (E)$
and such that $B$ has a complete system of matrix units each of
which has distance less than $N (N + 1) \dt$ from $C$, then there is
a unitary $u \in A$ with $\| u - 1 \| < \ep$
such that $u^* B u \S C$.
Apply this with $B = \ph (E)$ and
$C = \bigoplus_{j = 0}^N e_j A e_j$.
We note that if $s \in S$ then,
because $\| e_j \ph (s) - \ph (s) e_j \| \leq \dt$ for all $j$,
\[
{\mathrm{dist}} (\ph (s), \, C)
  \leq \left\| s - \ssum{j = 0}{N} e_j \ph (s) e_j \right\|
  \leq \sum_{j \neq k} \| e_j \ph (s) e_k \|
  < N (N + 1) \dt,
\]
and that every element of $C$ commutes with every $e_j$.
\end{proof}

\begin{lem}\label{StARPDfn}
Let $A$ be a stably finite \suca\  %
and let $\af \in \Aut (A)$ satisfy $\af^n = \id_A$.
Then the action of $\Zqn$ generated by $\af$ has the
tracial Rokhlin property \ifo\  the conditions of
Definition~\ref{ARPDfn} are satisfied, except that~(1) is replaced by:
\begin{itemize}
\item[(1$'$)]
$e = \sum_{j = 0}^{n - 1} e_j$ is $\af$-invariant, and
$\| \af (e_j) - e_{j + 1} \| < \ep$ for $0 \leq j \leq n - 1$,
where, following Convention~\ref{ARPConv}, we take $e_n = e_0$.
\end{itemize}
The analogous statement holds for the strict Rokhlin property of
Definition~\ref{ERPDfn}.
\end{lem}

\begin{proof}
That~(1$'$) implies~(1) is trivial.
For the reverse direction, we first prove~(1$'$)
without the condition on $\af$-invariance of $e$.
Apply Definition~\ref{ARPDfn} with
$\frac{1}{n} \ep$ in place of $\ep$.
Then, using $\af^n = \id_A$,
\[
\| e_0 - \af (e_{n - 1} ) \|
  = \| \af^{n - 1} (e_0) - e_{n - 1} \|
  \leq \sum_{j = 0}^{n - 2}
        \| \af^{n - j - 2} ( \af (e_j) - e_{j + 1} ) \|
  < (n - 1) \ts{ \frac{1}{n} } \ep < \ep.
\]

Now we arrange for $\af$-invariance of $\sum_{j = 0}^{n - 1} e_j$.
\Wolog\  $\| a \| \leq 1$ for all $a \in F$.

Let
\[
\dt_1 = \min \left( \frac{\ep}{5}, \, \frac{1}{4}, \,
              \frac{1}{4 n + 1} \right).
\]
Apply Lemma~\ref{FixRel} with $E = \C \oplus \C$,
with $N = n$, and with $\dt_1$ in place of $\ep$.
Let $\dt_2$ be the resulting value of $\dt$.
Then set
\[
\dt_3 = \min \left( \dt_1, \, \frac{\dt_2}{4 n^2}, \,
              \frac{1}{4 n^2} \right).
\]

Apply the version of~(1$'$) without $\af$-invariance of $e$,
and with $\dt_3$ in place of $\ep$.
Let $f_0, f_1, \dots, f_{n - 1}$ be the resulting \pj s.
Set $f = \sum_{j = 0}^{n - 1} f_j$.
Then
\[
\| \af (f) - f \| \leq \sum_{j = 0}^{n - 1} \| \af (f_j) - f_{j + 1} \|
  < n \dt_3.
\]
Inductively, we get $\| \af^k (f) - f \| \leq k n \dt_3$ for $k \in \N$.
In particular,
with $c = \ts{\frac{1}{n}} \sum_{j = 0}^{n - 1} \af^j (f)$,
we get
\[
\| c - f \|
 \leq \frac{1}{n} \sum_{j = 0}^{n - 1} \| \af^j (f) - f \|
 < \frac{1}{n} \sum_{j = 0}^{n - 1} j n \dt_3
 \leq n^2 \dt_3.
\]
Since $n^2 \dt_3 < \frac{1}{2}$, the \pj\  %
$e = \ch_{ (\frac{1}{2}, \I)} (c)$ is defined,
is $\af$-invariant (because $c$ is), and satisfies
$\| e - c \| \leq \| c - f \|$,
so that $\| e - f \| \leq 2 \| c - f \| < 2 n^2 \dt_3$.
It follows that
$\| f_j e - e f_j \| < 4 n^2 \dt_3 \leq \dt_2$
for $0 \leq j \leq n - 1$.
By the choice of $\dt_2$ and following Lemma~\ref{FixRel},
with the \hm\  $\ph \colon \C \oplus \C \to A$
being $\ph (\ld_0, \ld_1) = \ld_0 e + \ld_1 (1 - e)$
and with $f_0, \, f_1, \, \dots, \, f_{n - 1}, \, 1 - f$ in place of
$e_0, e_1, \dots, e_N$,
there is a unitary $u \in A$ with $\| u - 1 \| < \dt_1$ such that
$e_j = u f_j u^*$ commutes with $e$ for $0 \leq j \leq n - 1$.
We get $\| e_j - f_j \| < 2 \dt_1$.

Since $\sum_{j = 0}^{n - 1} f_j = f$ we get
\[
\left\| e - \ssum{j = 0}{n - 1} e_j \right\|
 \leq \| e - f \| + \sum_{j = 0}^{n - 1} \| f_j - e_j \|
 < 2 n^2 \dt_3 + 2 n \dt_1 < 1.
\]
Since $e$ commutes with $\sum_{j = 0}^{n - 1} e_j$, this
implies that $e = \sum_{j = 0}^{n - 1} e_j$.
Since $e$ is unitarily equivalent to $f$, we now have
Conditions~(3) and~(4) of Definition~\ref{ARPDfn}.

{}From $\| e_j - f_j \| < 2 \dt_1$, we get
\[
\| e_j a - a e_j \|
  \leq \| f_j a - a f_j \| + 2 \| e_j - f_j \|
  < \dt_3 + 4 \dt_1 \leq 5 \dt_1 \leq \ep.
\]
This is Condition~(2) of Definition~\ref{ARPDfn}.
Moreover, for $0 \leq j \leq n - 1$ we get
\begin{align*}
\| \af (e_j) - e_{j + 1} \|
&  \leq \| \af (f_j) - f_{j + 1} \|
      + \| e_j - f_j \| + \| e_{j + 1} - f_{j + 1} \|  \\
&  < \dt_3 + 2 \dt_1 + 2 \dt_1 \leq 5 \dt_1 \leq \ep.
\end{align*}
This is the estimate in Condition~(1$'$) of the lemma,
and finishes the proof of the part about the \aRp.

The proof for the strict Rokhlin property is the same as the
first part of the proof for the \aRp,
since in this case $e = 1$ is automatically $\af$-invariant.
\end{proof}

We finish this section by proving that crossed products by
actions with the \aRp\  are still simple.

\begin{lem}\label{ARPImpOuter}  
Let $A$ be a stably finite \suca,
and let $\af \in \Aut (A)$ be an automorphism which
satisfies $\af^n = \id_A$ and such that the action of $\Zqn$
generated by $\af$ has the tracial Rokhlin property.
Then $\af^k$ is outer for $1 \leq k \leq n - 1$.
\end{lem}

\begin{proof}
Let $1 \leq k \leq n - 1$ and let $u \in A$ be unitary.
We prove that $\af^k \neq \Ad (u)$.
Apply Definition~\ref{ARPDfn} with $F = \{ u \}$,
with $\ep = \frac{1}{2 n}$, with $N = 1$, and with $x = 1$.
Then for $1 \leq k \leq n - 1$, we have
\[
\| \af^{k} (e_0) - e_{k} \|
  \leq \sum_{j = 0}^{k - 1}
        \| \af^{k - j - 1} ( \af (e_j) - e_{j + 1} ) \|
  < k \ep < \ts{\frac{1}{2}}.
\]
Also $\| e_0 u - u e_0 \| < \ep < \frac{1}{2}$,
whence $\| u e_0 u^* - e_0 \| < \frac{1}{2}$.
The choice $N = 1$ implies that $e_0 \neq 0$, so
orthogonality of $e_0$ and $e_k$ implies $\| e_k - e_0 \| = 1$.
It follows that
\[
\| \af^k (e_0) - u e_0 u^* \|
  \geq \| e_k - e_0 \| - \| \af^k (e_0) - e_k \| - \| u e_0 u^* - e_0 \|
  > 0.
\]
Therefore $\af^k \neq \Ad (u)$.
\end{proof}

\begin{cor}\label{CrPrIsSimple}
Let $A$ be a stably finite \suca,
and let $\af \in \Aut (A)$ be an automorphism which
satisfies $\af^n = \id_A$ and such that the action of $\Zqn$
generated by $\af$ has the tracial Rokhlin property.
Then $\CZnAa$ is simple.
\end{cor}

\begin{proof}
This follows from Lemma~\ref{ARPImpOuter} and Theorem~3.1 of~\cite{Ks1}.
\end{proof}

\section{Crossed products by actions on C*-algebras with
   tracial rank zero}\label{Sec:CrPrd}

\indent
In this section, we prove that the crossed product of a
\ca\  with tracial rank zero by an action with the \aRp\  %
again has tracial rank zero.

The following result gives the criterion we use for a
simple \ca\  to have tracial rank zero.
Note that, by Theorem~7.1(a) of~\cite{LnTTR},
tracial rank zero is the same as
tracially AF in the sense of Definition~2.1 of~\cite{LnTAF}.

\begin{prp}\label{TAFCond}
Let $A$ be a simple separable unital \ca.
Then $A$ has tracial rank zero in the sense of
Definition~3.1 of~\cite{LnTTR} \ifo\  the following holds.

For every finite set $S \S A$, every $\ep > 0$,
every nonzero positive element $x \in A$, and every $N \in \N$,
there is a \pj\  $p \in A$ and a finite
dimensional unital subalgebra $E \S p A p$ (that is,
$p$ is the identity of $E$) such that:
\bei
\item[(1)]
$\| p a - a p \| < \ep$ for all $a \in S$.
\item[(2)]
For every $a \in S$ there exists $b \in E$ such that
$\| p a p - b \| < \ep$.
\item[(3)]
$1 - p$ is \mvnt\  to a \pj\  in ${\overline{x A x}}$.
\item[(4)]
There are $N$ \mops\  in $p A p$, each of which is \mvnt\  to $1 - p$.
\eei
\end{prp}

\begin{proof}
It is immediate that tracial rank zero implies
the condition in the proposition.

The proof of the other direction
is similar to the proof of Proposition~3.8 of~\cite{LnTAF},
and we refer to parts of that proof for some of the argument.
However, there appears to be a gap there, so we give more detail
at the relevant point in our argument.

According to Definition~2.1 of~\cite{LnTAF}, we must prove that
if $S$, $\ep$, $x$, and $N$ are as in the hypotheses, and if in
addition $a_0 \in A$ is a given nonzero element, then $p$ and $E$
can in addition be chosen so that $\| p a_0 p \| > \| a_0 \| - \ep$.
\Wolog\  $\| a_0 \| = 1$ and $\ep < 1$.

As in the first paragraph
of the proof of Proposition~3.8 of~\cite{LnTAF},
we may assume that $A \not\cong M_n$ for any $n$,
and we conclude that $A$ is stably finite and has Property~(SP).

Now let $S$, $\ep$, $x$, and $N$ be as in the hypotheses, and
let $a_0 \in A$ satisfy $\| a_0 \| = 1$.
Define a \cfn\  $f \colon [0, \I) \to [0, 1]$ by
\[
f (t) = \left\{ \begin{array}{ll}
     0     & \hspace{3em}  0 \leq t \leq 1 - \ts{ \frac{1}{4} } \ep \\
     8 \ep^{-1} \left( t - 1 + \ts{ \frac{1}{4} } \ep \right)
           & \hspace{3em}  1 - \ts{ \frac{1}{4} } \ep \leq t
                                 \leq 1 - \ts{ \frac{1}{4} } \ep  \\
     1     & \hspace{3em}  1 - \ts{ \frac{1}{8} } \ep \leq t
    \end{array} \right..
\]
Then $b = f ( a_0^* a_0 ) \neq 0$,
so there is a nonzero \pj\  $q_0 \in {\overline{b A b}}$,
and by Lemma~3.1 of~\cite{LnTAF} there is a nonzero \pj\  $q \leq q_0$
such that $q$ is \mvnt\  to a \pj\  in ${\overline{x A x}}$.
Use Lemma~3.2 of~\cite{LnTAF} to choose nonzero
\mops\  $q_1, \, q_2 \leq q$ such that $q_1 \sim q_2$.

Choose $\ep_0 > 0$ so small that whenever $D$ is a unital \ca\  and
$e, \, f \in D$ are \pj s such that $\| e f \| < \ep_0$, then
$f \precsim 1 - e$.
Choose $\ep_1 > 0$ so small that whenever $D$ is a unital \ca\  and
$b \in D_{\sa}$ satisfies $\| b^2 - b \| < \ep_1$, then
there is a \pj\  $e \in D$ with
$\| e - b \|
  < \min \left( \frac{1}{12} \ep, \, \frac{1}{2} \ep_0^2 \right)$.
Clearly $\ep_1 < 1$.
Set
\[
\dt = \min \left( \ts{ \frac{1}{12} } \ep, \, \ts{ \frac{1}{2} } \ep_0^2,
          \, \ts{ \frac{1}{5} } \ep_1 \right).
\]
Apply the hypotheses
with $S \cup \{ q, \, a_0 \}$ in place of $S$,
with $\dt$ in place of $\ep$,
with $q_1$ in place of $x$,
and with $\max (N, 2)$ in place of $N$.
Let $p$ and $E \S p A p$ be the resulting \pj\  and \fd\  subalgebra.
We finish the proof by showing that $\| p a_0 p \| > \| a_0 \| - \ep$.

If $\| p q \| < \ep_0$ then $q \precsim 1 - p$, whence
\[
1 - p \precsim q_1 \leq 1 - p
\andeqn 1 - p \precsim q_2 \leq 1 - p.
\]
Since $q_1 q_2 = 0$, this contradicts stable finiteness of $A$.
So $\| p q \| \geq \ep_0$, whence $\| p q p \| \geq \ep_0^2$.
By the choice of $p$ and $E$,
there is $b \in E$ with
\[
\| b - p q p \|
 < \min \left( \ts{ \frac{1}{12} } \ep, \, \ts{ \frac{1}{2} } \ep_0^2,
                 \, \ts{ \frac{1}{5} } \ep_1 \right).
\]
We have
\[
\| (p q p)^2 - p q p \| = \| p q (p q - q p ) p \|
   <  \ts{ \frac{1}{5} } \ep_1,
\]
so
\begin{align*}
\| b^2 - b \|
 & \leq \| (p q p)^2 - p q p \| + \| b \| \cdot \| b - p q p \|
      + \| b - p q p \| \cdot \| p q p \| + \| b - p q p \|   \\
 & < \left( 4 + \ts{ \frac{1}{5} } \ep_1 \right) \ts{ \frac{1}{5} } \ep_1
   \leq \ep_1.
\end{align*}
Therefore there is a \pj\   $e \in E$ with
$\| e - b \|
  < \min \left( \frac{1}{2} \ep_0^2, \, \frac{1}{12} \ep \right)$.
This gives
$\| e - p q p \| < \frac{1}{2} \ep_0^2 + \frac{1}{2} \ep_0^2 = \ep_0^2$.
Since $\| p q p \| \geq \ep_0^2$, we have $e \neq 0$.
Furthermore, $\| e - p q p \| < \frac{1}{6} \ep$.

Similarly
$\| (q p q)^2 - q p q \| < \ts{ \frac{1}{5} } \ep_1 < \ep_1$,
so there is a
\pj\  $e_0 \in q A q$ such that $\| e_0 - q p q \| < \frac{1}{12} \ep$.
Since $\| p q p - q p q \| \leq 2 \| p q - q p \| < \frac{1}{6} \ep$,
we get
$\| e - e_0 \| < \frac{1}{6} \ep + \frac{1}{6} \ep + \frac{1}{12} \ep
  < \frac{1}{2} \ep$.
Because $\ep < 1$ this implies $e \sim e_0$, so $e_0 \neq 0$.
Also, since $\| (a_0^* a_0)^{1/2} \| = 1$ we immediately get
\[
\| (a_0^* a_0)^{1/2} e (a_0^* a_0)^{1/2}
   - (a_0^* a_0)^{1/2} e_0 (a_0^* a_0)^{1/2} \| < \ts{ \frac{1}{2} } \ep.
\]

Since $e_0 \neq 0$ and $e_0 \leq q \in {\overline{b A b}}$,
the definition of $b$ implies
$\| e_0 (a_0^* a_0) e_0 \| \geq 1 - \frac{1}{4} \ep$.
Therefore
\begin{align*}
1 & \leq \| [(a_0^* a_0)^{1/2} e_0]^* [(a_0^* a_0)^{1/2} e_0] \|
           + \ts{ \frac{1}{4} } \ep
    = \| (a_0^* a_0)^{1/2} e_0 (a_0^* a_0)^{1/2} \| + \ts{ \frac{1}{4} } \ep
                \\
  & < \| (a_0^* a_0)^{1/2} e (a_0^* a_0)^{1/2} \|
         + \ts{ \frac{1}{2} } \ep + \ts{ \frac{1}{4} } \ep
    \leq \| (a_0^* a_0)^{1/2} p (a_0^* a_0)^{1/2} \|
         + \ts{ \frac{1}{2} } \ep + \ts{ \frac{1}{4} } \ep
                \\
  & = \| p (a_0^* a_0) p \| + \ts{ \frac{3}{4} } \ep.
\end{align*}
So $\| p (a_0^* a_0) p \| > 1 - \ts{ \frac{3}{4} } \ep$.

Using $\| p a_0 - a_0 p \| < \frac{1}{12} \ep$, we now get
\[
\| p a_0 p \|^2 = \| (p a_0^* p) (p a_0 p) \|
  \geq \| p a_0^* a_0 p \| - \ts{ \frac{1}{12} } \ep
  > 1 - \ts{ \frac{3}{4} } \ep - \ts{ \frac{1}{12} } \ep
  > 1 - \ep.
\]
So $\| p a_0 p \| > \sqrt{1 - \ep} \geq 1 - \ep$,
as desired.
\end{proof}

We also recall the properties of simple unital \ca s
with tracial rank zero.

\begin{thm}\label{TAFProp}
(H.\  Lin.)
Let $A$ be a simple unital \ca\  with tracial rank zero.
Then $A$ has real rank zero, stable rank one, and cancellation of \pj s
(Definition~\ref{CancD}).
Moreover, the order on \pj s over $A$ is determined by traces
(Definition~\ref{OrdDetD}).
\end{thm}

\begin{proof}
In view of Theorem~7.1(a) of~\cite{LnTTR},
real rank zero and stable rank one are Theorem~3.4 of \cite{LnTAF},
and cancellation of \pj s is Lemma~3.3 of \cite{LnTAF}.
That the order is determined by traces is
Theorems~6.8 and~6.13 of~\cite{LnTTR}.
\end{proof}

\begin{lem}\label{PreserveFD}
Let $A$ be a \suca, and let $\af \in \Aut (A)$ be an automorphism which
satisfies $\af^n = \id_A$ and such that the action of $\Zqn$
generated by $\af$ has the tracial Rokhlin property.
Then for every finite set
$F \S A$, every \fd\  subalgebra $E \S A$, every $\ep > 0$,
every $N \in \N$,
and every nonzero positive element $x \in A$,
there are \mops\  $e_0, e_1, \dots, e_{n - 1} \in A$
and a unitary $v \in A$ such that:
\begin{itemize}
\item[(1)]
$\| \af (e_j) - e_{j + 1} \| < \ep$ for $0 \leq j \leq n - 1$,
where, following Convention~\ref{ARPConv},  we take $e_n = e_0$.
\item[(2)]
$\| e_j a - a e_j \| < \ep$ for $0 \leq j \leq n - 1$ and all $a \in F$.
\item[(3)]
With $e = \sum_{j = 0}^{n - 1} e_j$, the \pj\  $1 - e$
is \mvnt\  to a \pj\  in ${\overline{x A x}}$.
\item[(4)]
For every $j$ with $0 \leq j \leq n - 1$,
there are $N$ \mops\  $f_1, f_2, \dots, f_N \leq e_j$,
each of which is \mvnt\  to the \pj\  $1 - e$ of~(3).
\item[(5)]
$\| v - 1 \| < \ep$, and $e_j$ commutes with all elements of
$v E v^*$ for $0 \leq j \leq n - 1$.
\item[(6)]
$\af (e) = e$.
\end{itemize}
\end{lem}

\begin{proof}
Let $F$, $E$, $\ep$, $N$, and $x$ be as in the hypotheses.
\Wolog\  $\| a \| \leq 1$ for all $a \in F$.
Let $S$ be a complete system of matrix units for $E$.
Apply Lemma~\ref{FixRel} with this $E$, $S$, and $\ep$, and with
$n$ in place of $N$, obtaining $\dt_0 > 0$.
Define
$\dt = \min \left( \ep, \dt_0 \right)$.

Apply Lemma~\ref{StARPDfn} with $F \cup S$ in place of $F$,
with $\dt$ in place of $\ep$,
and with $N$ and $x$ as given.
Let $e_0, \, e_1, \, \dots, \, e_{n - 1}$ and
$e = \sum_{j = 0}^{n - 1} e_j$ be the resulting \pj s.
Note that $\af (e) = e$.
Parts~(3), (4), and~(6) of the conclusion are immediate.
Part~(1) follows from
$\| \af (e_j) - e_{j + 1} \| < \dt \leq \ep$
for $0 \leq j \leq n - 1$.
It also follows that for all $a \in F \cup S$ we have
$\| e_j a - a e_j \| < \dt \leq \ep$, which gives Part~(2).

Since $\| e_j a - a e_j \| < \dt \leq \dt_0$ for $a \in S$,
the choice of $\dt_0$ from Lemma~\ref{FixRel}, taking
the \hm\  $\ph$ there to be the inclusion of $E$ in $A$,
provides a unitary $v_0 \in A$ with $\| v_0 - 1 \| < \ep$ and such that
$v_0 e_j v_0^*$ commutes with every element of $E$.
We obtain Part~(5) of the conclusion by taking $v = v_0^*$.
\end{proof}

\begin{lem}\label{PertOfCrPrd}
Let $A$ be a unital \ca\  and let $\af \in \Aut (A)$
satisfy $\af^n = \id_A$.
Let $w \in A$ be a unitary such that
\[
w \af (w) \af^2 (w) \cdots \af^{n - 1} (w) = 1.
\]
Then the automorphism $\bt = \Ad (w) \circ \af \in \Aut (A)$
satisfies $\bt^n = \id_A$.
Moreover, letting $u \in \CZnAa$
and $v \in \Cs{n}{A}{\bt}$ be the canonical
unitaries implementing the automorphisms $\af$ and $\bt$,
there is an isomorphism
\[
\ph \colon
 \Cs{n}{A}{\bt} \to \CZnAa
\]
such that $\ph (a) = a$ for all $a \in A$ and
$\| \ph (v) - u \| = \| w - 1 \|$.
\end{lem}

\begin{proof}
As is implicit in the statement of the lemma, we identify
$A$ with its image in each of the crossed products.

That $\bt^n = \id_A$ is easy to check.

The unitary $w u$ is in $\CZnAa$.
For $a \in A$, we have
$(w u) a (w u)^* = w \af (a) w^* = \bt (a)$, by the definition of $\bt$.
Moreover, using $u^n = 1$ and $u a u^* = \af (a)$ for $a \in A$, we get
\begin{align*}
(w u)^n & = w (u w u^{-1}) (u^2 w u^{-2})
            \cdots (u^{n - 1} w u^{-(n - 1)}) u^n  \\
    & = [w \af (w) \af^2 (w) \cdots \af^{n - 1} (w)] u^n = 1.
\end{align*}
The universal property of crossed products therefore provides a \hm\  %
\[
\ph \colon
 \Cs{n}{A}{\bt} \to \CZnAa
\]
such that $\ph (a) = a$ for all $a \in A$ and $\ph (v) = w u$.

A similar argument shows that there is a \hm\  %
\[
\ps \colon
 \CZnAa \to \Cs{n}{A}{\bt}
\]
such that $\ps (a) = a$ for all $a \in A$ and $\ph (u) = w^* v$.
One checks that $\ps \circ \ph (b) = b$ for all
$b \in \Cs{n}{A}{\bt}$ by checking this for $b \in A$
and for $b = v$.
Similarly $\ph \circ \ps = \id_{\CZnAa}$.
Therefore $\ph$ is an isomorphism.

Finally, $\| \ph (v) - u \| = \| w u - u \| = \| w - 1 \|$, as desired.
\end{proof}

\begin{lem}\label{SPAndInv}
Let $A$ be a simple \ca\  with Property~(SP), and
let $\af \in \Aut (A)$ satisfy $\af^n = \id_A$.
Let $p \in A$ be a nonzero \pj.
Then there exists a nonzero \pj\  $q \in A$ such that
$\af^j (q) \precsim p$ for all $j$.
\end{lem}

\begin{proof}
Using Property~(SP) and Lemma~3.1 of~\cite{LnTAF},
find a nonzero \pj\  $e_1 \leq p$ such that $e_1 \precsim \af (p)$.
In the same way, find a
nonzero \pj\  $e_2 \leq e_1$ such that $e_2 \precsim \af^2 (p)$.
Proceed inductively.
Set $q = e_{n - 1}$.
Then $q$ is a nonzero \pj\  such that
$q \precsim \af^j (p)$ for $0 \leq j \leq n - 1$, hence for all $j$.
\end{proof}

\begin{lem}\label{P}
Let $A$ be a unital \ca,
let $\af \in \Aut (A)$ satisfy $\af^n = \id_A$, and let $\ep > 0$.
Let $e_0, \, e_1, \, \dots, \, e_{n - 1} \in A$ be \mops,
let $e = \sum_{j = 1}^{n - 1} e_j$, and assume that $\af (e) = e$.
Let $w_1, \, w_2, \, \dots, \, w_{n - 1} \in A$
be partial isometries, satisfying
\[
w_j w_j^* = e_j, \,\,\,\,\,\, w_j^* w_j = \af ( e_{j - 1} ), \andeqn
\| w_j - e_j \| < \ep
\]
for $1 \leq j \leq n - 1$.
Then the element
\[
w = 1 - e + w_1 + w_2 + \cdots w_{n - 1}
  + \af^{n - 1} (w_1^*) \af^{n - 2} (w_2^*)  \cdots \af (w_{n - 1}^*)
\]
is a unitary in $A$ satisfying:
\begin{itemize}
\item[(1)]
$\| w - 1 \| < 2 n^2 \ep$.
\item[(2)]
$w \af (w) \af^2 (w) \cdots \af^{n - 1} (w) = 1$.
\item[(3)]
Following Convention~\ref{ARPConv},
$w \af ( e_{j - 1} ) w^* = e_j$ for $1 \leq j \leq n$.
\item[(4)]
$\| (\Ad (w) \circ \af)^k (a) - \af^k (a) \| \leq 4 k n^2 \ep \| a \|$
for $a \in A$ and $k \in \N$.
\item[(5)]
$(\Ad (w) \circ \af)^n = \id_A$.
\end{itemize}
\end{lem}

\begin{proof}
Define
\[
z = \af^{n - 1} (w_1^*) \af^{n - 2} (w_2^*)  \cdots \af (w_{n - 1}^*).
\]
We claim that
\[
z z^* = e_0, \,\,\,\,\,\, z^* z = \af ( e_{n - 1} ), \andeqn
\| z - e_0 \| < n^2 \ep.
\]
For the first, we observe that
\[
\af (w_{n - 1}^*) \af (w_{n - 1}) = \af^2 (e_{n - 2}),
\]
\[
\af^2 (w_{n - 2}^*) \af (w_{n - 1}^*) \af (w_{n - 1}) \af^2 (w_{n - 2})
  = \af^2 (w_{n - 2}^* e_{n - 2} w_{n - 2}) = \af^3 (e_{n - 3}),
\]
etc., ending with
\[
z z^* = \af^n (e_0) = e_0.
\]
For the second, we observe that
\[
\af^{n - 1} (w_1) \af^{n - 1} (w_1^*) = \af^{n - 1} (e_1),
\]
\[
\af^{n - 2} (w_2) \af^{n - 1} (w_1)
              \af^{n - 1} (w_1^*) \af^{n - 2} (w_2^*)
    = \af^{n - 2} (w_2 \af (e_1) w_2^*) = \af^{n - 2} (e_2),
\]
etc., ending with
\[
z^* z = \af (e_{n - 1}).
\]
For the third, first observe that
\[
\| \af^{n - k} (e_k) - e_0 \| = \| \af^{- k} (e_k) - e_0 \|
   = \| e_k - \af^{k} (e_0) \| < k \ep,
\]
whence
\[
\| \af^{n - k} (w_k) - e_0 \|
  \leq \| \af^{n - k} (w_k - e_k) \| + \| \af^{n - k} (e_k) - e_0 \|
  < (k + 1) \ep.
\]
Therefore
\[
\| z - e_0 \| < \sum_{k = 1}^{n - 1} (k + 1) \ep < n^2 \ep.
\]
This completes the proof of the claim.

The element $w$ of the statement is now defined by
\[
w = 1 - e + w_1 + w_2 + \cdots w_{n - 1} + z.
\]
Since $\af (e) = e$, we have
$\sum_{j = 0}^{n - 1} \af (e_{j - 1}) = e = \sum_{j = 0}^{n - 1} e_j$,
so that $w$ is in fact unitary.
Clearly $w \af (e_{j - 1}) w^* = e_j$.
This is Part~(3) of the conclusion.

Part~(1) of the conclusion is the estimate
\[
\| w - 1 \| \leq \| z - e_0 \| + \sum_{j = 1}^{n - 1} \| w_j - e_j \|
     < n^2 \ep + (n - 1) \ep < 2 n^2 \ep.
\]

To prove Part~(2), we simplify the notation by setting $w_0 = z$
and interpreting all subscripts as elements of $\Zqn$.
{}From $w_j w_j^* = e_j$ and  $w_j^* w_j = \af (e_{j - 1})$
we get $\af^k (w_j) \af^k (w_j)^* = \af^k (e_j)$ and
$\af^k (w_j)^* \af^k (w_j) = \af^{k + 1} (e_{j - 1})$.
It follows that
\[
w \af (w) \af^2 (w) \cdots \af^{n - 1} (w)
 = 1 - e + \sum_{j = 0}^{n - 1}
   w_j \af (w_{j - 1}) \af^2 (w_{j - 2})
              \cdots \af^{n - 1} (w_{j - n + 1}).
\]
Because $\af^n = \id_A$, we have
\begin{align*}
w_0
 & = \af^{n - 1} (w_1^*) \af^{n - 2} (w_2^*)  \cdots \af (w_{n - 1}^*)
           \\
 & = [\af^{- 1} (w_1^*) \cdots \af^{- j} (w_j^*)]
       [\af^{n - j - 1} (w_{j - n + 1}^*) \cdots  \af (w_{- 1}^*)],
\end{align*}
and we get, by substituting for the term $\af^j (w_0)$,
\[
w_j \af (w_{j - 1}) \af^2 (w_{j - 2})
              \cdots \af^{n - 1} (w_{j - n + 1})
   = e_j \af^n (e_{j - n}) = e_j^2 = e_j
\]
for all $j$.
This gives Part~(2) of the conclusion.

Define $\bt = \Ad (w) \circ \af \in \Aut (A)$.
Then $\bt^n = \id_A$ by Part~(2) and Lemma~\ref{PertOfCrPrd}.
This is Part~(5).
Also, for $a \in A$ we have
\[
\| \bt (a) - \af (a) \| \leq 2 \| w - 1 \| \cdot \| a \|
   \leq 4 n^2 \ep \| a \|,
\]
and Part~(4) of the conclusion follows by induction.
\end{proof}

\begin{thm}\label{RokhTAF}
Let $A$ be a \suca, and let $\af \in \Aut (A)$ be an automorphism which
satisfies $\af^n = \id_A$ and such that the action of $\Zqn$
generated by $\af$ has the tracial Rokhlin property.
Suppose that $A$ has tracial rank zero.
Then $\CZnAa$ has tracial rank zero.
\end{thm}

\begin{proof}
It suffices to verify the condition of Proposition~\ref{TAFCond},
for a finite set $S$ of the form $S = F \cup \{ u \}$, where
$F$ is a finite subset of the unit ball of $A$ and
$u \in \CZnAa$
is the canonical unitary implementing the automorphism $\af$.
So let $F \S A$ be a finite subset
with $\| a \| \leq 1$ for all $a \in F$, let $\ep > 0$,
let $N \in \N$, and let
$x \in \CZnAa$ be a nonzero positive element.

The \ca\  $A$ has Property~(SP) by Theorem~\ref{TAFProp}.
So we can apply Theorem~4.2 of \cite{JO}, with $N = \{ 1 \}$,
to find a nonzero \pj\  $p_0 \in A$ which is \mvnt\  in
$\CZnAa$ to a \pj\  in
${\overline{x \CZnAa x}}$.

Since $A$ has real rank zero (by Theorem~\ref{TAFProp}), we may
use Theorem~1.1(i) of \cite{Zh7} to find $(N + 1) (n + 1)$ nonzero
\mvnt\  \mops\  in $p_0 A p_0$.
Call one of them $p_1$.
Use Lemma~\ref{SPAndInv} to find a nonzero \pj\  $p \leq p_1$
such that $\af^j (p) \precsim p_1$ for all $j$.

Set
\[
\ep_0 = \frac{\ep}{12 (n + 1)^5}.
\]
Choose $\dt > 0$ with $\dt < \ep_0$,
and so small that whenever $e$ and $f$ are \pj s in a \ca\  $C$
such that $\| e - f \| < \dt$, then there is a partial isometry
$s \in C$ such that
\[
s s^* = e, \,\,\,\,\,\, s^* s = f, \andeqn
\| s - e \| < \ep_0.
\]

Apply the condition for tracial rank zero in Proposition~\ref{TAFCond}
with $\dt$ in place of $\ep$, with the finite set $S$ there being
\[
S_0 = F \cup \af (F) \cup \cdots \cup \af^{n - 1} (F),
\]
and with $p$ in place of $a$.
We obtain a \pj\  $q_0$ such that
$1 - q_0 \precsim p$,
and a \fd\  subalgebra $E_0$ with $q_0 \in E_0 \S q_0 A q_0$,
such that for every $a \in S_0$ we have $\| q_0 a - a q_0 \| < \dt$
and there exists $b \in E_0$ such that $\| q_0 a q_0 - b \| < \dt$.

Apply Lemma~\ref{PreserveFD} with $S_0$ in place of $F$,
with $E_0 + \C (1 - q_0)$ in place of the \fd\  subalgebra $E$,
with $\dt$ in place of $\ep$, with $N = 1$,
and with $p$ in place of $x$.
We obtain a unitary
$y \in A$ and \mops\  $e_0, \, e_1, \, \dots, \, e_{n - 1} \in A$ which
commute with all elements of $y (E_0 + \C q_0 ) y^*$,
such that $\| e_j a - a e_j \| < \dt$ for all $a \in S_0$,
such that $\| y - 1 \| < \dt$,
such that $\| \af (e_j) - e_{j + 1} \| < \dt$,
and such that $e =\sum_{j = 0}^{n - 1} e_j$ is $\af$-invariant
and $1 - e \precsim p$.

According to the choice of $\dt$, for $1 \leq j \leq n - 1$
there are partial isometries $w_j \in A$ such that
\[
w_j w_j^* = e_j, \,\,\,\,\,\, w_j^* w_j = \af ( e_{j - 1} ), \andeqn
\| w_j - e_j \| < \ep_0.
\]
Apply Lemma~\ref{P} to the $e_j$ and $w_j$,
with $\ep_0$ in place of $\ep$.
We obtain a unitary $w$ as there such that $\| w - 1 \| < 2 n^2 \ep_0$,
and such that the automorphism $\bt = \Ad (w) \circ \af$
satisfies $\bt^n = \id_A$,
\[
\| \bt^k (a) - \af^k (a) \| \leq 4 k n^2 \ep_0 \| a \|
\]
for $k \in \N$ and $a \in A$, and $\bt (e_j) = e_{j + 1}$ for all $j$.

Define
\[
q = \sum_{k = 0}^{n - 1} \bt^k (e_0 y q_0 y^* e_0)
\andeqn
E = \bigoplus_{k = 0}^{n - 1} \bt^k (e_0 y E_0 y^* e_0).
\]
By construction, $e_0$ commutes with $y q_0 y^*$, and
the \pj s
\[
e_0, \, \bt (e_0), \, \dots, \, \bt^{n - 1} (e_0), \, 1 - e
\]
are orthogonal, so that $q$ is a $\bt$-invariant \pj.
Similarly, $e_0$ commutes with every element of $y E_0 y^*$,
so $E$ is a $\bt$-invariant \fd\  subalgebra of $A$.

Let $a \in F$.
We estimate $\| q a - a q \|$
and the distance from $q a q$ to $E$.
We begin by estimating $\| [ \bt^k (e_0 y q_0 y^* e_0), \, a] \|$.
Recall that $a \in F$ implies $\| a \| \leq 1$.
Using
\[
[ e_0 q_0 e_0, \, \af^{n - k} (a) ]
  = e_0 q_0 [e_0, \, \af^{n - k} (a) ]
     + e_0 [q_0, \, \af^{n - k} (a) ] e_0
     + [e_0, \, \af^{n - k} (a) ] q_0 e_0,
\]
and because $\af^{n - k} (a) \in S_0$, we get
\[
\| [ e_0 q_0 e_0, \, \af^{n - k} (a) ] \|
  \leq \| [q_0, \, \af^{n - k} (a) ] \|
         + 2 \| [e_0, \, \af^{n - k} (a) ] \|
    < \dt + 2 \dt = 3 \dt.
\]
So
\begin{align*}
\| [ \bt^k (e_0 q_0 e_0), \, a] \|
 & = \| [ e_0 q_0 e_0, \, \bt^{n - k} (a) ] \|  \\
 & \leq 2 \| \bt^{n - k} (a) - \af^{n - k} (a) \|
        + \| [ e_0 q_0 e_0, \, \af^{n - k} (a) ] \|  \\
 & < 8 (n - k) n^2 \ep_0 + 3 \dt
   < (8 n^3 + 3) \ep_0.
\end{align*}
Therefore
\begin{align*}
\| [ \bt^k (e_0 y q_0 y^* e_0), \, a] \|
 & = \| [ e_0 y q_0 y^* e_0, \, \bt^{n - k} (a) ] \|
   \leq 4 \| y - 1 \| + \| [ \bt^k (e_0 q_0 e_0), \, a] \|  \\
 & < 4 \dt + (8 n^3 + 3) \ep_0
   < (8 n^3 + 7) \ep_0.
\end{align*}
Now
\[
\| [ q, a ] \|
  \leq \sum_{k = 0}^{n - 1} \| [ \bt^k (e_0 y q_0 y^* e_0), \, a ] \|
  < n (8 n^3 + 7) \ep_0 < \ep.
\]

We next estimate the distance from $q a q$ to $E$.
We begin by estimating
\begin{align*}
& \left\| q a q
   - \ssum{k = 0}{n - 1}
    [\bt^k (e_0 y q_0 y^* e_0)] a [\bt^k (e_0 y q_0 y^* e_0)]
                    \right\|  \\
& \hspace*{3em} \mbox{}
  \leq \sum_{k = 0}^{n - 1} \sum_{l \neq k}
     \| [\bt^k (e_0 y q_0 y^* e_0)] a [\bt^l (e_0 y q_0 y^* e_0)] \|  \\
& \hspace*{3em} \mbox{}
  \leq \sum_{k = 0}^{n - 1} \sum_{l \neq k}
     \left( \| [\bt^k (e_0 y q_0 y^* e_0)] [\bt^l (e_0 y q_0 y^* e_0)] a \|
         + \| [ \bt^l (e_0 y q_0 y^* e_0), \, a ] \|
                          \rule{0em}{2.3ex} \right).
\end{align*}
In each summand in the last expression, the first term contains the
expression $\bt^k (e_0) \bt^l (e_0) = e_k e_l$, which is zero because
$k \neq l$.
Therefore this term is zero.
The second term was estimated above by $(8 n^3 + 7) \ep_0$,
and there are fewer than $n^2$ summands, so we conclude that
\[
\left\| q a q
   - \ssum{k = 0}{n - 1}
    [\bt^k (e_0 y q_0 y^* e_0)] a [\bt^k (e_0 y q_0 y^* e_0)]
                    \right\|
  < n^2 (8 n^3 + 7) \ep_0.
\]

By construction, there exists $b_k \in E_0$ such that
$\| b_k - q_0 \af^{n - k} (a) q_0 \| < \dt$.
Set $c_k = e_0 y b_k y^* e_0$.
Then set $c = \sum_{k = 0}^{n - 1} \bt^k (c_k) \in E$.
Since $e_0$ commutes with $y q_0 y^*$, we have,
recalling at the fifth step that
$\| \bt^k (a) - \af^k (a) \| \leq 4 k n^2 \ep_0 \| a \|$ and
$\| y - 1 \| < \dt$,
\begin{align*}
\| [e_0 y q_0 y^* e_0] \bt^{n - k} (a) [e_0 y q_0 y^* e_0]
   - c_k \|
& = \| [e_0 y q_0 y^*] \bt^{n - k} (a) [ y q_0 y^* e_0]
   - e_0 y b_k y^* e_0 \|   \\
& \leq \| q_0 y^* \bt^{n - k} (a) y q_0 - b_k \|  \\
& \leq 2 \| y - 1 \| + \| q_0 \bt^{n - k} (a) q_0 - b_k \|  \\
& < 2 \dt + 4 (n - k) n^2 \ep_0
    + \| q_0 \af^{n - k} (a) q_0 - b_k \|  \\
& < 2 \dt + 4 n^3 \ep_0 + \dt
  \leq (4 n^3 + 3) \ep_0.
\end{align*}
It follows that
\begin{align*}
& \left\| c - \ssum{k = 0}{n - 1}
    \bt^k (e_0 y q_0 y^* e_0) a \bt^k (e_0 y q_0 y^* e_0)
                    \right\|  \\
& \hspace*{3em} \mbox{}
  \leq \sum_{k = 0}^{n - 1} \| c_k -
      [e_0 y q_0 y^* e_0] \bt^{n - k} (a) [e_0 y q_0 y^* e_0] \|
  < n (4 n^3 + 3) \ep_0.
\end{align*}
Therefore
\[
\| c - q a q \| < [n^2 (8 n^3 + 7) + n (4 n^3 + 3)] \ep_0
   \leq 12 (n + 1)^5 \ep_0 < \ep.
\]

We regard $q$ as a \pj\  in $\Cs{n}{A}{\bt}$.
We further let
\[
D = \Csw{n}{E}{\bt |_E},
\]
which is a \fd\  subalgebra of $\Cs{n}{A}{\bt}$.
Let
\[
\ph \colon
 \Cs{n}{A}{\bt} \to \CZnAa
\]
be the isomorphism of Lemma~\ref{PertOfCrPrd}.
We take the \pj\  required in Lemma~\ref{TAFCond} to be $q$,
and the \fd\  subalgebra to be $\ph (D)$.
{}From what we just did, every element $a \in F$ satisfies
$\| [ q, a ] \| < \ep$ in $\Cs{n}{A}{\bt}$,
and there is $c \in E \S D$ such that $\| c - q a q \| < \ep$.
Since $\ph (q) = q$ and $\ph (a) = a$ for all $a \in F$,
in $\CZnAa$ every element $a \in F$ satisfies
$\| [ q, a ] \| < \ep$,
and there is $c \in E \S \ph (D)$ such that $\| c - q a q \| < \ep$.

Letting $v \in \Cs{n}{A}{\bt}$ be the
canonical unitary implementing the automorphism $\bt$,
we have $[q, v] = 0$, because $\bt (q) = q$, and $q v q \in D$.
Therefore in $\CZnAa$ we have
$[q, \, \ph (v)] = 0$ and $q \ph (v) q \in \ph (D)$.
Lemma~\ref{PertOfCrPrd} gives
\[
\| \ph (v) - u \| = \| w - 1 \| < 2 n^2 \ep_0.
\]
Therefore
\[
\| [q, u] \| < 4 n^2 \ep_0 < \ep
\andeqn
\| q u q - q \ph (v) q \| < 2 n^2 \ep_0 < \ep.
\]

We next show that $1 - q$ is \mvnt\  to a \pj\  in
${\overline{x \CZnAa x}}$.
Recall that
\[
1 - e \precsim p, \,\,\,\,\,\, 1 - q_0 \precsim p,
\andeqn (y q_0 y^*) e_0 = e_0 (y q_0 y^*).
\]
Furthermore,
\[
1 - q = 1 - \sum_{k = 0}^{n - 1} \bt^k (e_0 y q_0 y^* e_0)
      = 1 - e + \sum_{k = 0}^{n - 1} \bt^k (e_0 y [1 - q_0] y^* e_0).
\]
Now, with the \mvnc\  in $\Cs{n}{A}{\bt}$, we have
\[
\bt^k (e_0 y [1 - q_0] y^* e_0)
  \leq \bt^k (y [1 - q_0] y^*)
  \precsim 1 - q_0 \precsim p.
\]
Because
$\Cs{n}{A}{\bt} \cong \CZnAa$
via an isomorphism $\ph$ which fixes every element of $A$,
we get $\bt^k (e_0 y [1 - q_0] y^* e_0) \precsim p$
in $\CZnAa$ as well.
Thus, in $\CZnAa$,
the \pj\  $1 - q$ is the orthogonal sum of
$n + 1$ \pj s, each of which is \mvnt\  to a sub\pj\  of $p$.
By the choice of $p$, there are $n + 1$ \mops\  in
${\overline{x \CZnAa x}}$,
each \mvnt\  to $p$.
Therefore $1 - q$ is \mvnt\  to a \pj\  in
${\overline{x \CZnAa x}}$.

It remains to prove Condition~(4) of Proposition~\ref{TAFCond},
that is, that there are $N$ \mops\  in
$q \CZnAa q$,
each \mvnt\  in $\CZnAa$ to $1 - q$.
It suffices prove this in $A$ instead.
Since $A$ has cancellation of \pj s (by Theorem~\ref{TAFProp}),
it suffices to show that $N [1 - q] \leq [q]$ in $K_0 (A)$;
in fact, it suffices to show that $(N + 1) [1 - q] \leq [1]$.
We saw in the previous paragraph that
\[
1 - q = 1 - e + \sum_{k = 0}^{n - 1} \bt^k (e_0 y [1 - q_0] y^* e_0).
\]
By construction, we have $[1 - e] \leq [p] \leq [p_1]$ in $K_0 (A)$.
We also have
\begin{align*}
[\bt^k (e_0 y [1 - q_0] y^* e_0)]
 & \leq [\bt^k (y [1 - q_0] y^*)]
   = [\af^k (y [1 - q_0] y^*)]    \\
 & = [\af^k (1 - q_0)]
   \leq [\af^k (p)]
   \leq [p_1].
\end{align*}
Therefore $[1 - q] \leq (n + 1) [p_1]$.
Since $(N + 1)(n + 1) [p_1] \leq [p_0] \leq [1]$,
this gives $(N + 1) [1 - q] \leq [1]$, as desired.
\end{proof}

\section{Tracially approximately inner automorphisms}\label{Sec:TAI}

\indent
In this section, we introduce the notion of a
tracially approximately inner automorphism.
This condition is needed to prove that the dual action of an
action with the \aRp\  again has the \aRp.

Here, we only prove those results of immediate use.
Some further results are found in Section~\ref{Sec:More},
and some examples are in Section~\ref{Sec:Qst}.

\begin{dfn}\label{TAInnDfn}
Let $A$ be a stably finite \suca\  and let $\af \in \Aut (A)$.
We say that $\af$ is
{\emph{tracially approximately inner}} if for every finite set
$F \S A$, every $\ep > 0$, every $N \in \N$,
and every nonzero positive element $x \in A$,
there exist a \pj\  $e \in A$ and a unitary $v \in e A e$
such that:
\begin{itemize}
\item[(1)]
$\| \af (e) - e \| < \ep$.
\item[(2)]
$\| e a - a e \| < \ep$ for all $a \in F$.
\item[(3)]
$\| v e a e v^* - \af (e a e) \| < \ep$ for all $a \in F$.
\item[(4)]
$1 - e$ is \mvnt\  to a \pj\  in ${\overline{x A x}}$.
\item[(5)]
There are $N$ \mops\  $f_1, f_2, \dots, f_N \leq e$,
each of which is \mvnt\  to $1 - e$.
\end{itemize}
\end{dfn}

As in Definition~\ref{ARPDfn},
we allow $e = 1$, in which case conditions~(4) and~(5) are vacuous.

The motivation for the terminology is the same as that for
the tracial Rokhlin property (Definition~\ref{ARPDfn}).
As there,
the condition does seem useful outside the stably finite case,
so we include stable finiteness in the definition.
As with the \aRp, we also do not attempt to formulate the correct
version in the nonsimple case;
at the very least such an extension should include
the requirement that $\| e a_0 e \| > \| a_0 \| - \ep$ for a
predetermined nonzero $a_0 \in A$.

This condition seems to be appropriate for use with the \aRp.
In particular, when we strengthen it in the presence of the \aRp\  below,
we must allow $e \neq 1$ even if we start with an approximately inner
automorphism.

\begin{rmk}\label{InnAndTAInn}
Let $A$ be a stably finite \suca\  and let $\af \in \Aut (A)$.
If $\af$ is approximately inner
then $\af$ is tracially approximately inner.
If $\af$ is tracially approximately inner
and $A$ does not have Property~(SP),
then $\af$ is approximately inner.
\end{rmk}

Example~\ref{CAR2} shows that a tracially approximately inner
automorphism need not be approximately inner,
even on a simple AF algebra.

When $A$ has cancellation of \pj s, the automorphism has finite order,
and the action it generates has the \aRp, we can strengthen
Condition~(1) in Definition~\ref{TAInnDfn} to true invariance,
and we can require that $v$ have the same order as $\af$.
In Proposition~\ref{FinOrdTAInn2} below, we will further strengthen
this result, requiring for example $\af (v) = v$.
We do not know whether cancellation of \pj s is really necessary.

\begin{lem}\label{FinOrdTAInn}
Let $A$ be a stably finite \suca\  with cancellation of \pj s,
and let $\af \in \Aut (A)$
be tracially approximately inner and satisfy $\af^n = \id_A$.
Suppose that the action of $\Zqn$ generated by $\af$ has the
tracial Rokhlin property.
Then for every finite set
$F \S A$, every $\ep > 0$, every $N \in \N$,
and every nonzero positive element $x \in A$,
there exist a \pj\  $e \in A$ and a unitary $v \in e A e$
such that:
\begin{itemize}
\item[(1)]
$\af (e) = e$.
\item[(2)]
$\| e a - a e \| < \ep$ for all $a \in F$.
\item[(3)]
$v^n = e$,
and $\| v e a e v^* - \af (e a e) \| < \ep$ for all $a \in F$.
\item[(4)]
$1 - e$ is \mvnt\  to a \pj\  in ${\overline{x A x}}$.
\item[(5)]
There are $N$ \mops\  $f_1, f_2, \dots, f_N \leq e$,
each of which is \mvnt\  to $1 - e$.
\end{itemize}
\end{lem}

The proof requires a fair amount to technical work to set up a
rather short punch line, so we explain the basic idea before we begin.
Assume for simplicity that $F$ is $\af$-invariant,
that we can use the strict Rokhlin property to
obtain \pj s $e_0, \, e_1, \, \dots, \, e_{n - 1}$ which exactly
commute with every element of $F$ and such that $\af (e_j) = e_{j - 1}$
for $0 \leq j \leq n - 1$ and $\sum_{j = 0}^{n - 1} e_j = 1$,
and further that we can find a unitary $v_0 \in A$ such that
$v_0 a v_0^* = \af (a)$ for
$a \in F \cup \{ e_0, \, e_1, \, \dots, \, e_{n - 1} \}$.
Then
\[
v = e_1 v_0 e_0 + e_2 v_0 e_1 + \cdots + e_{n - 1} v_0 e_{n - 2}
      + e_0 v_0^{- (n - 1)} e_{n - 1}
\]
is a unitary in $A$ which satisfies $v^n = 1$ and $v a v^* = \af (a)$
for all $a \in F$.

We also note that there is some relation between
this proof and that of Lemma~3.3 of~\cite{Iz}.

\smallskip

\noindent
{\emph{Proof of Lemma~\ref{FinOrdTAInn}.}}
Let $F \S A$ be a finite set, let $\ep > 0$, let $N \in \N$,
and let $x \in A$ be a nonzero positive element.

\Wolog\  $n > 1$ and every $a \in F$ satisfies $\| a \| \leq 1$.
Then $A \not\cong M_m$ for any $m$, because all automorphisms of
$M_m$ are inner and therefore don't have the \aRp.
By Remark~\ref{InnAndTAInn} and Lemma~\ref{RImpSP}, there are
two cases:
either $A$ has Property~(SP) and $A \not\cong M_m$ for any $m$,
or $\af$ is approximately inner and has the strict Rokhlin property.
We write the proof in the first case.
In the second case, we take the \pj s
\[
1 - f_0, \, p, \, g_1, \, g_2, \, \dots, \, g_{2 N + 2}
\]
appearing below to be all zero, and we obtain the conclusion of the
lemma with in addition $e = 1$.

Lemma~3.2 of~\cite{LnTAF} provides nonzero \mvnt\   orthogonal
\pj s $g_1, \, g_2, \, \dots, \, g_{2 N + 2} \in {\overline{x A x}}$.

Define $\ep_1 = \frac{1}{10} n^{-1} \ep$.
Choose $\ep_2 > 0$ with
\[
\ep_2 \leq \min \left( \frac{\ep_1}{4 n + 1},
        \, \frac{\ep_1}{2 n (n + 1)} \right),
\]
and also so small that whenever $D$ is a unital \ca\  and
$c \in D$ satisfies $\| c c^* - 1 \| < 2 n \ep_2$ and
$\| c^* c - 1 \| < 2 n \ep_2$, then the unitary
$u = c (c^* c)^{-1/2} \in D$ satisfies $\| u - c \| < \ep_1$.

Apply Lemma~\ref{StARPDfn}
with $G_0 = \bigcup_{j = 0}^{n - 1} \af^j (F)$
in place of $F$,
with ${\ts{ \frac{1}{5} }} \ep_2$ in place of $\ep$,
with $1$ in place of $N$, and
with $g_1$ in place of $x$.
Let
\[
p_0^{(0)}, \, p_1^{(0)}, \, \dots, \, p_{n - 1}^{(0)} \in A
\]
be the resulting \pj s,
and let $p = \sum_{j = 0}^{n - 1} p_j^{(0)}$.
Note that $\af (p) = p$.

Apply Lemma~\ref{FixRel} with $E = \C^{n}$,
with $N = 1$, and with ${\ts{ \frac{1}{5} }} \ep_2$ in place of $\ep$.
Let $\ep_3$ be the resulting value of $\dt$.
Set
\[
\ep_4 = \min \left( {\ts{ \frac{1}{125} }} \ep_2,
         \, {\ts{ \frac{1}{25} }} \ep_3 \right).
\]
Apply Lemma~\ref{FixRel} with $E = \C$,
with $N = 1$, and with $\ep_4$ in place of $\ep$.
Let $\ep_5$ be the resulting value of $\dt$.
Choose $\ep_6 > 0$ with
\[
\ep_6 \leq \min \left( \frac{\ep_5}{4 n + 1}, \, \frac{\ep_4}{n},
  \, \frac{1}{2 n} \right),
\]
and also so small that whenever $D$ is a unital \ca\  and
$p, \, q \in D$ are \pj s such that $\| p - q \| < 2 n \ep_6$,
then there is a unitary $u \in D$ such that $u p u^* = q$
and $\| u - 1 \| < \ep_4$.

Apply Definition~\ref{TAInnDfn}
with
\[
G = \left\{ \rsz{ p, \, p_0^{(0)}, \, p_1^{(0)},
                \, \dots, \, p_{n - 1}^{(0)}  } \right\}
    \cup G_0
\]
in place of $F$,
with $\ep_6$ in place of $\ep$,
with $1$ in place of $N$, and
with $g_1$ in place of $x$.
Let $f_0$ be the resulting \pj\  and let $w_0 \in f_0 A f_0$
be the resulting unitary.

Set
\[
b = \frac{1}{n} \sum_{j = 0}^{n - 1} \af^j (f_0).
\]
Then
\[
\| b - f_0 \|
 \leq \frac{1}{n} \sum_{k = 0}^{n - 1} \sum_{j = 0}^{k - 1}
       \| \af^j ( \af (f_0) - f_0 ) \| < n \ep_6 < \ts{ \frac{1}{2} }.
\]
Therefore we can define $f_1 = \ch_{ (\frac{1}{2}, \I)} (b)$,
and $f_1$ is a \pj\  with $\af (f_1) = f_1$ and
\[
\| f_1 - f_0 \| \leq \| f_1 - b \| + \| b - f_0 \| \leq 2 \| b - f_0 \|
 < 2 n \ep_6.
\]
By the choice of $\ep_6$, there is a unitary $z_1 \in A$
such that $z_1 f_0 z_1^* = f_1$ and $\| z_1 - 1 \| < \ep_4$.
We now have
\[
\| f_1 p - p f_1 \| \leq 2 \| f_1 - f_0 \| + \| f_0 p - p f_0 \|
  < 4 n \ep_6 + \ep_6 = (4 n + 1) \ep_6 \leq \ep_5.
\]
By the choice of $\ep_5$ using Lemma~\ref{FixRel}, and
applying this lemma in the fixed point algebra $A^{\af}$
with $\ph (\ld) = \ld f_1$ and $p$ and $1 - p$ in place of
$e_0$ and $e_1$,
we obtain a unitary $z_2 \in A^{\af}$
such that $f = z_2 f_1 z_2^*$ commutes with $p$ and
$\| z_2 - 1 \| < \ep_4$.
Note that $\af (f) = f$.
Further define $w = z_2 z_1 w_0 (z_2 z_1)^*$.
We now have
\[
\| f - f_0 \| \leq \| f - f_1 \| + \| f_1 - f_0 \|
  \leq 2 \| z_2 - 1 \| + \| f_1 - f_0 \|
  < 2 \ep_4 + 2 n \ep_6 \leq 4 \ep_4.
\]
and
\[
\| w - w_0 \| \leq 2 \| z_2 z_1 - 1 \| < 4 \ep_4.
\]
Since $\ep_6 \leq \ep_4$, we now obtain the following in place of the
conditions from Definition~\ref{TAInnDfn}.
(We give the proofs for~(2) and~(3) afterwards.)
\begin{itemize}
\item[(1)]
$\af (f) = f$.
\item[(2)]
$\| f a - a f \| < 25 \ep_4$ for all $a \in G$, and $f p = p f$.
\item[(3)]
The unitary $w \in f A f$ satisfies
$\| w f a f w^* - \af (f a f) \| < 25 \ep_4$ for all $a \in G$.
\item[(4)]
$1 - f \precsim g_1$.
\end{itemize}
For~(2), we observe that, because $\| a \| \leq 1$,
\[
\| f a - a f \| \leq \| f_0 a - a f_0 \| + 2 \| f - f_0 \|
  < \ep_6 + 8 \ep_4 \leq  25 \ep_4.
\]
For~(3), we observe that
\begin{align*}
\| w f a f w^* - \af (f a f) \|
 & \leq \| w_0 f_0 a f_0 w_0^* - \af (f_0 a f_0) \|
           + 2 \| w - w_0 \| + 4 \| f - f_0 \|  \\
 & < \ep_6 + 8 \ep_4 + 16 \ep_4 \leq 25 \ep_4.
\end{align*}

We use the choice of $\ep_3$ from Lemma~\ref{FixRel} and the
inequality $25 \ep_4 \leq \ep_3$, and apply this lemma in $p A p$
with $\ph \colon \C^{n} \to p A p$ being
\[
\ph ( \ld_0, \dots, \ld_{n - 1})
   = \sum_{j = 0}^{n - 1} \ld_j p_j^{(0)}
\]
and with $p f$ and $p (1 - f)$ in place of $e_0$ and $e_1$.
We obtain a unitary $y_0 \in p A p$ such that the
unitary $y = y_0 + 1 - p \in A$ has the property that
$p_j = y p_j^{(0)} y^*$ commutes with $f$
for $0 \leq j \leq n - 1$,
and that $\| y - 1 \| < {\ts{ \frac{1}{5} }} \ep_2$
and $\sum_{j = 0}^{n - 1} p_j = p$.
Then
$\left\| \rsz{ p_j - p_j^{(0)} } \right\| \!\! \rule{0em}{2.5ex}
    < {\ts{ \frac{2}{5} }} \ep_2$.
{}From
$\left\| \rsz{ a p_j^{(0)} - p_j^{(0)} a} \right\| \!\!\rule{0em}{2.5ex}
         < {\ts{ \frac{1}{5} }} \ep_2$
for $a \in G_0$,
we now get $\| a p_j - p_j a \| < \ep_2$,
and from
$\left\| \af \left( \rsz{ p_j^{(0)} } \right)
      - \rsz{ p_{j + 1}^{(0)} } \right\| \!\! \rule{0em}{2.5ex}
              < {\ts{ \frac{1}{5} }} \ep_2$
we get $\| \af (p_j) - p_{j + 1} \| < \ep_2$.
Moreover,
\begin{align*}
\| w f p_j f w^* - \af (f p_j f) \|
&  \leq \ts{ \left\| \rsz{ w f p_j^{(0)} f w^* }
             - \af \left( \rsz{ f p_j^{(0)} f } \right) \right\| }
     + 2 \ts{ \left\| \rsz{  p_j - p_j^{(0)} }  \right\| }    \\
&  < 25 \ep_4 + {\ts{ \frac{4}{5} }} \ep_2 \leq \ep_2.
\end{align*}
Therefore $\| w f p_j f w^* - f p_{j + 1} f \| < 2 \ep_2$
for $0 \leq j \leq n - 1$.
Furthermore, for $a \in G_0$ we have
\[
\| f p_j a - a f p_j \|
  \leq \| f \| \cdot \| p_j a - a p_j \|
               + \| f a - a f \| \cdot \| p_j \|
  < \ep_2 + 25 \ep_4 \leq 2 \ep_2.
\]
We also estimate
\begin{align*}
\| w f p - f p w \|
 & \leq \sum_{j = 0}^{n - 1} \| w f p_j - f p_{j + 1} w \|  \\
 & = \sum_{j = 0}^{n - 1} \| w f p_j f w^* - f p_{j + 1} f \|
   < 2 n \cdot 2 \ep_2 = 4 n \ep_2.
\end{align*}
Using $w \in f A f$, for $a \in G_0$ we have
\begin{align*}
\| w (p f) a (p f) w^* - \af ((p f) a (p f)) \|
& = \| w p f a f p w^* - p \af (f a f) p \| \\
& \leq \| p [w f a f w^* - \af (f a f)] p \| + 2 \| w f p - f p w \|  \\
& < 25 \ep_4 + 4 n \ep_2 \leq (4 n + 1) \ep_2.
\end{align*}

Define $e = f p$ and $e_j = f p_j$.
Then
\[
e = \sum_{j = 0}^{n - 1} e_j \andeqn \| w e - e w \| < 4 n \ep_2,
\]
and for $0 \leq j \leq n - 1$ and $a \in G_0$ we have
\[
\| w e_j w^* - e_{j + 1} \| < 2 \ep_2, \,\,\,\,\,\,
\| \af (e_j) - e_{j + 1} \| < \ep_2,
\]
\[
\| w e a e w^* - \af (e a e) \| < (4 n + 1) \ep_2 \leq \ep_1, \andeqn
\| e_j a - a e_j \| < 2 \ep_2.
\]

Since $p$ and $f$ are $\af$-invariant, so is $e$, which is
Condition~(1) of the conclusion.
Also, if $a \in F$ then
\[
\| e a - a e \| \leq \sum_{j = 0}^{n - 1} \| e_j a - a e_j \| < 2 n \ep_2
 \leq \ep_1 \leq \ep,
\]
so we have Condition~(2) of the conclusion.

Following Convention~\ref{ARPConv}, set
$c = \sum_{j = 0}^{n - 1} e_{j + 1} w e_j \in e A e$.
Then
\[
\| c^* c - e \|
  \leq \sum_{j = 0}^{n - 1} \| e_j w^* e_{j + 1} w e_j - e_j \|
  \leq \sum_{j = 0}^{n - 1} \| w^* e_{j + 1} w - e_j \|
  < 2 n \ep_2.
\]
Similarly $\| c c^* - e \| < 2 n \ep_2$.
By the choice of $\ep_2$,
the unitary $v_0 = c (c^* c)^{-1/2} \in e A e$
(functional calculus in $e A e$) satisfies
$v_0 v_0^* = v_0^* v_0 = e$ and $\| v_0 - c \| < \ep_1$.
Also,
\begin{align*}
\| c - w e \|
 & \leq \| w e - e w \| + \| c - e w e \|  \\
 & \leq \| w e - e w \|
   + \sum_{j = 0}^{n - 1} \,\, \sum_{1 \leq k \leq n, \, k \neq j + 1}
                           \| e_k w e_j \|  \\
 & < 4 n \ep_2 + n (n - 1) \cdot 2 \ep_2 = 2 n (n + 1) \ep_2 \leq \ep_1.
\end{align*}
Therefore $\| v_0 - w e \| < \ep_1 + \ep_1 = 2 \ep_1$.
So, if $a \in G_0$ then
\[
\| v_0 e a e v_0^* - \af (e a e) \|
   \leq 2 \| v_0 - w e \| + \| w e a e w^* - \af (e a e) \|
   < 4 \ep_1 + \ep_1 = 5 \ep_1.
\]
Moreover,
since $c^* c$ commutes with all $e_j$, so does $(c^* c)^{-1/2}$.
{}From the relation $c = \sum_{j = 0}^{n - 1} e_{j + 1} c e_j$,
we then get $v_0 = \sum_{j = 0}^{n - 1} e_{j + 1} v_0 e_j$.
Since $v_0$ is unitary, this implies
$v_0 e_j v_0^* = e_{j + 1}$ for $0 \leq j \leq n - 1$.

Now define
\[
v = e_1 v_0 e_0 + e_2 v_0 e_1 + \cdots + e_{n - 1} v_0 e_{n - 2}
      + e_0 v_0^{- (n - 1)} e_{n - 1}.
\]
Then $v$ is a unitary in $e A e$ such that $v^n = e$.
Let $a \in F$.
We need to estimate $\| v e a e v^* - \af (e a e) \|$.
Set $b = \sum_{j = 0}^{n - 1} e_j a e_j$.
Then
\begin{align*}
\| b - e a e \|
&  \leq \sum_{j = 0}^{n - 1} \,\, \sum_{0 \leq k \leq n - 1, \, k \neq j}
             \| e_k a e_j \|   \\
&  \leq \sum_{j = 0}^{n - 1} \,\, \sum_{0 \leq k \leq n - 1, \, k \neq j}
             \| e_k a - a e_k \| \cdot \| e_j \|
  < 2 n (n - 1) \ep_2.
\end{align*}
Next, we calculate
\begin{align*}
v b v^* & = \sum_{j = 0}^{n - 1} v e_j a e_j v^*
    = v_0^{- (n - 1)} e_{n - 1} a e_{n - 1} v_0^{n - 1}
      + \sum_{j = 0}^{n - 2} v_0 e_j a e_j v_0^*   \\
  & = e_0 v_0^{- (n - 1)} e a e v_0^{n - 1} e_0
      + \sum_{j = 0}^{n - 2} e_{j + 1} v_0 e a e v_0^* e_{j + 1}.
\end{align*}
Now $\| v_0 e a e v_0^* - \af (e a e) \| < 5 \ep_1$ since $a \in G_0$.
Also, since $\af^n = \id_A$ it follows that all
$\af^k (a)$, for $k \in \Z$, are in $G_0$ as well, so that an inductive
argument gives
$\| v_0^k \af (e a e) v_0^{-k} - \af^{k + 1} (e a e) \| < 5 k \ep_1$
for $k \geq 1$.
Putting $k = n - 1$, using $\af^n = \id_A$, and conjugating by
$v_0^{- (n - 1)}$, we obtain
\[
\| \af (e a e) - v_0^{- (n - 1)} e a e v_0^{n - 1} \|
          < 5 (n - 1) \ep_1.
\]
Therefore
\begin{align*}
& \left\| v b v^* - \ssum{j = 0}{n - 1} e_j \af (e a e) e_j \right\| \\
& \hspace*{3em} \mbox{}
   \leq \|e_0 [ v_0^{- (n - 1)} e a e v_0^{n - 1} - \af (e a e) ] e_0 \|
       + \sum_{j = 0}^{n - 2}
          \| e_{j + 1} [v_0 e a e v_0^* - \af (e a e) ] e_{j + 1} \|  \\
& \hspace*{3em} \mbox{}
   < 5 (n - 1) \ep_1 + (n - 1) \cdot 5 \ep_1
   = 10 (n - 1) \ep_1.
\end{align*}
On the other hand,
\[
\af (b) = \sum_{j = 0}^{n - 1} \af (e_j a e_j)
   = \sum_{j = 0}^{n - 1} \af (e_j) \af (e a e) \af (e_j).
\]
Therefore
\[
\left\| \af (b) - \ssum{j = 0}{n - 1} e_j \af (e a e) e_j \right\|
  \leq \sum_{j = 0}^{n - 1} 2 \| \af (e_j) - e_{j + 1} \|
  < 2 n \ep_2 \leq \ep_1.
\]
Putting everything together, we get
$\| v b v^* - \af (b) \| < (10 n - 9) \ep_1$, so
\begin{align*}
\| v e a e v^* - \af (e a e) \|
 & \leq \| v b v^* - \af (b) \| + 2 \| b - e a e \| \\
 & < (10 n - 9) \ep_1 + 4 n (n - 1) \ep_2 \leq (10 n - 7) \ep_1 < \ep.
\end{align*}
This, together with the relation $v^n = 1$ from above,
is Condition~(3) of the conclusion.

To prove Condition~(4) of the conclusion,
write $1 - e = 1 - f + f (1 - p)$.
By construction,
\[
1 - f \precsim g_1, \,\,\,\,\,\, 1 - p \precsim g_1, \andeqn
g_1 \sim g_2.
\]
Therefore $1 - e \precsim g_1 + g_2 \in {\overline{x A x}}$,
as required.

We prove Condition~(5) of the conclusion.
First note that
\[
1 - f \precsim g_1 \leq 1 - \sum_{j = 3}^{2 N + 2} g_j
\andeqn
1 - p \precsim g_2 \leq 1 - \sum_{j = 3}^{2 N + 2} g_j,
\]
so
\[
1 - e = 1 - f + f (1 - p)
 \precsim g_1 + g_2 \leq 1 - \sum_{j = 3}^{2 N + 2} g_j.
\]
Therefore, because $A$ has cancellation of \pj s,
$\sum_{j = 3}^{2 N + 2} g_j \precsim e$.
Choose a partial isometry $s \in A$ such that
\[
s^* s = \sum_{j = 3}^{2 N + 2} g_j \andeqn s s^* \leq e.
\]
For $1 \leq j \leq N$ let $q_j^{(1)} \leq g_{j + 2}$
be a \pj\  with $q_j^{(1)} \sim 1 - f$ and let
$q_j^{(2)} \leq g_{j + N + 2}$
be a \pj\  with $q_j^{(2)} \sim 1 - p$.
Set $q_j = q_j^{(1)} + q_j^{(2)}$.
Then the \pj s $s q_j s^*$ are $N$ \mops\  in $e A e$,
each of which satisfies $1 - e \precsim s q_j s^*$.
We have shown that there are $N$ \mops\  in $e A e$, each of which is
\mvnt\  to $1 - e$, as desired.
\QED

\smallskip

Our further strengthening of Definition~\ref{TAInnDfn}
requires two preliminary lemmas.

\begin{lem}\label{StabOfFinOrdU}
Let $D$ be a unital \ca, let $n \in \N$, and let $\ep > 0$.
Then there is $\dt > 0$ such that whenever $v \in D$ is a unitary such
that $v^n = 1$, whenever $B$ is a unital C*-subalgebra of $D$, and
whenever $c \in B$ satisfies $\| c - v \| < \dt$, then there
is a unitary $w \in B$ such that $w^n = 1$ and $\| w - v \| < \ep$.
\end{lem}

\begin{proof}
This is semiprojectivity of $\C^n$,
which is the universal \ca\  generated by a unitary
$v$ with $v^n = 1$.
(See Chapter~14 of~\cite{Lr}.)
\end{proof}

\begin{lem}\label{StabOfFinOrdU2}
Let $D$ be a unital \ca, let $n \in \N$, let $\af \in \Aut (D)$
satisfy $\af^n = \id_D$, and let $\ep > 0$.
Then there is $\dt > 0$ such that whenever $v \in D$ is a unitary such
that $v^n = 1$, and whenever $c \in D$ satisfies $\| c - v \| < \dt$
and $\af (c) = \exp (- 2 \pi i / n) c$, then there
is a unitary $w \in D$ such that
\[
\| w - v \| < \ep, \,\,\,\,\,\, w^n = 1, \andeqn
\af (w) = \exp (- 2 \pi i / n) w.
\]
\end{lem}

\begin{proof}
Let $\om = \exp (2 \pi i / n)$.
Define an open set $U \S S^1$ by
\[
U = S^1 \SM \{ \exp (\pi i / n) \om^j \colon 0 \leq j \leq n - 1 \}.
\]
Let $f \colon U \to \C$ be the \cfn\  which takes the constant value
$\om^j$ on the open arc from
$\exp (- \pi i / n) \om^j$ to $\exp (\pi i / n) \om^j$.
Choose $\ep_0 > 0$ with $\ep_0 \leq \frac{1}{2} \ep$ and so small that
whenever $z \in D$ is a unitary such
that $\| z^n - 1 \| < \ep_0$, then $\spec (z) \S U$ and
$\| f (z) - z \| < \frac{1}{2} \ep$.
Choose $\dt > 0$ with $\dt \leq \frac{1}{2} \ep$ and so small that
whenever $v \in D$ is a unitary and whenever
$c \in D$ satisfies $\| c - v \| < \dt$, then
$\| c (c^* c)^{-1/2} - v \| < {\ts{ \frac{1}{n} }} \ep_0$.

Now let $v \in D$ be a unitary such
that $v^n = 1$, and let $c \in D$ satisfy $\| c - v \| < \dt$
and $\af (c) = \om^{-1} c$.
By the choice of $\dt$,
the unitary $z = c (c^* c)^{-1/2} \in D$ satisfies
$\| z - v \| < {\ts{ \frac{1}{n} }} \ep_0$.
Therefore
\[
\| z^n - 1 \|
    = \| z^n - v^n \|
    \leq \sum_{k = 1}^n
       \| z \|^{n - k} \cdot \| z - v \| \cdot \| v \|^{k - 1}
    < n \cdot \left( {\ts{ \frac{1}{n} }} \ep_0 \right) = \ep_0.
\]
By the choice of $\ep_0$, the unitary $w = f (z)$
satisfies $\| w - z \| < \frac{1}{2} \ep$.
Therefore
\[
\| w - v \| \leq \| w - z \| + \| z - v \|
  < \ts{ \frac{1}{2} } \ep + {\ts{ \frac{1}{2 n} }} \ep \leq \ep,
\]
as desired.

It is clear that $w^n = 1$.

It remains only to show that $\af (w) = \om^{-1} w$.
Since $\af (c^*) = \om c^*$, it follows that
$c^* c$ is in the fixed point algebra $D^{\af}$.
Therefore $\af (z) = \om^{-1} z$.
Since $f ( \om^{-1} \zt) = \om^{-1} f (\zt)$ for every $\zt \in U$,
we get $f (\om^{-1} z) = \om^{-1} f (z)$.
Therefore
\[
\af (w) = \af ( f (z)) = f ( \af (z))
    = f (\om^{-1} z) = \om^{-1} f (z) = \om^{-1} w,
\]
as desired.
\end{proof}

\begin{prp}\label{FinOrdTAInn2}
Let $A$ be a stably finite \suca\  with cancellation of \pj s,
and let $\af \in \Aut (A)$
be tracially approximately inner and satisfy $\af^n = \id_A$.
Suppose that the action of $\Zqn$ generated by $\af$ has the
tracial Rokhlin property.
Then for every finite set
$F \S A$, every $\ep > 0$, every $N \in \N$,
and every nonzero positive element $x \in A$,
there exist a \pj\  $e \in A$ and unitaries $v_1, \, v_2 \in e A e$
such that the Conditions (1)--(5) of Lemma~\ref{FinOrdTAInn}
are satisfied for both $v_1$ and $v_2$ in place of $v$, and in addition
\[
\af (v_1) = v_1, \,\,\,\,\,\,
\af (v_2) = \exp (- 2 \pi i / n) v_2, \andeqn
\| v_1 v_2 - v_2 v_1 \| < \ep.
\]
\end{prp}

As for Lemma~\ref{FinOrdTAInn},
the proof requires a fair amount to technical work to set up a
rather short punch line, and we explain the basic idea before we begin.
Assume for simplicity that $F$ is $\af$-invariant,
that we got from Lemma~\ref{FinOrdTAInn} a unitary $v_0 \in A$
such that $v_0^n = 1$ and $v_0 a v_0^* = \af (a)$ for all $a \in F$,
and further that we can use the strict Rokhlin property to
obtain \pj s $e_0, \, e_1, \, \dots, \, e_{n - 1}$ which exactly
commute with $v_0$ and every element of $F$
and such that $\af (e_j) = e_{j + 1}$
for $0 \leq j \leq n - 1$ and $\sum_{j = 0}^{n - 1} e_j = 1$.
(We apply the approximate innerness condition and the
Rokhlin property in the opposite order from the proof
of Lemma~\ref{FinOrdTAInn}.)
Then
\[
v_1 = \sum_{j = 0}^{n - 1} \af^j (e_0 v_0 e_0)
\andeqn
v_2 = \sum_{j = 0}^{n - 1} \exp (- 2 \pi i j / n) \af^j (e_0 v_0 e_0)
\]
are unitaries in $A$ which satisfy the conclusion.

\smallskip

\noindent
{\emph{Proof of Proposition~\ref{FinOrdTAInn2}.}}
Let $F \S A$ be a finite set, let $\ep > 0$, let $N \in \N$,
and let $x \in A$ be a nonzero positive element.
\Wolog\  $\af (F) = F$ and $\| a \| \leq 1$ for all $a \in F$.

As in the proof of Lemma~\ref{FinOrdTAInn}, there are two cases:
either $A$ has Property~(SP) and $A \not\cong M_m$ for any $m$,
or $\af$ is approximately inner and has the strict Rokhlin property.
We again write the proof in the first case, and the second case
is handled the same way as there.

Lemma~3.2 of~\cite{LnTAF} provides nonzero \mvnt\   orthogonal
\pj s $g_1, \, g_2, \dots, g_{2 N + 2} \in {\overline{x A x}}$.

Choose $\ep_1 > 0$ with $n^2 \ep_1 \leq \ts{ \frac{1}{13} } \ep$
and so small that $2 n \ep_1$ will serve as $\dt$ in
Lemmas~\ref{StabOfFinOrdU} and~\ref{StabOfFinOrdU2}
with the given value of $n$ and with $\ts{ \frac{1}{13} } \ep$
in place of $\ep$.
Choose $\ep_2 > 0$ with
$\ep_2 \leq \min ( (n + 1)^{-1} \ep, \, \ep_1)$
and so small that whenever $D$ is a unital \ca\  and
$p, \, q \in D$ are \pj s such that $\| p - q \| < n \ep_2$,
then there is a unitary $u \in D$ such that $u p u^* = q$
and $\| u - 1 \| < \ep_1$.

Apply Lemma~\ref{FinOrdTAInn}
with $F$ as given,
with $\ts{ \frac{1}{6} } \ep_2$ in place of $\ep$,
with $1$ in place of $N$,
and with $g_1$ in place of $x$.
We obtain an $\af$-invariant \pj\  $f_0$
and a unitary $w_0 \in f_0 A f_0$ such that, in particular,
$w_0^n = f_0$ and
\[
\| a f_0 - f_0 a \| < \ts{ \frac{1}{6} } \ep_2
\andeqn
\| w_0 (f_0 a f_0) w_0^* - \af (f_0 a f_0) \| < \ts{ \frac{1}{6} } \ep_2
\]
for all $a \in F$.

We want to use Lemma~\ref{StARPDfn} for the automorphism $\af$,
but we also want the resulting \pj s
to commute with $f_0$ and $w_0$,
and we want to retain $\af$-invariance of $f_0$ and
the property $w_0^n = 1$.
This requires perturbation both of the Rokhlin \pj s
and of $f_0$ and $w_0$.
To this end, we apply Lemma~\ref{FixRel} three times.
The first time, we take $E = \C^n$ and use $n - 1$ in place of $N$
and $\ts{ \frac{1}{30} } \ep_2$ in place of $\ep$.
Let $\rh$ be the resulting value of $\dt$, and set
$\ep_3 = \min \left( {\ts{ \frac{1}{n} }} \rh, \,
    {\ts{ \frac{1}{30} }} \ep_2 \right)$.
The second time, we take $E = \C^2$ and use $n - 1$ in place of $N$
and ${\ts{ \frac{1}{5} }} \ep_3$ in place of $\ep$.
Let $\ep_4$ be the resulting value of $\dt$, and also require
$\ep_4 \leq {\ts{ \frac{1}{5} }} \ep_3$.
The third time, we take $E = \C^2$ and use $1$ in place of $N$
and ${\ts{ \frac{1}{5} }} \ep_4$ in place of $\ep$.
Let $\ep_5$ be the resulting value of $\dt$, and also require
$\ep_5 \leq {\ts{ \frac{1}{5} }} \ep_4$.

Apply Lemma~\ref{StARPDfn} to $\af$,
with $F \cup \{ f_0, w_0 \}$ in place of $F$,
with $\frac{1}{n} \ep_5$ in place of $\ep$,
with $1$ in place of $N$,
and with $g_1$ in place of $x$.
We obtain \mops\  %
$p^{(0)}_0, \, p^{(0)}_1, \, \dots, \, p^{(0)}_{n - 1} \in A$.
Set $p = \sum_{j = 0}^{n - 1} p^{(0)}_j$.
Then in particular
\[
\| p f_0 - f_0 p \|
 \leq \sum_{j = 0}^{n - 1} \| p^{(0)}_j f_0 - f_0 p^{(0)}_j \|
  < n \left( \ts{ \frac{1}{n} } \ep_5 \right) = \ep_5.
\]

By the choice of $\ep_5$ using Lemma~\ref{FixRel}, and
applying this lemma in the fixed point algebra $A^{\af}$
with $\ph (\ld_1, \ld_2) = \ld_1 p + \ld_2 (1 - p)$
and $f_0$ and $1 - f_0$ in place of $e_0$ and $e_1$,
we obtain a unitary $z_1 \in A^{\af}$
such that $f = z_1 f_0 z_1^*$ commutes with $p$ and
$\| z_1 - 1 \| < {\ts{ \frac{1}{5} }} \ep_4$.
Note that $\af (f) = f$.
Moreover, $\| f - f_0 \| \leq 2 \| z_1 - 1 \|$, so
\begin{align*}
\| f p p^{(0)}_j - p^{(0)}_j f p \|
 & = \| (f p^{(0)}_j - p^{(0)}_j f) p \|
   \leq \| f_0 p^{(0)}_j - p^{(0)}_j f_0 \| + 4 \| z_1 - 1 \| \\
 & < {\ts{ \frac{1}{n} }} \ep_5 + {\ts{ \frac{4}{5} }} \ep_4
   \leq \ep_4.
\end{align*}
Also set $w_1 = z_1 w_0 z_1^* \in f A f$, so that $w_1^n = f$ and
\[
\| w_1 - w_0 \| \leq 2 \| z_1 - 1 \| < {\ts{ \frac{2}{5} }} \ep_4.
\]

By the choice of $\ep_4$ using Lemma~\ref{FixRel}, and
applying this lemma in the corner $p A p$
with $\ph (\ld_1, \ld_2) = \ld_1 f p + \ld_2 (1 - f) p$
and $p^{(0)}_0, \, p^{(0)}_1, \, \dots, \, p^{(0)}_{n - 1}$
in place of $e_0, e_1, \dots, e_N$,
we obtain a unitary $z_2 \in p A p$
such that $p_j = z_2 p^{(0)}_j z_2^*$ commutes with $f p$ and
$\| z_2 - p \| < {\ts{ \frac{1}{5} }} \ep_3$.
Because $p_j \leq p$ and $f$ commutes with $p$,
it now follows that
$p_j$ commutes with $f$ for $0 \leq j \leq n - 1$.
Moreover, $\| p_j - p^{(0)}_j \| \leq 2 \| z_2 - p \|$, so
\begin{align*}
\| w_1 p_j - p_j w_1 \|
&  \leq \| w_0 p^{(0)}_j - p^{(0)}_j w_0 \| + 2 \| p_j - p^{(0)}_j \|
                 + 2 \| w_1 - w_0 \|   \\
&  < \ep_5 + {\ts{ \frac{4}{5} }} \ep_3 + {\ts{ \frac{4}{5} }} \ep_4
   \leq {\ts{ \frac{1}{25} }} \ep_3 + {\ts{ \frac{4}{5} }} \ep_3
                              + {\ts{ \frac{4}{25} }} \ep_3
   = \ep_3.
\end{align*}
It follows that $\| w_1^k p_j - p_j w_1^k \| < k \ep_3$ for $k \geq 1$.
Therefore, with $\om = \exp (2 \pi i / n)$, the \pj s
\[
q_l = \frac{1}{n} \sum_{k = 0}^{n - 1} \om^{-l k} w_1^k
\]
satisfy $\| q_l p_j - p_j q_l \| < n \ep_3$.
Since $f$ commutes with $w_1$ and the $p_j$,
we get $\| q_l f p_j - f p_j q_l \| < n \ep_3$.
Let $\ph \colon \C^n \to f A f$ be the \hm\  which sends
$(1, \, \om, \, \dots, \, \om^{n - 1})$ to $w_1$.
Then the $q_l$ are the images under $\ph$
of the matrix units of $\C^n$.
By the choice of $\ep_3$ using Lemma~\ref{FixRel}, and
applying this lemma in the corner $f A f$ with this $\ph$
and with $f p_0, \, f p_1, \, \dots, \, f p_{n - 1}$
in place of $e_0, e_1, \dots, e_N$,
we obtain a unitary $z_3 \in f A f$
(which is $u^*$ in the conclusion of the lemma)
such that $w = z_3 w_1 z_3^*$ commutes with
$f p_0, f p_1, \dots, f p_{n - 1}$ and
$\| z_3 - f \| < {\ts{ \frac{1}{30} }} \ep_2$.
Since $w \in f A f$ and $f p_j = p_j f$,
we also get $w p_j = p_j w$.

We now have the algebraic relations
\[
\af (f) = f, \,\,\,\,\,\,
w \in f A f, \,\,\,\,\,\,
w^n = f, \,\,\,\,\,\,
p = \sum_{j = 0}^{n - 1} p_j,
\]
and, for $0 \leq j \leq n - 1$,
\[
p_j w = w p_j \andeqn p_j f = f p_j.
\]
We further claim that if $a \in F$ then
\[
\| f a - a f \| < \ep_2 \andeqn
\| w (f a f) w^* - \af (f a f) \| < \ep_2,
\]
and, for $0 \leq j \leq n - 1$,
\[
\| p_j a - a p_j \| < \ep_2 \andeqn
\| \af (p_j) - p_{j + 1} \| < \ep_2.
\]

We prove the claim.
For $a \in F$, we have
\begin{align*}
\| f a - a f \|
& \leq 2 \| f - f_0 \| + \| f_0 a - a f_0 \|
  \leq 4 \| z_1 - 1 \| + \| f_0 a - a f_0 \|  \\
& < \ts{ \frac{4}{5}} \ep_4 + \ts{ \frac{1}{6}} \ep_2 \leq \ep_2.
\end{align*}
This is the first estimate.
For the second,
\begin{align*}
\| w (f a f) w^* - \af (f a f) \|
& \leq 2 \| w - w_0 \| + 4 \| f - f_0 \|
   + \| w_0 (f_0 a f_0) w_0^* - \af (f_0 a f_0) \|  \\
& < 4 \| z_3 - f \| + 2 \| w_1 - w_0 \| + \ts{ \frac{8}{5}} \ep_4
   + \ts{ \frac{1}{6}} \ep_2  \\
& < \ts{ \frac{4}{30}} \ep_2 + \ts{ \frac{12}{5}} \ep_4
     + \ts{ \frac{5}{6}} \ep_2
  < \ep_2.
\end{align*}
For the third,
\begin{align*}
\| p_j a - a p_j \|
& \leq 2 \| p_j - p_j^{(0)} \| + \| p_j^{(0)} a - a p_j^{(0)} \|
  \leq 4 \| z_2 - p \| + \| p_j^{(0)} a - a p_j^{(0)} \|   \\
& < \ts{ \frac{4}{5}} \ep_3 + \ts{ \frac{1}{n}} \ep_5
  \leq \ep_3 \leq \ep_2.
\end{align*}
Finally,
\begin{align*}
\| \af (p_j) - p_{j + 1} \|
& \leq \| p_j - p_j^{(0)} \| + \| p_{j + 1} - p_{j + 1}^{(0)} \|
     + \| \af (p_j^{(0)}) - p_{j + 1}^{(0)} \|   \\
& \leq 4 \| z_2 - p \| + \| \af (p_j^{(0)}) - p_{j + 1}^{(0)} \|
  < \ts{ \frac{4}{5}} \ep_3 + \ts{ \frac{1}{n}} \ep_5
  \leq \ep_3 \leq \ep_2.
\end{align*}

We now define $e = f p$.
Then $\af (e) = e$, which is Condition~(1) of the conclusion of
Lemma~\ref{FinOrdTAInn}.
Furthermore, for $a \in F$ we have
\[
\| e a - a e \|
  \leq \| f \| \sum_{j = 0}^{n - 1} \| p_j a - a p_j \|
          + \| f a - a f \| \cdot \| p \|
  < (n + 1) \ep_2 \leq \ep.
\]
This is Condition~(2) of the conclusion of
Lemma~\ref{FinOrdTAInn}.

Next, define $v_0 = w e$,
which is a unitary in $e A e$ with $v_0^n = e$,
and $e_j = e p_j$ for $0 \leq j \leq n - 1$,
which are \pj s in $e A e$ which commute with $v_0$.
Since $\| \af (p_j) - p_{j + 1} \| < \ep_2$ and $\af (e) = e$,
we get $\| \af^j (e_j) - e_{j + 1} \| < \ep_2$
for $0 \leq j \leq n - 1$.
So
\[
\| \af^j (e_0) - e_j \| < j \ep_2 \leq n \ep_2
\]
for $0 \leq j \leq n - 1$.
By the choice of $\ep_2$, there are unitaries $y_j \in e A e$
with $\| y_j - e \| < \ep_1$ such that $y_j \af^j (e_0) y_j^* = e_j$
for $0 \leq j \leq n - 1$.
Define
\[
c_1 = \sum_{j = 0}^{n - 1} y_j \af^j (e_0 v_0 e_0) y_j^*
\andeqn
d_1 = \sum_{j = 0}^{n - 1} \af^j (e_0 v_0 e_0).
\]
Because the $e_j$ are orthogonal, and because $v_0$ commutes
with the $e_j$ and satisfies $v_0^n = e$,
it follows that $c_1$ is a unitary in $e A e$ with $c_1^n = e$.
Moreover, $d_1$ is $\af$-invariant and satisfies
\[
\| c_1 - d_1 \| \leq \sum_{j = 0}^{n - 1} 2 \| y_j - e \| < 2 n \ep_1.
\]
By the choice of $\ep_1$ using Lemma~\ref{StabOfFinOrdU}, there is
a unitary $v_1$ in the fixed point algebra $(e A e)^{\af}$ such that
$v_1^n = 1$ and $\| v_1 - c_1 \| < \ts{ \frac{1}{13} } \ep$.
Further define
\[
c_2 = \sum_{j = 0}^{n - 1} \om^{j} y_j \af^j (e_0 v_0 e_0) y_j^*
\andeqn
d_2 = \sum_{j = 0}^{n - 1} \om^{j} \af^j (e_0 v_0 e_0).
\]
By the same argument as above,
$c_2$ is a unitary in $e A e$ with $c_2^n = e$,
and
\[
\af (d_2) = \om^{- 1} d_2  \andeqn  \| c_2 - d_2 \| < 2 n \ep_1.
\]
So from Lemma~\ref{StabOfFinOrdU2} we get a unitary $v_2 \in e A e$
such that
\[
v_2^n = 1, \,\,\,\,\, \| v_2 - c_2 \| < \ts{ \frac{1}{13} } \ep,
\andeqn \af (v_2) = \om^{-1} v_2.
\]
We also observe that $c_1 c_2 = c_2 c_1$.
{}From $\| v_1 - c_1 \| < \ts{ \frac{1}{13} } \ep$
and $\| v_2 - c_2 \| < \ts{ \frac{1}{13} } \ep$ we therefore
get $\| v_1 v_2 - v_2 v_1 \| < \ts{ \frac{4}{13} } \ep < \ep$.
We have proved the new conditions on both $v_1$ and $v_2$.

We now estimate
\[
\| v_1 (e a e) v_1^* - \af (e a e) \| \andeqn
\| v_2 (e a e) v_2^* - \af (e a e) \|
\]
for $a \in F$.
We begin by observing that
\[
\| e_j a - a e_j \|
 \leq \| f a - a f \| \cdot \| p_j \|
        + \| f \| \cdot \| p_j a - a p_j \|
 < 2 \ep_2.
\]
Then set $b = \sum_{j = 0}^{n - 1} e_j a e_j$.
Since the $e_j$ are orthogonal,
\[
\| e a e - b \|
  \leq \sum_{j = 0}^{n - 1} \,\, \sum_{0 \leq k \leq n - 1, \, k \neq j}
            \| e_j a - a e_j \|
  < 2 n (n - 1) \ep_2 \leq 2 n^2 \ep_1.
\]
Also,
\[
\left\| \af (b) - \ssum{j = 0}{n - 1} e_j \af (a) e_j \right\|
  \leq \sum_{j = 0}^{n - 1} 2 \| e_{j + 1} - \af (e_j) \|
  < 2 n \ep_2 \leq 2 n \ep_1.
\]
The next step is to estimate
$\left\| c_l b c_l^* - \sum_{j = 0}^{n - 1} e_j \af (a) e_j \right\|$.
First, use $y_j \af^j (e_0) y_j^* = e_j$ and $e_0 v_0 = v_0 e_0$ to
calculate
\[
\af^j (e_0 v_0 e_0) y_j^* e_j = \af^j (e_0 v_0 e_0) \af^j (e_0) y_j^* e
  = \af^j (e_0 v_0) y_j^* e.
\]
Also $\af (e) = e$ and
\[
\| v_0 e a e v_0^* - \af (e a e) \|
  = \| e w f a f w^* e - e \af (f a f ) e \| < \ep_2.
\]
Now $\af^{-j} (a) \in F$, so the above applies with
$\af^j (a)$ in place of $a$.
Using this fact at the last step, we get
\begin{align*}
& \| [y_j \af^j (e_0 v_0 e_0) y_j^*] e_j a e_j
                        [y_j \af^j (e_0 v_0 e_0) y_j^*]^*
               - e_j \af (a) e_j \|  \\
& \hspace*{2em} \mbox{}
  = \| \af^j (e_0 v_0) y_j^* e a e y_j \af^j (v_0^* e_0)
               - y_j^* e_j \af (e a e) e_j y_j \|  \\
& \hspace*{2em} \mbox{}
  \leq \| y_j^* e a e y_j - e a e \|  \\
& \hspace*{4em} \mbox{}
     + \| \af^j (e_0) \af^j (v_0 e \af^{-j} (a) e v_0^*) \af^j (e_0)
              - \af^j (e_0) y_j^* \af (e a e) y_j \af^j (e_0)  \|  \\
& \hspace*{2em} \mbox{}
  \leq 2 \| y_j - e \|
    + \| v_0 e \af^{-j} (a) e v_0^* - \af^{- j + 1} (e a e) \|
    + \| \af (e a e) - y_j^* \af (e a e) y_j \|  \\
& \hspace*{2em} \mbox{}
  \leq 2 \| y_j - e \|
    + \| v_0 [e \af^{-j} (a) e] v_0^* - \af (e \af^{-j} (a) e) \|
    + 2 \| y_j - e \|  \\
& \hspace*{2em} \mbox{}
  < 4 \ep_1 + \ep_2 \leq 5 \ep_1.
\end{align*}
Now for $l = 1, \, 2$ we get
\begin{align*}
& \left\| c_l b c_l^* - \ssum{j = 0}{n - 1} e_j \af (a) e_j \right\| \\
& \hspace*{3em} \mbox{}
  \leq \sum_{j = 0}^{n - 1}
    \| [y_j \af^j (e_0 v_0 e_0) y_j^*] e_j a e_j
                        [y_j \af^j (e_0 v_0 e_0) y_j^*]^*
               - e_j \af (a) e_j \|  \\
& \hspace*{3em} \mbox{}
  < 5 n \ep_1.
\end{align*}
Putting everything together, for $l = 1, \, 2$ we now get
\begin{align*}
& \| v_l (e a e) v_l^* - \af (e a e) \|  \\
& \hspace*{3em} \mbox{}
 \leq 2 \| v_l - c_l \| + 2 \| e a e - b \|
   + \| c_l b c_l^* - \af (b) \|  \\
& \hspace*{3em} \mbox{}
 \leq 2 \| v_l - c_l \| + 2 \| e a e - b \|  \\
& \hspace*{7em} \mbox{}
   + \left\| c_l b c_l^* - \ssum{j = 0}{n - 1} e_j \af (a) e_j \right\|
   + \left\| \af (b) - \ssum{j = 0}{n - 1} e_j \af (a) e_j \right\|
             \\
& \hspace*{3em} \mbox{}
 < \ts{ \frac{2}{13} } \ep + 4 n^2 \ep_1 + 5 n \ep_1 + 2 n \ep_1 < \ep.
\end{align*}
We now have Condition~(3) of Lemma~\ref{FinOrdTAInn}, for both
$v_1$ and $v_2$.

The proof of Conditions~(4) and~(5) of Lemma~\ref{FinOrdTAInn}
is just like in the proof of that lemma.
We have
\[
1 - e = 1 - f + f (1 - p) \precsim g_1 + g_2 \in {\overline{x A x}},
\]
which gives Condition~(4).
Moreover,
\[
1 - e \precsim g_1 + g_2 \leq 1 - \sum_{k = 3}^{2 N + 2} g_k,
\]
so $\sum_{k = 3}^{2 N + 2} g_k \precsim e$ by cancellation,
and there are $N$ \mops\  dominated by $\sum_{k = 3}^{2 N + 2} g_k$,
each \mvnt\  to $1 - e$.
\QED

\section{Duality}\label{Sec:Dual}

\indent
In this section we show that if an action of $\Zqn$ on a
\suca\  has the \aRp\  and and its generator is tracially
approximately inner, then the dual action also has the \aRp.
We need to impose some conditions on comparison of \pj s
in order to make the proofs work.
These will be automatically satisfied in the cases we are
interested in, because in those cases both the original algebra
and the crossed product satisfy the much stronger condition
that the order on $K_0$ is determined by traces.
The result is related to Lemma~3.8 of~\cite{Iz}.

At the end of this section,
we use these results to prove a version of Theorem~\ref{RokhTAF}
for the fixed point algebra in place of the crossed product.

The following cancellation condition is needed
in the main theorem of this section.

\begin{dfn}\label{WDivDfn}
Let $D$ be a unital \ca.
We say that the {\emph{weak divisibility property}} holds for
\pj s in $D$ if whenever
$N \in \N$ and $e_1, e_2, \dots, e_N$ and $f_1, f_2, \dots, f_{N + 1}$
are two sets of mutually \mvnt\  orthogonal \pj s in $D$ such that
$\sum_{k = 1}^{N + 1} f_k \precsim \sum_{k = 1}^N e_k$,
then $f_1 \precsim e_1$.
\end{dfn}

\begin{lem}\label{GetWDiv}
Let $A$ be a \suca\  with tracial rank zero.
Then the \pj s in $A$ have the weak divisibility property.
\end{lem}

\begin{proof}
\Wolog\  $e_1 \neq 0$.
For any tracial state $\ta$ on $A$, we have $\ta (e_1) > 0$ and
$(N + 1) \ta (f_1) \leq N \ta (e_1)$,
so $\ta (f_1) < \ta (e_1)$.
By Theorem~\ref{TAFProp},
the order on \pj s over $A$ is determined by traces,
so $f_1 \precsim e_1$.
\end{proof}

\begin{thm}\label{DualHasRokh}
Let $A$ be a \suca\  and let $\af \in \Aut (A)$
be tracially approximately inner and satisfy $\af^n = \id_A$.
Suppose that the action of $\Zqn$ generated by $\af$ has the
tracial Rokhlin property.
Assume that $A$ has cancellation of \pj s,
and that the \pj s in $\CZnAa$
have the weak divisibility property (Definition~\ref{WDivDfn}).
Then the dual action of $\Zqn$ on $\CZnAa$
has the tracial Rokhlin property and its generator
is tracially approximately inner.
\end{thm}

\begin{proof}
Let $B = \CZnAa$, and let $u \in B$ be the
standard unitary (satisfying $u a u^* = \af (a)$ for $a \in A$).
Let ${\widehat{\af}} \in \Aut (B)$ be the automorphism which generates
the dual action.
Thus, with $\om = \exp (2 \pi i / n)$,
we have ${\widehat{\af}} (u) = \om u$ and
${\widehat{\af}} (a) = a$ for $a \in A$.

When we verify the definition of the \aRp\  and
of tracial approximate innerness for ${\widehat{\af}}$,
we may clearly assume that the finite set has the form
$F = F_0 \cup \{ u \}$ for some $\af$-invariant finite subset $F_0$
of the unit ball of $A$.
We now claim that we may take the nonzero positive element $x \in B$
to in $A$.
Theorem~4.2 of \cite{JO} shows that if $A$ has Property~(SP),
then so does $B$.
Moreover, $B$ is simple by Corollary~\ref{CrPrIsSimple},
so Theorem~4.2 of \cite{JO} also shows that if $B$ has Property~(SP),
then so does $\Cs{n}{B}{{\widehat{\af}}} \cong M_n (A)$.
Property~(SP) obviously passes to \hsa s, so we see that
$A$ has Property~(SP) \ifo\  $B$ has Property~(SP).
If $A$ and $B$ do not have Property~(SP), then by Lemma~\ref{RImpSP}
and Remark~\ref{InnAndTAInn} we may as well take $x = 0$.
If $A$ and $B$ have Property~(SP),
and if $x \in B$ is a nonzero positive element,
then we can choose a nonzero \pj\  $q_0 \in {\overline{x B x}}$,
and use Theorem~4.2 of \cite{JO} to find a nonzero \pj\  $p \in A$
such that $p \precsim q_0$ in $B$.
It clearly suffices to use $p$ in place of $x$.
This proves the claim.

Accordingly, let $F_0$ be an $\af$-invariant finite subset
of the unit ball of $A$ and set $F = F_0 \cup \{ u \}$,
let $\ep > 0$, let $N \in \N$, and
let $x \in A$ be a nonzero positive element.
\Wolog\  $\ep < 1$.

Choose a \pj\  $e \in A$ and unitaries $v_1, \, v_2 \in e A e$
following Proposition~\ref{FinOrdTAInn2},
with $F_0$ in place of $F$,
with $(n + 2)^{-1} \ep$ in place of $\ep$,
with $N (n + 1)$ in place of $N$,
and with $x$ as given.
In $B$, since
\[
u^n = 1, \,\,\,\,\,\, v_1^n = e, \,\,\,\,\,\, u e u^* = \af (e) = e,
\andeqn u v_1 u^* = \af (v_1) = v_1,
\]
we have $(u^* v_1)^n = e$.
Moreover, for $a \in F_0$ we have
\[
\| v_1 (e a e) v_1^* - \af (e a e) \| < \frac{\ep}{n + 2} \andeqn
u (e a e) u^* = \af (e a e),
\]
so
\[
\| (u^* v_1) (e a e) (u^* v_1)^* - e a e \| < \frac{\ep}{n + 2}.
\]
Therefore, for $0 \leq j \leq n - 1$ we have
\[
\| (u^* v_1)^j (e a e) - (e a e) (u^* v_1)^j \|
  < j \left( \frac{\ep}{n + 2} \right)
  \leq \frac{n \ep}{n + 2}.
\]
Define \pj s $e_j \in B$ by
\[
e_j = \frac{1}{n} \sum_{k = 0}^{n - 1} (\om^{-j} u^* v_1)^k
\]
for $0 \leq j \leq n - 1$.
Then $\sum_{j = 0}^{n - 1} e_j = e$.
Moreover,
\[
{\widehat{\af}} (e_j)
  = \frac{1}{n} \sum_{k = 0}^{n - 1}
             ( \om^{-j} {\overline{\om}} u^* v_1)^k
  = \frac{1}{n} \sum_{k = 0}^{n - 1} (\om^{-(j + 1)} u^* v_1)^k
  = e_{j + 1}.
\]
Also, $u$ commutes with $v_1$, hence with $u^* v_1$,
hence with all $e_j$, and for $a \in F_0$ we have
\begin{align*}
\| a e_j - e_j a \|
& = \| a e e_j - e_j a e \|
  \leq 2 \| e a - a e \| + \| (e a e) e_j - e_j (e a e) \|  \\
& \leq 2 \| e a - a e \|
     + \frac{1}{n} \sum_{k = 0}^{n - 1}
              \| (e a e) (u^* v_1)^k - (u^* v_1)^k (e a e) \|
               \\
& < (2 + n) \left( \frac{\ep}{n + 2} \right) = \ep.
\end{align*}
We have verified Conditions~(1) and~(2) of Definition~\ref{ARPDfn}.
For Condition~(3), we observe that $1 - e$ is \mvnt\  to a \pj\  in
${\overline{x A x}}$.
Condition~(4) will be verified below.

We now set $v = u^* v_2$ and
verify the first four conditions of Definition~\ref{TAInnDfn},
using the same \pj\  $e$ as before.
The first two and the fourth have already been done.
For the third, first consider $a \in F_0$.
We have
\[
\| (u^* v_2) (e a e) (u^* v_2)^* - e a e \| < \frac{\ep}{n + 2}
\]
for the same reason that this holds for $v_1$ in place of $v_2$.
So
\[
\| v (e a e) v^* - {\widehat{\af}} (e a e) \|
   = \| v (e a e) v^* - e a e \|
   < \frac{\ep}{n + 2} \leq \ep.
\]
Furthermore,
\[
v u v^* = u^* v_2 u v_2^* u = \af^{-1} (v_2) v_2^* u = \om v_2 v_2^* u
  = {\widehat{\af}} (u).
\]

It remains to verify the last condition of both definitions.
We have
\[
v v_1 v^* = u^* v_2 v_1 v_2^* u = \af^{-1} (v_2 v_1 v_2^*)
 = (\om v_2) v_1 (\om^{-1} v_2^*) = v_2 v_1 v_2^*.
\]
Combining this with the computation of $v u v^*$ and the
estimate $\| v_2 v_1 v_2^* - v_1 \| < (n + 2)^{-1} \ep$, we get
\[
\| v (u^* v_1) v^* - \om^{-1} u^* v_1 \| < \frac{\ep}{n + 2}.
\]
By induction,
\[
\| v (u^* v_1)^k v^* - (\om^{-1} u^* v_1)^k \| < \frac{k \ep}{n + 2}
\]
for $k \geq 1$.
Therefore
\begin{align*}
\| v e_j v^* - e_{j + 1} \|
& \leq \frac{1}{n} \sum_{k = 0}^{n - 1}
    \| v ( \om^{-j} u^* v_1)^k v^* - (\om^{-(j + 1)} u^* v_1)^k \| \\
& < \frac{1}{n} \sum_{k = 0}^{n - 1} \frac{k \ep}{n + 2} \leq \ep.
\end{align*}
Since $\ep < 1$,
it follows that $e_j$ is unitarily equivalent to $e_{j + 1}$.

By construction there are $N (n + 1)$ \mops\  %
\[
g_1, g_2, \dots, g_{N (n + 1)} \leq e,
\]
each of which is \mvnt\  to
$1 - e$.
For $0 \leq m \leq n$ set
\[
h_m = \sum_{k = 1}^N g_{m N + k},
\]
which are $n + 1$ \mops\  in $e A e$, each of which is \mvnt\  to $h_0$.
The weak divisibility property for \pj s in $B$ implies that
$h_0$ is \mvnt\  in $B$ to a sub\pj\  of $e_0$.
Thus there are $N$ \mops\  in $e_0 B e_0$,
each of which is \mvnt\  to $1 - e$.
Since $e_0, \, e_1, \, \dots, \, e_{n - 1}$ are all
unitarily equivalent in $B$, the same is true for $e_j B e_j$
for $0 \leq j \leq n - 1$.
This is Condition~(4) of Definition~\ref{ARPDfn}.
Since $e_0 \leq e$,
we also have Condition~(5) of Definition~\ref{TAInnDfn}.
\end{proof}

\begin{cor}\label{Iff}
Let $A$ be a \suca\  and let $\af \in \Aut (A)$
be tracially approximately inner and satisfy $\af^n = \id_A$.
Suppose that the action of $\Zqn$ generated by $\af$ has the
tracial Rokhlin property.
Assume that $A$ has cancellation of \pj s.
Then $A$ has tracial rank zero \ifo\  $\CZnAa$
has tracial rank zero.
\end{cor}

\begin{proof}
If $A$ has tracial rank zero, the conclusion is
Theorem~\ref{RokhTAF}.
So suppose that $\CZnAa$ has tracial rank zero.
Then the \pj s in $\CZnAa$ have
the weak divisibility property by Lemma~\ref{GetWDiv}.
Theorem~\ref{DualHasRokh} applies,
so that the dual action of $\Zqn$
on $\CZnAa$ has the tracial Rokhlin property.
Therefore Theorem~\ref{RokhTAF} implies that the crossed
product by the dual action has tracial rank zero.
Since this crossed product is isomorphic to $M_n \otimes A$,
it follows from
Theorem~3.12(1) of~\cite{LnTAF} that $A$ has tracial rank zero.
\end{proof}

Now we turn to the fixed point algebra.

\begin{prp}\label{FPIsSimple}
Let $A$ be a \ca, let $G$ be a compact group, and let
$\af \colon G \to \Aut (A)$ be a \ct\  action of $G$ on $A$.
Suppose $C^* (G, A, \af)$ is simple.
Then the fixed point algebra $A^{\af}$ is simple,
is isomorphic to a full \hsa\  of $C^* (G, A, \af)$, and is
strongly Morita equivalent to $C^* (G, A, \af)$.
\end{prp}

\begin{proof}
See the Proposition, Corollary, and proof of the Corollary in~\cite{Rs}.
\end{proof}

\begin{cor}\label{FPIsTAF}
Let $A$ be a \suca, and let $\af \in \Aut (A)$ be an automorphism which
satisfies $\af^n = \id_A$ and such that the action of $\Zqn$
generated by $\af$ has the tracial Rokhlin property.
Suppose that $A$ has tracial rank zero.
Then the fixed point algebra $A^{\af}$ has tracial rank zero.
\end{cor}

\begin{proof}
Theorem~\ref{RokhTAF} implies that $\CZnAa$
has tracial rank zero.
This algebra is simple by Proposition~\ref{CrPrIsSimple}.
It therefore follows from Proposition~\ref{FPIsSimple} that
$A^{\af}$ is isomorphic to a \hsa\  in $\CZnAa$.
So Theorem~3.12(1) of~\cite{LnTAF} implies that
$A^{\af}$ has tracial rank zero.
\end{proof}

\begin{cor}\label{WhenFPIsTAF}
Let $A$ be a \suca, and let $\af \in \Aut (A)$ be a
tracially approximately inner automorphism which
satisfies $\af^n = \id_A$ and such that the action of $\Zqn$
generated by $\af$ has the tracial Rokhlin property.
Assume that $A$ has cancellation of \pj s, and that
the fixed point algebra $A^{\af}$ has tracial rank zero.
Then $A$ has tracial rank zero.
\end{cor}

\begin{proof}
The crossed product $\CZnAa$ is simple
by Proposition~\ref{CrPrIsSimple}.
It therefore follows from Proposition~\ref{FPIsSimple} that
$A^{\af}$ and $\CZnAa$ are strongly Morita
equivalent.
Since everything is unital, $\CZnAa$ is
isomorphic to a \hsa\  in some matrix algebra over
$A^{\af}$.
Therefore Theorems~3.10 and~3.12(1) of~\cite{LnTAF} imply that
$\CZnAa$ has tracial rank zero.
Now Corollary~\ref{Iff} shows that $A$ has tracial rank zero.
\end{proof}

\section{Higher dimensional noncommutative toruses}\label{Sec:NCT}

\indent
In this section and the next two, we prove that every
simple higher dimensional noncommutative torus is an AT algebra.

The first result in this direction was the
Elliott-Evans Theorem~\cite{EE} for the ordinary irrational
rotation algebras, which we use here as the initial step of an
induction argument.
Without giving a complete list of later work, we mention four
highlights.
All simple three dimensional noncommutative toruses were shown to be
AT algebras in~\cite{LQ1}.
In arbitrary dimension,
``most'' simple higher dimensional noncommutative toruses were
shown to be AT algebras in~\cite{Bc}.
Corollary~6.6 of~\cite{Ks4}
gives this result in all cases in which,
in the skew symmetric matrix giving the commutation relations,
the entries above the diagonal are rationally independent,
as well as some others.
Theorem~3.14 of~\cite{Ln14} shows that the crossed product of
$(S^1)^{d}$ by a minimal rotation is an AT algebra; in this case,
most of the entries of the relevant skew symmetric matrix are zero.

The proof of Corollary~6.6 of~\cite{Ks4} is an induction argument:
if, when the noncommutative torus is written as a successive crossed
product by actions of $\Z$, all the intermediate crossed products
are simple, then the main results of~\cite{Ks4} reduce the
problem to the Elliott-Evans Theorem.
One has some choice here: different choices of the commutation
relations may well give the same \ca.
As a very simple example, one might simply write the generators
in a different order.
Unfortunately, it seems not to be possible in general to choose
commutation relations to give the same algebra, or even a
Morita equivalent algebra (see~\cite{RS}), and in such a way
that the method of Corollary~6.6 of~\cite{Ks4} applies, or even
in such a way as to get a tensor product of algebras to which
this method applies.
However, if one allows one more kind of modification, namely the
replacements of unitary generators by integer powers of themselves,
then it is always possible to replace a noncommutative torus by a
tensor product of algebras covered by Corollary~6.6 of~\cite{Ks4}.
The new algebra isn't isomorphic, or even Morita equivalent, to
the original.
But if one replaces only one generator, the new algebra is the fixed
point algebra of a tracially approximately inner action of a finite
cyclic group which has the \aRp.
This operation thus preserves tracial rank zero.
Because $K_0$ and $K_1$ are torsion free, the classification
theorem for simple nuclear \ca s with tracial rank zero,
Theorem~5.2 of~\cite{Ln15}, shows
that this operation preserves the property of being an AT algebra.

In this section we present, in a form convenient for our purposes,
some standard facts about higher dimensional noncommutative toruses.
In Section~\ref{Sec:ARokhOnNCT},
we prove that the relevant actions have the \aRp,
and in Section~\ref{Sec:NCTInd} we show how to combine this
result with~\cite{EE}, \cite{Ks4}, and~\cite{Ln15} to prove that
all simple higher dimensional noncommutative toruses are AT algebras.

\begin{ntn}\label{StdNtn}
Let $\te$ be a skew symmetric real $d \times d$ matrix.
The {\emph{noncommutative torus}} $A_{\te}$ is by definition~\cite{Rf2}
the universal \ca\  generated by unitaries $u_1, u_2, \dots, u_d$
subject to the relations
\[
u_k u_j = \exp (2 \pi i \te_{j, k} ) u_j u_k
\]
for $1 \leq j, k \leq d$.
(Of course, if all $\te_{j, k}$ are integers, it is not really
noncommutative.)
\end{ntn}

Some authors use $\te_{k, j}$ in the commutation relation instead.
See for example~\cite{Ks2}.

\begin{rmk}\label{CoordFree}
We note (see the beginning of Section~4 of~\cite{Rf1} and
the introduction to~\cite{RS}) that $A_{\te}$ is the universal
\ca\  generated by unitaries $u_x$, for $x \in \Z^d$,
subject to the relations
\[
u_y u_x = \exp (\pi i \langle x, \te (y) \rangle ) u_{x + y}
\]
for $x, \, y \in \Z^d$.

It follows that if $B \in {\mathrm{GL}}_d (\Z)$, and if
$B^{\mathrm{t}}$ denotes the transpose of $B$, then
$A_{B^{\mathrm{t}} \te B} \cong A_{\te}$.
That is, $A_{\te}$ is unchanged if $\te$ is rewritten in terms of
some other basis of $\Z^d$.
\end{rmk}

\begin{rmk}\label{AlgFromBichar}
Let $\af$ be a skew symmetric real bicharacter on $\Z^d$, that is,
a $\Z$-bilinear function $\af \colon \Z^d \times \Z^d \to \R$
such that $\af (x, y) = - \af (y, x)$ for all $x, \, y \in \Z^d$.
For any basis $(b_1, b_2, \dots, b_d)$ of $\Z^d$,
there is a unique skew symmetric real $d \times d$ matrix $\te$
such that
\[
\af \left( \ssum{k = 1}{d} x_k b_k, \, \ssum{k = 1}{d} y_k b_k \right)
  = \ssum{j, k = 1}{d}  x_j \te_{j, k} y_k
\]
for all $x, \, y \in \Z^d$.
We define $A_{\af} = A_{\te}$.
Remark~\ref{CoordFree} shows that this \ca\  is independent of the
choice of basis.
\end{rmk}

\begin{rmk}\label{Restrict}
Let $\af$ be a skew symmetric real bicharacter on $\Z^d$, and let
$H \S \Z^d$ be a subgroup.
Then $H \cong \Z^m$ for some $m \leq d$.
By abuse of notation, we write $\af |_H$ for the restriction of
$\af$ to $H \times H \S \Z^d \times \Z^d$.
There is a noncommutative torus $A_{\af |_H}$ by
Remark~\ref{AlgFromBichar}, which does not depend on the choice
of the isomorphism $H \cong \Z^m$.

For a skew symmetric real $d \times d$ matrix $\te$ and a subgroup
$H \S \Z^d$ with a specified ordered basis, we write
$\te |_H$ for the matrix in that basis of the restriction to $H$ of
the real bicharacter $(x, y) \mapsto \langle x, \te y \rangle$.
For subgroups such as $\Z^m \times \{ 0 \}$ or
$\Z^m \times \{ 0 \} \times \Z^l$, we use without comment the
obvious basis.
\end{rmk}

We formalize a remark made in 1.7 of~\cite{El0}, according to which all
noncommutative toruses can be obtained as successive crossed
product by $\Z$.

\begin{lem}\label{ItCrPrd}
Let $\af$ be a skew symmetric real bicharacter on $\Z^d$.
Then there is an automorphism
$\ph$ of $A_{\af |_{\Z^{d - 1} \times \{ 0 \} }}$ which is homotopic
to the identity and such that
\[
A_{\af}
 \cong C^* ( \Z, \, A_{\af |_{\Z^{d - 1} \times \{ 0 \} }}, \, \ph).
\]
\end{lem}

\begin{proof}
Let $\te$ be the matrix of $\af$ in the standard basis.
Let $\bt = \af |_{\Z^{d - 1} \times \{ 0 \} }$.
Then the matrix of $\bt$ is
$( \te_{j, k} )_{1 \leq j, k \leq d - 1}$.
Let $u_1, \, u_2, \, \dots, \, u_{d - 1}$ be the standard
generators of $A_{\bt}$.
Then $\ph$ is determined by
$\ph (u_j) = \exp (2 \pi i \af_{j, d}) u_j$.
It is clear that $\ph$ is homotopic to the identity.
\end{proof}

The following definition is essentially from Section~1.1 of~\cite{Sl}.

\begin{dfn}\label{NondegDfn}
The skew symmetric real $d \times d$ matrix $\te$ is
{\emph{nondegenerate}} if whenever $x \in \Z^d$ satisfies
$\exp (2 \pi i \langle x, \te y \rangle ) = 1$ for all $y \in \Z^d$,
then $x = 0$.
Otherwise, $\te$ is {\emph{degenerate}}.
We similarly refer to degeneracy and nondegeneracy of
a skew symmetric real bicharacter on $\Z^d$.
\end{dfn}

\begin{lem}\label{NondegCond}
Let $\te$ be a skew symmetric real $d \times d$ matrix.
Then $\te$ is degenerate \ifo\  there
exists $x \in \Q^d \SM \{ 0 \}$ such that
$\langle x, \, \te y \rangle \in \Q$ for all $y \in \Q^d$.
\end{lem}

\begin{proof}
If $\te$ is degenerate, choose $w \neq 0$ such that
$\exp (2 \pi i \langle w, \te y \rangle ) = 1$ for all $y \in \Z^d$.
Then $\langle w, \te y \rangle \in \Z$ for all $y \in \Z^d$.
If now $y \in \Q^d$ is arbitrary,
then there exists $m \in \Z \SM \{ 0 \}$ such that $m y \in \Z^d$.
So
\[
\langle w, \te y \rangle
 = \ts{\frac{1}{m}} \langle w, \te (m y) \rangle
 \in \ts{\frac{1}{m}} \Z \S \Q.
\]

Conversely, assume $x \in \Q^d \SM \{ 0 \}$ and
$\langle x, \, \te y \rangle \in \Q$ for all $y \in \Q^d$.
Choose $m \in \Z$ with $m > 0$ such that
$m \langle x, \, \te \dt_k \rangle \in \Z$ for $1 \leq k \leq d$.
Then $m x \neq 0$ and
$\exp (2 \pi i \langle m x, \te y \rangle ) = 1$ for all $y \in \Z^d$.
\end{proof}

\begin{lem}\label{ConjByRat}
Let $\te$ be a skew symmetric real $d \times d$ matrix.
Let $B \in {\mathrm{GL}}_d (\Q)$.
Then $B^{\mathrm{t}} \te B$ is nondegenerate
\ifo\  $\te$ is nondegenerate.
\end{lem}

\begin{proof}
It suffices to prove one direction.
Suppose $\te$ is degenerate.
By Lemma~\ref{NondegCond}, there is
$x \in \Q^d \SM \{ 0 \}$ such that
$\langle x, \, \te y \rangle \in \Q$ for all $y \in \Q^d$.
Then $B^{-1} x \in \Q^d \SM \{ 0 \}$ and
\[
\langle B^{-1} x, \, B^{\mathrm{t}} \te B y \rangle
  = \langle x, \, \te B y \rangle \in \Q
\]
for all $y \in \Q^d$.
So $B^{\mathrm{t}} \te B$ is degenerate.
\end{proof}

The following result is well known.

\begin{thm}\label{Simplicity}
The \ca\  $A_{\te}$ of Notation~\ref{StdNtn} is simple \ifo\  $\te$
is nondegenerate.
Moreover, if $A_{\te}$ is simple it has a unique tracial state.
\end{thm}

\begin{proof}
If $\te$ is nondegenerate, then $A_{\te}$ is simple by Theorem~3.7
of~\cite{Sl}.
(Note the standing assumption of nondegeneracy throughout
Section~3 of~\cite{Sl}.)

When $A_{\te}$ is simple, the proof of Lemma~3.1 of~\cite{Sl}
shows that $A_{\te}$ can have at most one tracial state.
Existence of a tracial state is well known, or can be
obtained from Lemma~\ref{ItCrPrd} by induction on $n$.

If $\te$ is degenerate, then we follow 1.8 of~\cite{El0}.
Choose $n \in \Z^d \SM \{ 0 \}$ such that
$\exp (2 \pi i \langle n, \te y \rangle ) = 1$ for all $y \in \Z^d$.
Then $v = u_1^{n_1} u_2^{n_2} \cdots u_d^{n_d}$
is a nontrivial element of the center of $A_{\te}$,
which is therefore not simple.
\end{proof}

\section{The tracial Rokhlin property and
  higher dimensional noncommutative toruses}\label{Sec:ARokhOnNCT}

\indent
In this section, we prove that if $\te$ is nondegenerate, then the
action of $\Zqn$ which multiplies one of the standard generators
of $A_{\te}$ by a primitive $n$-th root of~$1$ has the \aRp.
We note for comparison the related result in Section~6 of~\cite{Ks2},
that if $\af \in \Aut (A_{\te})$ is of the form
$\af (u_j) = \ld_j u_j$, with $\ld_1, \ld_2, \dots, \ld_n \in S^1$,
and if all positive powers of $\af$ are outer,
then $\af$ has the Rokhlin property.

As is done in the proof in~\cite{BDR}
that $A_{\te}$ has real rank zero,
and analogously to Section~6 of~\cite{Ks2},
we will reduce to a construction in the ordinary
irrational rotation algebras.
We therefore begin with several facts about these algebras.

For reference, and to establish notation,
we state the following theorem.
Its proof is contained in Theorem~1.1 and Proposition~1.3 of~\cite{AP}.
Also see Corollary~3.6 and Definition~3.3 of~\cite{RfF}.
We refer to~\cite{Dx} for information on \ct\  fields of \ca s.
See especially Sections~10.1 and~10.3.
We use $v$ and $w$
for the generators of the irrational rotation algebra,
to avoid confusion with the generators $u_1, u_2, \dots, u_d$
of a higher dimensional noncommutative torus $A_{\te}$.

\begin{thm}\label{FieldOfRotAlg}
For $\et \in \R$ let $A_{\et}$ be the rotation algebra, the
universal \ca\  %
generated by unitaries $v_{\et}$ and $w_{\et}$ satisfying
$w_{\et} v_{\et} = \exp (2 \pi i \et) v_{\et} w_{\et}$.
Let $A$ be the \ca\  of the discrete Heisenberg group,
which is the universal \ca\  generated by unitaries $v, w, z$
subject to the relations
\[
w v = z v w, \,\,\,\,\,\, z v = v z, \andeqn z w = w z.
\]
Then there is a \ct\  field of \ca s over $S^1$ whose fiber
over $\exp (2 \pi i \et)$ is $A_{\et}$, whose \ca\  of \ct\  sections
is $A$, and such that the
evaluation map $\ev_{\et} \colon A \to A_{\et}$ of sections
at $\exp (2 \pi i \et)$ is determined by
\[
\ev_{\et} (v) = v_{\et}, \,\,\,\,\,\, \ev_{\et} (w) = w_{\et},
\andeqn \ev_{\et} (z) = \exp (2 \pi i \et) \cdot 1.
\]
\end{thm}

Since we will only formally deal with one \ct\  field in this section,
the following notation is unambiguous.

\begin{ntn}\label{SectionNot}
For a subset $E \S S^1$, we let $\Gm (E)$ be the set of \ct\  sections
of the \ct\  field of Theorem~\ref{FieldOfRotAlg} over $E$.
(See 10.1.6 of~\cite{Dx}.)

For any such section $a$, we further write
$a (\et)$ for $a (\exp (2 \pi i \et))$.
No confusion should arise.
\end{ntn}

\begin{lem}\label{ContOfTrace}
Let the notation be as in Theorem~\ref{FieldOfRotAlg}
and Notation~\ref{SectionNot}.
Let $\ta_{\et}$ be the standard trace on $A_{\et}$, satisfying
$\ta (1) = 1$ and
$\ta_{\et} (v_{\et}^m w_{\et}^n) = 0$ unless $m = n = 0$.
Let $U \S S^1$ be an open set, and let $a \in \Gm (U)$.
Then $\et \mapsto \ta_{\et} (a (\et))$ is \ct.
\end{lem}

\begin{proof}
We check continuity at $\et_0$.
Choose a \cfn\  $h \colon S^1 \to [0, 1]$ such that $\supp (h) \S U$
and such that $h = 1$ on a \nbhd\  of $\et_0$.
Then it suffices to consider the \ct\  section $h a$ in place of $a$.
Now $h a$ is the restriction to $U$ of a \ct\  section $b$ defined on
all of $S^1$, satisfying $b (\zt) = 0$ for $\zt \not\in U$.
Accordingly, we may restrict to the case $U = S^1$.
Then $a \in A$.

{}From the formulas
\[
\ev_{\et} (v) = v_{\et}, \,\,\,\,\,\, \ev_{\et} (w) = w_{\et},
\andeqn \ev_{\et} (z) = \exp (2 \pi i \et) \cdot 1
\]
and the definition of $\ta_{\et}$,
it is immediate that if $b$ is any (noncommutative) monomial in
$v$, $w$, $z$, and their adjoints,
then $\et \mapsto \ta_{\et} (b (\et))$ is \ct.
Therefore the same holds for any noncommutative polynomial,
and hence for any norm limit of noncommutative polynomials,
including $a$.
\end{proof}

\begin{lem}\label{2D}
Let the notation be as in Theorem~\ref{FieldOfRotAlg}
and Lemma~\ref{ContOfTrace}.
Let $\et \in \R \SM \Q$.
Let $n \in \N$, let $\om = \exp (2 \pi i / n)$, and let
$\af \colon A_{\et} \to A_{\et}$ be the unique automorphism
satisfying $\af (v_{\et}) = \om v_{\et}$ and $\af (w_{\et}) = w_{\et}$.
Then for every $\ep > 0$ there exist
\mops\  $e_0, e_1, \dots, e_{n - 1}$ such that (with $e_n = e_0$)
we have $\af (e_j) = e_{j + 1}$ for $0 \leq j \leq n - 1$, and such that
$1 - n \ta_{\et} (e_0) < \ep$.
\end{lem}

\begin{proof}
Set $\ep_0 = \frac{1}{4 n} \ep$.
Let $f \colon S^1 \to [0, 1]$ be a \cfn\  such that
$\supp (f)$ is contained in the open arc from $1$ to $\om$,
and such that $f (\zt) = 1$ for all $\zt$ in the closed arc
from $\exp (2 \pi i \ep_0)$ to
$\exp \left( 2 \pi i \left[ \frac{1}{n} - \ep_0 \right] \right)$.
Then $f (v_{\et})$ is a positive element of $A_{\et}$ with
$\| f (v_{\et}) \| \leq 1$ and
$\ta_{\et} (f (v_{\et})) \geq \frac{1}{n} - 2 \ep_0$.
Since $A_{\et}$ has real rank zero (see Remark~6 of~\cite{EE},
or Theorem~1.5 of~\cite{BKR}),
there is a \pj\  $e_0$ in the \hsa\  $B$ of $A_{\et}$ generated by
$f (v_{\et})$ such that $\| e_0 f (v_{\et}) - f (v_{\et}) \| < \ep_0$.
Therefore
$\| e_0 f (v_{\et}) e_0 - f (v_{\et}) \| < 2 \ep_0$.
Since $e_0 f (v_{\et}) e_0 \leq e_0$, it follows that
\[
\ta_{\et} (e_0) \geq \ta_{\et} (e_0 f (v_{\et}) e_0)
   > \ta_{\et} (f (v_{\et}) ) - 2 \ep_0
   \geq \ts{ \frac{1}{n} } - 4 \ep_0.
\]

We have $\af^k (f (v_{\et})) \af^l (f (v_{\et})) = 0$ for
$0 \leq k, \, l \leq n - 1$ and $k \neq l$.
Therefore $\af^k (B) \af^l (B) = \{ 0 \}$ for such $k$ and $l$,
whence also $\af^k (e_0) \af^l (e_0) = 0$.
Define $e_k = \af^k (e_0)$ for $0 \leq k \leq n - 1$.
Then $e_0, e_1, \dots, e_{n - 1}$ are \mops\  such that
$\af (e_j) = e_{j + 1}$ for $0 \leq j \leq n - 1$.
Moreover,
\[
1 - n \ta_{\et} (e_0)
  < 1 - n \left( \ts{ \frac{1}{n} } - 4 \ep_0 \right)
  = 4 n \ep_0 = \ep,
\]
as desired.
\end{proof}

We now return to the higher dimensional noncommutative toruses.
The idea is to find an approximately central copy
of an ordinary irrational rotation algebra $A_{\et}$,
such that the restriction to it of our action
is the one in Lemma~\ref{2D}.
Since the \pj s in $A_{\et}$ must be chosen ahead of time,
at least approximately, we must require
that $\et$ be arbitrarily close to some fixed $\et_0$.
Nondegeneracy enters through Lemma~\ref{DensityCond} below.
To obtain the correct restricted action,
we use the condition~(3) in Lemma~\ref{MakeSubalgStep2} below.
From then on, we roughly follow the argument used in~\cite{BKR}
to prove approximate divisibility.
We vary the arrangement slightly to make part of the argument
easily available for later use.

\begin{dfn}\label{HomToInnD}
Let $\te$ be a skew symmetric real $d \times d$ matrix.
Let
\[
n = (n_1, n_2, \dots, n_d) \in \Z^d \andeqn
v = u_1^{n_1} u_2^{n_2} \cdots u_d^{n_d} \in A_{\te}.
\]
We write $\gm_n$ for the inner automorphism $\Ad (v)$ of the
noncommutative torus $A_{\te}$.
We further define a \hm\  $\sm \colon \Z^d \to (S^1)^d$ by
the formula $\sm (n)_j = \exp (2 \pi i (\te n)_j)$ for
$1 \leq j \leq d$.
(Here the expression $\te n$
is the usual action of a $d \times d$ matrix on an element of $\R^d$.)
\end{dfn}

\begin{lem}\label{HomToInnW}
Let $\te$ be a skew symmetric real $d \times d$ matrix.
With $\gm$ and $\sm$ as in Definition~\ref{HomToInnD},
we have $\gm_n (u_j) = \sm (n)_j u_j$ for $n \in \Z^d$ and
$1 \leq j \leq d$.
Moreover, if $m \in \Z^d$, then
\[
\gm_n ( u_1^{m_1} u_2^{m_2} \cdots u_d^{m_d} )
  = \exp ( 2 \pi i \langle m, \te n \rangle )
               u_1^{m_1} u_2^{m_2} \cdots u_d^{m_d}
\]
for all $n \in \Z^d$.
\end{lem}

\begin{proof}
The first formula is the special case of the second obtained
by setting $m = \dt_j$, the $j$-th standard basis vector of $\Z^d$.
By linearity, both formulas will follow if we check the first when
$m = \dt_j$ and $n = \dt_k$.
Since $(\te \dt_k)_j = \te_{j, k}$, this is just
the commutation relation
\[
u_k u_j u_k^* = \exp (2 \pi i \te_{j, k} ) u_j,
\]
which is the same as the one in from Notation~\ref{StdNtn}.
\end{proof}

\begin{lem}\label{DensityCond}
Let $\te$ be a skew symmetric real $d \times d$ matrix.
The \hm\  $\sm \colon \Z^d \to (S^1)^d$ of Definition~\ref{HomToInnD}
has dense range \ifo\  $\te$ is nondegenerate.
\end{lem}

\begin{proof}
Assume $\sm$ does not have dense range.
Let $H = {\overline{\sm (\Z^d) }}$, which is a proper closed subgroup
of $(S^1)^d$.
Choose a nontrivial character $\mu \colon (S^1)^d \to S^1$ whose
kernel contains $H$.
By the identification of the dual group of $(S^1)^d$,
there is $r \in \Z^d \SM \{ 0 \}$ such that
\[
\mu (\zt_1, \zt_2, \dots, \zt_d)
   = \zt_1^{r_1} \zt_2^{r_2} \cdots \zt_d^{r_d}
\]
for all $\zt \in (S^1)^d$.
Because $H \S \Ker (\mu)$, for all $n \in \Z^d$ we have
\begin{align*}
1 & = \mu ( \sm (n))
    = \exp ( 2 \pi i (\te n)_1)^{r_1}
      \exp ( 2 \pi i (\te n)_2)^{r_2} \cdots
      \exp ( 2 \pi i (\te n)_d)^{r_d}   \\
  & = \exp ( 2 \pi i \langle r, \te n \rangle).
\end{align*}
Thus $\te$ is degenerate.

Now suppose that $\te$ is degenerate.
Then we may choose $r \in \Z^d \SM \{ 0 \}$ such that
$\exp ( 2 \pi i \langle r, \te n \rangle) = 1$ for all $n \in \Z^d$.
Reversing the above calculation,
we find that the nontrivial character
\[
\mu (\zt_1, \zt_2, \dots, \zt_d)
   = \zt_1^{r_1} \zt_2^{r_2} \cdots \zt_d^{r_d}
\]
satisfies $\mu (\sm (n)) = 1$ for all $n \in \Z^d$.
Therefore $\sm$ does not have dense range.
\end{proof}

\begin{cor}\label{DensityFinInd}
Let $\te$ be a nondegenerate skew symmetric real $d \times d$ matrix.
Let $G \S \Z^d$ be a subgroup with finite index.
Let $\sm \colon \Z^d \to (S^1)^d$
be the \hm\   of Definition~\ref{HomToInnD}.
Then $\sm (G)$ is dense in $(S^1)^d$.
\end{cor}

\begin{proof}
Let $H = {\overline{\sm (G) }}$.
Let $S$ be a set of coset representatives for $G$ in $\Z^d$.
Then the sets $\sm (m) H$, for $m \in S$, are closed
and are pairwise equal or disjoint.
By Lemma~\ref{DensityCond}, their union is $(S^1)^d$.
Since there are finitely many of them, and since $(S^1)^d$ is
connected, it follows that all are equal to $(S^1)^d$.
\end{proof}

\begin{cor}\label{GaugeAppInn}
Let $\te$ be a nondegenerate skew symmetric real $d \times d$ matrix.
Let $\zt_1, \zt_2, \dots, \zt_d \in S^1$.
Let $\af \in A_{\te}$ be the automorphism determined by
$\af (u_j) = \zt_j u_j$ for $1 \leq j \leq d$.
Then $\af$ is approximately inner.
\end{cor}

\begin{proof}
It suffices to find, for all $\ep > 0$,
a unitary $v \in A_{\te}$ such that
$\| \af (u_j) - v u_j v^* \| < \ep$ for $1 \leq j \leq d$.
Choose $\dt > 0$ small enough that if
$(\om_1, \om_2, \dots, \om_d) \in (S^1)^d$ satisfies
\[
d ( (\om_1, \om_2, \dots, \om_d),
   \, (\zt_1, \zt_2, \dots, \zt_d) ) < \dt,
\]
then $| \om_j - \zt_j | < \ep$ for $1 \leq j \leq d$.
Then use Lemma~\ref{DensityCond} to choose $n \in \Z^d$ such that
$d ( \sm (n), \, (\zt_1, \zt_2, \dots, \zt_d) ) < \dt$.
Take $v = u_1^{n_1} u_2^{n_2} \cdots u_d^{n_d}$
and use Lemma~\ref{HomToInnW}.
\end{proof}

The following corollary seems to be of interest,
but will not be used here.

\begin{cor}\label{AInnByFixedPt}
Let $\te$ be a nondegenerate skew symmetric real $d \times d$ matrix.
Let $n \in \N$ and let $\om = \exp ( 2 \pi i / n)$.
Let $1 \leq k \leq d$, and
let $\af \in A_{\te}$ be the automorphism determined by
$\af (u_k) = \om u_k$ and $\af (u_j) =  u_j$ for $j \neq k$.
Then for every finite set $F \S A_{\te}$ and every $\ep > 0$
there is a unitary $v \in A_{\te}$ such that $\af (v) = v$ and
$\| \af (a) - v a v^* \| < \ep$ for $a \in F$.
\end{cor}

\begin{proof}
\Wolog\  $k = 1$.
The proof is the same as for Corollary~\ref{GaugeAppInn}, but using
the finite index subgroup $n \Z \oplus \Z^{d - 1}$ in place of $\Z^d$.
This substitution is valid by Corollary~\ref{DensityFinInd}.
\end{proof}

\begin{lem}\label{MakeSubalgStep1}
Let $\te$ be a nondegenerate skew symmetric real $d \times d$ matrix.
Let $n, \, N \in \N$, and let $1 \leq k \leq d$.
Then for every $\ep > 0$ there exists
$l = (l_1, l_2, \dots, l_d) \in \Z^d$ such that:
\begin{itemize}
\item[(1)]
$v = u_1^{l_1} u_2^{l_2} \cdots u_d^{l_d}$ satisfies
$\| v u_j - u_j v \| < \ep$ for $1 \leq j \leq d$.
\item[(2)]
$l_k = 1 \pmod n$.
\item[(3)]
There is $j$ such that $| l_j | > N$.
\end{itemize}
\end{lem}

\begin{proof}
\Wolog\  $k = 1$.
Set $\af = \Ad (u_1^*)$.
There are $\zt_1, \zt_2, \dots, \zt_d \in S^1$ such that
$\af (u_j) = \zt_j u_j$ for $1 \leq j \leq d$.
Let $G = n \Z \oplus \Z^{d - 1}$, which is a
finite index subgroup of $\Z^d$.
According to Corollary~\ref{DensityFinInd}, the subgroup
$\sm (G)$ is dense in $(S^1)^d$.
Let
\[
F = \{ l \in \Z^d
 \colon {\mbox{$| l_j | \leq N + 1$ for $1 \leq j \leq d$}} \}.
\]
Since $F$ is finite, $\sm (G \SM F)$ is also dense in $(S^1)^d$.
Choose $\dt > 0$ small enough that if
$(\om_1, \om_2, \dots, \om_d) \in (S^1)^d$ satisfies
\[
d ( (\om_1, \om_2, \dots, \om_d),
   \, (\zt_1, \zt_2, \dots, \zt_d) ) < \dt,
\]
then $| \om_j - \zt_j | < \ep$ for $1 \leq j \leq d$.
Then use density of $\sm (G \SM F)$ to choose $r \in G \SM F$ such that
$d ( \sm (r), \, (\zt_1, \zt_2, \dots, \zt_d) ) < \dt$.
So with $v_0 = u_1^{r_1} u_2^{r_2} \cdots u_d^{r_d}$,
we get $\| v_0 u_j v_0^* - u_1^* u_j u_1 \| < \ep$ for $1 \leq j \leq d$.
Define
\[
l = (r_1 + 1, \, r_2, \, \dots, \, r_d) \in \Z^d
\andeqn
v = u_1^{l_1} u_2^{l_2} \cdots u_d^{l_d} = u_1 v_0 \in A_{\te}.
\]
Clearly $\| v u_j v^* - u_j \| < \ep$ for $1 \leq j \leq d$.
We have $l_1 = 1 \pmod n$ because $r_1 \in n \Z$.
We have $| l_j | > N$ for some $j$, because
$| r_j | > N + 1$ for some $j$.
\end{proof}

The next lemma is the analog in our context of Lemma~4.6
of~\cite{BKR}.

\begin{lem}\label{MakeSubalgStep2}
Let $\te$ be a nondegenerate skew symmetric real $d \times d$ matrix.
Let $n \in \N$, let $1 \leq k \leq d$, and let $\et_0 \in \R \SM \Q$.
Then for every $\ep > 0$ there exist
\[
l = (l_1, l_2, \dots, l_d) \in \Z^d \andeqn
m = (m_1, m_2, \dots, m_d) \in \Z^d
\]
such that:
\begin{itemize}
\item[(1)]
$v = u_1^{l_1} u_2^{l_2} \cdots u_d^{l_d}$
and $w = u_1^{m_1} u_2^{m_2} \cdots u_d^{m_d}$ satisfy
$\| v u_j - u_j v \| < \ep$ and $\| w u_j - u_j w \| < \ep$
for $1 \leq j \leq d$.
\item[(2)]
There is $\et \in \R \SM \Q$
such that $| \exp ( 2 \pi i \et) - \exp ( 2 \pi i \et_0) | < \ep$
and the unitaries $v$ and $w$ of Part~(1) satisfy
$w v = \exp ( 2 \pi i \et) v w$.
\item[(3)]
$l_k = 1 \pmod n$ and $m_k = 0 \pmod n$.
\end{itemize}
\end{lem}

\begin{proof}
\Wolog\  $k = 1$ and $\et_0 \in \left[- \frac{1}{2}, \frac{1}{2} \right]$.
Choose $N \in \N$ so large that $2 \pi / N < \ep$.
Use Lemma~\ref{MakeSubalgStep1} with $\te$, $n$, and $\ep$ as given,
with $k = 1$, and with this value of $N$,
obtaining
\[
l \in \Z^d
\andeqn
v = u_1^{l_1} u_2^{l_2} \cdots u_d^{l_d}.
\]
Note in particular that $\| v u_j v^* - u_j \| < \ep$
for $1 \leq j \leq d$ and $l_1 = 1 \pmod n$.
Let $s$ be an index such that $| l_s | > N$.

Let
\[
T = \left\{ \et \in \R \colon
 {\mbox{$\left( u_1^{r_1} u_2^{r_2} \cdots u_d^{r_d} \right) v
        \left( u_1^{r_1} u_2^{r_2} \cdots u_d^{r_d} \right)^*
             = \exp ( 2 \pi i \et) v$ for some $r \in \Z^d$}} \right\}.
\]
Then $T$ is a subgroup of $\R$ which is generated by $d + 1$ elements,
namely $1$ and elements
corresponding to letting $r$ run through the standard basis vectors
of $\Z^d$.
So $T \cap \Q$ is also finitely generated, and is therefore discrete.
Since $\et_0 \not\in \Q$, we have
$\dist (\et_0, \, T \cap \Q) > 0$.
Set $\ep_0 = \min ( \ep, \, \dist (\et_0, \, T \cap \Q))$.

Set
\[
M = \sum_{j = 1}^d | l_j | \andeqn
\dt = \min \left( \ts{ \frac{1}{2} } \ep_0, \,  M^{-1} \ep_0 \right).
\]
Let $G$ be the finite index subgroup
$G = n \Z \oplus \Z^{d - 1} \S \Z^d$.
Let
\[
\ld = \left(1, \, \dots, \, 1, \, \exp ( 2 \pi i \et_0 / l_s),
 \, 1, \, \dots, \, 1 \right) \in (S^1)^d,
\]
where $\exp ( 2 \pi i \et_0 / l_s)$ is in position $s$.
Use Corollary~\ref{DensityFinInd} and Lemma~\ref{HomToInnW}
to choose $m \in G$ such that
$\sm (m)$, as in Definition~\ref{HomToInnD}, is so close to $\ld$
that $w = u_1^{m_1} u_2^{m_2} \cdots u_d^{m_d}$ satisfies
$\| w u_j w^* - u_j \| < \dt$ for $j \neq s$, and
$\| w u_s w^* - \exp ( 2 \pi i \et_0 / l_s) u_s \| < \dt$.

Since $\dt \leq \ep$, it is clear that
$\| w u_j w^* - u_j \| < \ep$ for $j \neq s$.
Also
\[
\| w u_s w^* - u_s \|
 \leq \| w u_s w^* - \exp ( 2 \pi i \et_0 / l_s) u_s \|
   + | \exp ( 2 \pi i \et_0 / l_s) - 1 |.
\]
Using $\dt \leq \frac{1}{2} \ep$,
the first term is less than $\frac{1}{2} \ep$.
The second term satisfies
\[
| \exp ( 2 \pi i \et_0 / l_s) - 1 |
  < 2 \pi \left| \frac{\et_0}{l_s} \right|
  < 2 \pi \left( \frac{1}{2 N} \right)
  \leq {\textstyle{\frac{1}{2}}} \ep.
\]
Therefore $\| w u_j w^* - u_j \| < \ep$ for $j = s$ as well.
This completes the verification of Part~(1) of the conclusion.
Part~(3) holds because $m_1 \in n \Z$ by construction.

It remains to prove Part~(2).
For each $j$ with $1 \leq j \leq d$,
there is $\zt_j \in S^1$ such that $w u_j w^* = \zt_j u_j$.
Then
\[
w v w^* = \zt_1^{l_1} \zt_2^{l_2} \cdots \zt_d^{l_d} v.
\]
Thus $w v = \exp ( 2 \pi i \et) v w$ for some $\et \in \R$.
By construction we have $| \zt_j - 1 | < M^{-1} \ep_0$ for $j \neq s$,
and $| \zt_s - \exp ( 2 \pi i \et_0 / l_s) | < M^{-1} \ep_0$.
It follows that
\begin{align*}
\left| \rsz{ \zt_1^{l_1} \zt_2^{l_2} \cdots \zt_d^{l_d}
             - \exp ( 2 \pi i \et_0 / l_s)^{l_s} } \right|
  & \leq | l_s | \cdot | \zt_s - \exp ( 2 \pi i \et_0 / l_s) |
      + \sum_{j \neq s} | l_j | \cdot | \zt_j - 1 |  \\
  & < \sum_{j = 1}^d | l_j | M^{-1} \ep_0 \leq \ep_0.
\end{align*}
Therefore
\[
\| w v - \exp ( 2 \pi i \et_0) v w \|
  = \left| \rsz{ \zt_1^{l_1} \zt_2^{l_2} \cdots \zt_d^{l_d}
             - \exp ( 2 \pi i \et_0) } \right|
  < \ep_0,
\]
which is the same as
$| \exp ( 2 \pi i \et) - \exp ( 2 \pi i \et_0) | < \ep_0$.
In particular,
$| \exp ( 2 \pi i \et) - \exp ( 2 \pi i \et_0) | < \ep$,
as desired.
Moreover, $\et \in T$ and there is no $\rh \in T \cap \Q$
such that $| \exp ( 2 \pi i \rh) - \exp ( 2 \pi i \et_0) | < \ep_0$,
whence $\et \not\in \Q$.
\end{proof}

The proofs of the next two results together parallel the proof of
Theorem~1.5 of \cite{BKR}.
The first of them says, roughly,
that higher dimensional noncommutative toruses contain approximately
central copies of irrational rotation algebras,
constructed in a special way.
It will be used again later.
Unfortunately, the rotation parameter varies with the degree
of approximation.

\begin{lem}\label{MakeSubalgStep3}
Let $\te$ be a nondegenerate skew symmetric real $d \times d$ matrix,
let $n \in \N$, and let $1 \leq k \leq d$.
Then for every $\et_0 \in \R$,
every open set $U \S S^1$ containing $\exp ( 2 \pi i \et_0)$,
every finite subset $F \S A_{\te}$,
every finite subset $S \S \Gm (U)$
(following Notation~\ref{SectionNot}),
and every $\ep > 0$,
there exist $\et \in \R \SM \Q$ and
\[
l = (l_1, l_2, \dots, l_d) \in \Z^d \andeqn
m = (m_1, m_2, \dots, m_d) \in \Z^d
\]
such that:
\begin{itemize}
\item[(1)]
$| \et - \et_0 | < \ep$ and $\exp ( 2 \pi i \et) \in U$.
\item[(2)]
$x = u_1^{l_1} u_2^{l_2} \cdots u_d^{l_d}$
and $y = u_1^{m_1} u_2^{m_2} \cdots u_d^{m_d}$ satisfy
$y x = \exp ( 2 \pi i \et) x y$.
\item[(3)]
Following the notation of Theorem~\ref{FieldOfRotAlg},
and with $x$ and $y$ as in Part~(2),
let $\ph \colon A_{\et} \to A_{\te}$ be the \hm\  such that
$\ph (v_{\et}) = x$ and $\ph (w_{\et}) = y$.
Then $\| [a, \, \ph (b (\et)) ] \| < \ep$
for all $a \in F$ and all $b \in S$.
\item[(4)]
$l_k = 1 \pmod n$ and $m_k = 0 \pmod n$.
\end{itemize}
\end{lem}

\begin{proof}
Let the notation be as in Theorem~\ref{FieldOfRotAlg}
and Notation~\ref{SectionNot}.

\Wolog\  $\ep < 1$.
Then there is $\ep_0 > 0$ such that whenever $\zt \in S^1$ satisfies
$| \zt - \exp ( 2 \pi i \et_0) | < \ep_0$,
there is a unique $\et \in \R$ such that
$\exp ( 2 \pi i \et) = \zt$ and $| \et - \et_0 | < \ep$.

\Wolog\  $\| a \| \leq 1$ for all $a \in F$.
Replacing $U$ by an open set $V$ with
$\exp ( 2 \pi i \et_0) \in V \S {\overline{V}} \S U$,
we may assume every $b \in S$ is bounded.
Then \wolog\  $\| b (\et) \| \leq 1$ for all $b \in S$ and $\et \in U$.
Write $F = \{ a_1, a_2, \dots, a_s \}$
and $S = \{ b_1, b_2, \dots, b_t \}$.
Choose polynomials $g_1, g_2, \dots, g_t$
in four noncommuting variables such that
\[
\| g_r (v_{\et_0}, \, v_{\et_0}^*, \, w_{\et_0}, \, w_{\et_0}^*)
        - b_r (\et_0) \|
   < \ts{ \frac{1}{7} \ep }
\]
for $1 \leq r \leq t$.
Because the rotation algebras form a \ct\  field over $S^1$
(Theorem~\ref{FieldOfRotAlg}), there is $\dt > 0$ such that
whenever $| \et - \et_0 | < \dt$ we have $\exp ( 2 \pi i \et_0) \in U$,
and
\[
\| g_r (v_{\et}, \, v_{\et}^*, \, w_{\et}, \, w_{\et}^*) - b_r (\et) \|
   < \ts{ \frac{2}{7} \ep }
\]
for $1 \leq r \leq t$.

Choose polynomials $f_1, f_2, \dots, f_t$
in $2 d$ noncommuting variables such that
\[
\| f_r (u_1, \, u_1^*, \, \dots, \, u_d, \, u_d^*) - a_r \|
   < \frac{\ep}{7 (1 + \ep)}
\]
for $1 \leq r \leq s$.
Choose (see Proposition~4.3 of~\cite{BKR}) $\dt_0 > 0$ such that
whenever $D$ is a \ca\  and
\[
c_1, c_2, \dots, c_{2 d}, d_1, d_2, d_3, d_4 \in D
\]
are elements of norm~$1$
which satisfy $\| [c_r, d_j] \| < \dt_0$ for all $j$ and $r$, then
\[
\| [f_r (c_1, c_2, \dots, c_{2 d}), \,
               g_j (d_1, d_2, d_3, d_4)] \|
  < \ts{ \frac{1}{7} \ep }
\]
for $1 \leq r \leq s$ and $1 \leq j \leq t$.

Apply Lemma~\ref{MakeSubalgStep2} with $\te$, $n$, $\et_0$, and $k$
as given,
and with $\min (\ep_0, \dt, \dt_0)$ in place of $\ep$.
We obtain $\et \in \R \SM \Q$ and
\[
l = (l_1, l_2, \dots, l_d) \in \Z^d \andeqn
m = (m_1, m_2, \dots, m_d) \in \Z^d.
\]
Set
\[
x = u_1^{l_1} u_2^{l_2} \cdots u_d^{l_d} \andeqn
y = u_1^{m_1} u_2^{m_2} \cdots u_d^{m_d}.
\]
By the choice of $\ep_0$, we may assume that $| \et - \et_0 | < \ep$,
and by the choice of $\dt$ we have $\exp ( 2 \pi i \et_0) \in U$.
This is Part~(1) of the conclusion.
Parts~(2) and~(4) are immediate.

It remains to prove Part~(3).
Part~(1) of the conclusion of Lemma~\ref{MakeSubalgStep2} and
the choice of $\dt_0$ ensure that
\[
\| [f_r (u_1, \, u_1^*, \, \dots, \, u_d, \, u_d^*), \,
               g_j (x, x^*, y, y^*)] \|
  < \ts{ \frac{1}{7} \ep }
\]
for $1 \leq r \leq s$ and $1 \leq j \leq t$.
{}From the choice of $\dt$, we get
\[
\| g_j (x, x^*, y, y^*) \|
   < \| \ph (b_j (\et)) \| + \ts{ \frac{2}{7} \ep }
   < 1 + \ep
\]
for $1 \leq j \leq t$.
Using the choice of the polynomials $f_r$, we therefore get
\begin{align*}
\| [a_r, \, \ph (b_j (\et)) ] \|
& \leq 2 \| a_r \| \cdot \| \ph (b_j (\et))
            - g_j (x, x^*, y, y^*) \|  \\
& \hspace*{3em} \mbox{}
     + 2 \| a_r - f_r (u_1, \, u_1^*, \, \dots, \, u_d, \, u_d^*) \|
               \cdot \| g_j (x, x^*, y, y^*) \|  \\
& \hspace*{3em} \mbox{}
     + \| [f_r (u_1, \, u_1^*, \, \dots, \, u_d, \, u_d^*), \,
               g_j (x, x^*, y, y^*)] \|  \\
& < 2 \left( \frac{2 \ep}{7} \right)
     + 2 (1 + \ep) \left( \frac{\ep}{7 (1 + \ep)} \right)
     + \frac{\ep}{7} = \ep
\end{align*}
for $1 \leq r \leq s$ and $1 \leq j \leq t$,
as desired.
\end{proof}

\begin{prp}\label{IrratAPR}
Let $\te$ be a nondegenerate skew symmetric real $d \times d$ matrix.
Let $n \in \N$, let $\om = \exp (2 \pi i / n)$,
let $1 \leq k \leq d$, and, following Notation~\ref{StdNtn}, let
$\af \colon A_{\te} \to A_{\te}$ the unique automorphism
satisfying $\af (u_k) = \om u_k$ and $\af (u_r) = u_r$
for $r \neq k$.
Then the action of $\Zqn$
generated by $\af$ has the tracial Rokhlin property.
\end{prp}

\begin{proof}
Let $\ta$ be the unique tracial state on $A_{\te}$
(Theorem~\ref{Simplicity}).
We will show that for every $\ep > 0$ and
every finite subset $F \S A_{\te}$,
there are
\mops\  $e_0, e_1, \dots, e_{n - 1} \in A_{\te}$ such that:
\begin{itemize}
\item[(1)]
$\| \af (e_j) - e_{j + 1} \| < \ep$ for $0 \leq j \leq n - 1$.
\item[(2)]
$\| e_j a - a e_j \| < \ep$
for $0 \leq j \leq n - 1$ and $a \in F$.
\item[(3)]
$1 - n \ta (e_0) < \ep$.
\end{itemize}

We first argue this this is enough to deduce the \aRp.
We must prove Conditions~(3) and~(4) in Definition~\ref{ARPDfn}.
We first
recall that if $p, \, q \in A_{\te}$ are \pj s with $\ta (p) < \ta (q)$,
then $p \precsim q$.
This follows from Theorems~6.1 and~7.1 of~\cite{Rf1}, or from
Theorems~1.4(d) and~1.5 of~\cite{BKR}.
Also, $A_{\te}$ has Property~(SP) by Theorem~1.4(b) of~\cite{BKR}.
If now $N \in \N$ and a nonzero positive element $x \in A_{\te}$ are
given, then we may use Property~(SP)
to find a nonzero \pj\  $p \in {\overline{x A x}}$.
Let $e = \sum_{j = 0}^{n - 1} e_j$.
Assume, as we clearly may, that $\ep < 1$.
Then $\af (e_j) \sim e_{j + 1}$.
So $\ta (1 - e) = 1 - n \ta (e_0) < \ep$,
whence $\ta (e_0) > \frac{1}{n} (1 - \ep)$.
If
\[
\ep <
 \min \left( \frac{1}{2}, \, \frac{1}{2 n N}, \, \ta (p) \right),
\]
then $\ta (1 - e) < \ep$ implies
$\ta (1 - e) < \ta (p)$ and $N \ta (1 - e) < \ta (e_0)$,
so that Conditions~(3) and~(4) of Definition~\ref{ARPDfn}
follow from the comparison results above.

Now we prove Conditions (1)--(3) at the beginning of the proof.
Let the notation be as in Theorem~\ref{FieldOfRotAlg}
and Notation~\ref{SectionNot}.
Let $\ep > 0$.
Choose and fix $\et_0 \in \R \SM \Q$.
Choose $\ep_1 > 0$ such that whenever $a_0, a_1, \dots, a_{n - 1}$
are elements of a unital \ca\  $D$ with
\[
\| a_j a_r - \dt_{j, r} a_j \| < \ep_1 \andeqn
\| a_j^* - a_j \| < \ep_1
\]
for $0 \leq j, \, r \leq n - 1$, then there are \mops\  %
\[
q_0, q_1, \dots, q_{n - 1} \in D
\]
such that
$\| q_j - a_j \| < \frac{1}{3} n^{-1} \ep$ for $0 \leq j \leq n - 1$.
(For example, apply Definition~2.2 and Lemma~2.3 of~\cite{BKR}
with the \fd\  \ca\  $B$ taken to be $\C^{n + 1}$, using in
addition the element $a_n = 1 - \sum_{j = 0}^{n - 1} a_j$.)
Let $p_0, p_1, \dots, p_{n - 1} \in A_{\et_0}$ be the \pj s
$e_0, e_1, \dots, e_{n - 1}$ of Lemma~\ref{2D}
for $\et_0$ in place of $\et$ and $\frac{1}{3} \ep$ in place of $\ep$.
Because the rotation algebras form a \ct\  field over $S^1$
with section algebra $A$ (Theorem~\ref{FieldOfRotAlg}),
we may choose $c_0, c_1, \dots, c_{n - 1} \in A$
such that $\ev_{\et_0} (c_j) = p_j$ for $0 \leq j \leq n - 1$,
and we can furthermore find $\dt_0 > 0$ such that
$| \exp ( 2 \pi i \et) - \exp ( 2 \pi i \et_0) | < \dt_0$ implies
\[
\| \ev_{\et} (c_j) \ev_{\et} (c_r)
   - \dt_{j, r} \ev_{\et} (c_j) \| < \ep_1
\andeqn
\| \ev_{\et} (c_j)^* - \ev_{\et} (c_j) \| < \ep_1
\]
for $0 \leq j, \, r \leq n - 1$.
Let $V \S S^1$ be an open set such that $\exp ( 2 \pi i \et_0) \in V$
and such that $\zt \in {\overline{V}}$ implies
$| \zt - \exp ( 2 \pi i \et_0) | < \dt_0$.
Letting $c_j |_{\overline{V}}$ denote the restriction of $c_j$,
regarded as a section, to ${\overline{V}}$,
we get
\[
\left\| \ts{ \left( c_j |_{\overline{V}} \right) }
               \ts{ \left( c_r |_{\overline{V}} \right) }
        - \dt_{j, r} c_j |_{\overline{V}} \right\|
     < \ep_1 \andeqn
\left\| \ts{ \left( c_j |_{\overline{V}} \right) }^*
                - c_j |_{\overline{V}} \right\| < \ep_1
\]
for $0 \leq j, \, r \leq n - 1$,
so that there are \mops\  %
\[
q_0, q_1, \dots, q_{n - 1} \in \Gm ( {\overline{V}} )
\]
such that
$\left\| q_j - c_j |_{\overline{V}} \right\| < \frac{1}{3} n^{-1} \ep$
for $0 \leq j \leq n - 1$.
Since the restriction map $A = \Gm (S^1) \to \Gm ( {\overline{V}} )$
is surjective,
there exist $b_0, b_1, \dots, b_{n - 1} \in A$
such that $b_j |_{\overline{V}} = q_j$ for $0 \leq j \leq n - 1$.

Let the generators of $A$ be as in Theorem~\ref{FieldOfRotAlg},
and let $\bt \in \Aut (A)$ be the unique automorphism such that
\[
\bt (v) = \om v, \,\,\,\,\,\, \bt (w) = w, \andeqn \bt (z) = z.
\]
Let $\bt_{\et} \in \Aut (A_{\et})$ be defined by
$\bt_{\et} ( v_{\et}) = \om v_{\et}$
and $\bt_{\et} ( w_{\et}) = w_{\et}$.
Then $\ev_{\et} \circ \bt = \bt_{\et} \circ \ev_{\et}$.
Since $\bt$ sends continuous sections to continuous sections,
there is an open set $U_0 \S V$ such that $\et_0 \in U_0$
and if $\et \in U_0$
then for $0 \leq j \leq n - 1$ and with $b_n = b_0$,
\[
\| \bt_{\et} ( b_j (\et) ) - b_{j + 1} (\et) \|
\andeqn
\| \bt_{\et_0} ( b_j (\et_0) ) - b_{j + 1} (\et_0) \|
\]
differ by less than $\frac{1}{3} \ep$.
For such $\et$ we have $b_j (\et) = q_j (\et)$, so,
using $c_j (\et_0) = p_j$ and $\bt_{\et_0} ( p_j) = p_{j + 1}$
at the second last step,
\begin{align*}
\| \bt_{\et} ( q_j (\et) ) - q_{j + 1} (\et) \|
& < \| \bt_{\et_0} ( q_j (\et_0) ) - q_{j + 1} (\et_0) \|
        + \ts{ \frac{1}{3} } \ep     \\
& \hspace*{-2em}
  < \left\| q_j - c_j |_{\overline{V}} \right\|
        + \left\| q_{j + 1} - c_{j + 1} |_{\overline{V}} \right\|
        + \| \bt_{\et_0} ( c_j (\et_0) ) - c_{j + 1} (\et_0) \|
        + \ts{ \frac{1}{3} } \ep     \\
& \hspace*{-2em}
  < \ts{ \frac{1}{3} } n^{-1} \ep + \ts{ \frac{1}{3} } n^{-1} \ep
        + \ts{ \frac{1}{3} } \ep \leq \ep
\end{align*}
for $0 \leq j \leq n - 1$.

Using Lemma~\ref{ContOfTrace},
choose an open set $U \S U_0$ such that $\et_0 \in U$
and if $\et \in U$ then for $0 \leq j \leq n - 1$ we have
$| \ta_{\et} (q_j (\et)) - \ta_{\et_0} (q_j (\et_0)) |
                < \ts{ \frac{1}{3} } n^{-1} \ep$.

Apply Lemma~\ref{MakeSubalgStep3} with $\te$, $n$, $k$, $\et_0$,
$U$, and $F$ as given,
with $\min ( \ep, \dt)$ in place of $\ep$,
and with $S = \{ q_0, \, q_1, \, \dots, \, q_{n - 1} \}$.
We obtain $\et \in (\R \SM \Q) \cap U$ and
\[
l = (l_1, l_2, \dots, l_d) \in \Z^d \andeqn
m = (m_1, m_2, \dots, m_d) \in \Z^d.
\]
Set
\[
x = u_1^{l_1} u_2^{l_2} \cdots u_d^{l_d} \andeqn
y = u_1^{m_1} u_2^{m_2} \cdots u_d^{m_d},
\]
so that $y x = \exp ( 2 \pi i \et) x y$.
Let $\ph \colon A_{\et} \to A_{\te}$ be the \hm\  such that
$\ph (v_{\et}) = x$ and $\ph (w_{\et}) = y$,
and set $e_j = \ph ( q_j (\et))$ for $0 \leq j \leq n - 1$.
We verify Conditions~(1), (2), and~(3) at the beginning of the
proof for this choice of $e_0, e_1, \dots, e_{n - 1}$.

We do Condition~(1).
Because $l_k = 1 \pmod n$ and $m_k = 0 \pmod n$, we have
$\af (x) = \om x$ and $\af (y) = y$.
It follows that $\af \circ \ph = \ph \circ \bt_{\et}$.
Therefore
\[
\| \af (e_j) - e_{j + 1} \|
   \leq \| \bt_{\et} ( q_j (\et) ) - q_{j + 1} (\et) \| < \ep
\]
for $0 \leq j \leq n - 1$, as desired.

Condition~(2) is immediate from Part~(3) of Lemma~\ref{MakeSubalgStep3}.

Finally, we check Condition~(3).
By uniqueness of the tracial states, we have
$\ta \circ \ph = \ta_{\et}$.
Therefore, using the choice of $U$ at the second step and
$\| q_j (\et_0) - p_j \| < \ts{ \frac{1}{3} } n^{-1} \ep$
at the third step, we get
\[
\ta (e_j)
   = \ta_{\et} (q_j (\et))
   > \ta_{\et_0} (q_j (\et_0)) - \ts{ \frac{1}{3} } n^{-1} \ep
   > \ta_{\et_0} (p_j) - \ts{ \frac{2}{3} } n^{-1} \ep.
\]
Therefore
\[
1 - n \ta (e_0)
  < 1 - n \ta (p_0) + \ts{ \frac{2}{3} } \ep
  < \ts{ \frac{1}{3} } \ep + \ts{ \frac{2}{3} } \ep = \ep.
\]
This completes the proof of~(3).
\end{proof}

As Hanfeng Li pointed out,
the following lemma also holds when $\te$ is degenerate.

\begin{lem}\label{FPOfMult}
Let $\te$ be a nondegenerate skew symmetric real $d \times d$ matrix.
Let $n \in \N$, let $\om = \exp (2 \pi i / n)$,
let $1 \leq l \leq d$, and, following Notation~\ref{StdNtn}, let
$\af \colon A_{\te} \to A_{\te}$ the unique automorphism
satisfying $\af (u_l) = \om u_l$ and $\af (u_k) = u_k$
for $k \neq l$.
Let $B \in {\mathrm{GL}}_d (\Q)$ be the matrix
$B = \diag (1, \dots, 1, n, 1, \dots, 1)$, where $n$ is in the
$l$-th position.
Then the fixed point algebra $A_{\te}^{\af}$ is isomorphic to
$A_{B^{\mathrm{t}} \te B}$.
\end{lem}

\begin{proof}
We observe that
\[
(B^{\mathrm{t}} \te B)_{j, k} = \left\{ \begin{array}{ll}
   n \te_{j, k} & \hspace{3em}  {\mbox{$j = l$ or $k = l$}}   \\
   \te_{j, k}   & \hspace{3em}  {\mbox{otherwise}}
    \end{array} \right..
\]
(Note that $(B^{\mathrm{t}} \te B)_{l, l} = \te_{l, l} = 0$.)
Moreover, $B^{\mathrm{t}} \te B$ is nondegenerate by
Lemma~\ref{ConjByRat}.
Therefore
\[
D = C^* (u_1, \, \dots, \, u_{l - 1}, \, u_l^n,
           \, u_{l + 1}, \, \dots, \, u_d)
  \S A_{\te}
\]
is isomorphic to $A_{B^{\mathrm{t}} \te B}$.

We claim that $A_{\te}^{\af} = D$.
That $D \S A_{\te}^{\af}$ is clear.
For the reverse inclusion, define $E \colon A_{\te} \to A_{\te}^{\af}$
by $E (a) = \frac{1}{n} \sum_{j = 0}^{n - 1} \af^j (a)$.
Then $E$ is a surjective continuous linear map, so it suffices to
show that $E ( u_1^{m_1} u_2^{m_2} \cdots u_d^{m_d} ) \in D$
for all $m = (m_1, m_2, \dots, m_d) \in \Z^d$.
If $m_l$ is divisible by $n$ then
$u_1^{m_1} u_2^{m_2} \cdots u_d^{m_d}$
is a fixed point of $E$ and is in $D$,
and otherwise $E ( u_1^{m_1} u_2^{m_2} \cdots u_d^{m_d} ) = 0 \in D$.
\end{proof}

\begin{cor}\label{IrrFPIsTAF}
Let $\te$ be a nondegenerate skew symmetric real $d \times d$ matrix.
Let $n \in \N$, let $1 \leq l \leq d$, and let
\[
B = \diag (1, \dots, 1, n, 1, \dots, 1) \in {\mathrm{GL}}_d (\Q),
\]
where $n$ is in the $l$-th position.
Then $A_{B^{\mathrm{t}} \te B}$ has tracial rank zero
\ifo\  $A_{\te}$ has tracial rank zero.
\end{cor}

\begin{proof}
The \ca\  $A_{\te}$ has cancellation of \pj s, by
Theorems~6.1 and~7.1 of~\cite{Rf1}, or by
\cite{Hg} and Theorems~1.4(d) and~1.5 of~\cite{BKR}.
The corollary therefore follows from Lemma~\ref{FPOfMult},
Proposition~\ref{IrratAPR}, Corollary~\ref{GaugeAppInn},
Corollary~\ref{FPIsTAF}, and Corollary~\ref{WhenFPIsTAF}.
\end{proof}

\section{Direct limit decomposition for simple noncommutative
            toruses}\label{Sec:NCTInd}

\indent
In this section, we use the results of the previous two sections
to prove that every simple higher dimensional noncommutative torus
is an AT algebra.
(See the introduction to Section~\ref{Sec:NCT} for a discussion
of previous work on this problem.)

The following result is essentially Corollary~6.6 of~\cite{Ks4}.

\begin{prp}\label{SingleStep}
Let $\af$ be a nondegenerate skew symmetric real bicharacter on $\Z^n$.
Suppose that $A_{\af |_{\Z^{n - 1} \times \{ 0 \} }}$ is a
simple AT algebra with real rank zero.
Then $A_{\af}$ is a simple AT algebra with real rank zero.
\end{prp}

\begin{proof}
Let $\bt = \af |_{\Z^{n - 1} \times \{ 0 \} }$.
We note that
\[
K_0 (A_{\bt}) \cong K_1 (A_{\bt}) \cong \Z^{2^{n - 1}}
\]
by Lemma~\ref{ItCrPrd} and by repeated application of the
Pimsner-Voiculescu exact sequence~\cite{PV}.
In particular, both groups are finitely generated.
Further write $A_{\af} = C^* ( \Z, A_{\bt}, \ph)$ as in
Lemma~\ref{ItCrPrd}, with $\ph$ homotopic to the identity.
Thus, in the notation of~\cite{Ks4}
(see the introduction to~\cite{Ks4}),
$\ph \in {\mathrm{HInn}} (A_{\bt})$.
So the proof of Corollary~6.5 of~\cite{Ks4} shows that the
hypotheses of Theorem~6.4 of~\cite{Ks4} hold.
We know from Lemma~\ref{Simplicity} that
$A_{\af} = C^* ( \Z, A_{\bt}, \ph)$ has a unique tracial state.
Therefore Theorem~6.4 of~\cite{Ks4} implies that
$A_{\af} = C^* ( \Z, A_{\bt}, \ph)$
is a simple AT algebra with real rank zero.
\end{proof}

\begin{lem}\label{GenOfGLQ}
The group ${\mathrm{GL}}_d (\Q)$ is generated as a group by
${\mathrm{GL}}_d (\Z)$ and all matrices of the form
$\diag (1, \dots, 1, n, 1, \dots, 1)$, where $n \in \N$ is nonzero
and is in an arbitrary position.
\end{lem}

\begin{proof}
Let $G$ be the subgroup of ${\mathrm{GL}}_d (\Q)$ generated
by ${\mathrm{GL}}_d (\Z)$ and the matrices
$\diag (1, \dots, 1, n, 1, \dots, 1)$.
It suffices to show that $G$ contains all of the following three
kinds of elementary matrices:
\[
E_j^{(1)} (r) = \diag (1, \dots, 1, r, 1, \dots, 1),
\]
where $r \in \Q \SM \{ 0 \}$ and is the $j$-th diagonal entry in the
matrix;
the transposition matrix $E_{j, k}^{(2)}$, for $1 \leq j < k \leq d$,
which acts on the standard basis vectors by
\[
E_{j, k}^{(2)} (\dt_l) = \left\{ \begin{array}{ll}
     \dt_l   & \hspace{3em}  l \neq j, \, k \\
     \dt_k   & \hspace{3em}  l = j  \\
     \dt_j   & \hspace{3em}  l = k
    \end{array} \right.;
\]
and the matrix $E_{j, k}^{(3)} (r)$
for $1 \leq j, \, k \leq n$ with $j \neq k$ and $r \in \Q$,
given by
\[
E_{j, k}^{(3)} (r) (\dt_l) = \left\{ \begin{array}{ll}
     \dt_l             & \hspace{3em}  l \neq k \\
     \dt_k + r \dt_j   & \hspace{3em}  l = k
    \end{array} \right..
\]

If $r = (-1)^m p / q$ with $m = 0$ or $m = 1$ and with $p$ and $q$
positive integers, then
\[
E_j^{(1)} (r) = E_j^{(1)} ((-1)^m) E_j^{(1)} (p) [E_j^{(1)} (q)]^{-1},
\]
where the first factor is in ${\mathrm{GL}}_d (\Z)$ and
$E_j^{(1)} (p)$ and $E_j^{(1)} (q)$ are among the other generators
of $G$.
The matrix $E_{j, k}^{(2)}$ is already in ${\mathrm{GL}}_d (\Z)$.
For $E_{j, k}^{(3)} (r)$, we may conjugate by a permutation matrix,
which is in ${\mathrm{GL}}_d (\Z)$, and split off as a direct summand
a $(d - 2) \times (d - 2)$ identity matrix, and thus reduce to the
case $d = 2$, $j = 1$, and $k = 2$.
Write $r = p / q$ with $p \in \Z$ and $q \in \N$.
Then the factorization
\[
E_{1, 2}^{(3)} (r)
 = \left( \begin{array}{cc} 1 & p/q \\ 0 & 1 \end{array} \right)
 = \left( \begin{array}{cc} q & 0 \\ 0 & 1 \end{array} \right)^{-1}
   \left( \begin{array}{cc} 1 & p \\ 0 & 1 \end{array} \right)
   \left( \begin{array}{cc} q & 0 \\ 0 & 1 \end{array} \right)
\]
shows that $E_{1, 2}^{(3)} (r) \in G$.
\end{proof}

\begin{cor}\label{TAFConjByRat}
Let $\te$ be a nondegenerate skew symmetric real $d \times d$ matrix.
Let $B \in {\mathrm{GL}}_d (\Q)$.
Then $A_{B^{\mathrm{t}} \te B}$ has tracial rank zero
\ifo\  $A_{\te}$ has tracial rank zero.
\end{cor}

\begin{proof}
Combine Lemma~\ref{GenOfGLQ}, Remark~\ref{CoordFree}, and
Corollary~\ref{IrrFPIsTAF}.
\end{proof}

\begin{lem}\label{StrOfNondeg}
Let $\te$ be a nondegenerate skew symmetric real $d \times d$ matrix,
with $d > 2$.
Suppose that there is no subgroup $H$ of $\Z^d$ of rank $d - 1$
such that $\te |_H$
(in the sense of Remark~\ref{Restrict}) is nondegenerate.
Let $r < d - 1$ be the maximal rank of a proper subgroup $H$ of $\Z^d$
such that $\te |_H$ is nondegenerate.
Then there exists $B \in {\mathrm{GL}}_d (\Q)$ such that
$B^{\mathrm{t}} \te B$ has the block form
\[
B^{\mathrm{t}} \te B
   = \left( \begin{array}{cc} \rh_{1, 1} & \rh_{1, 2} \\
           (\rh_{1, 2})^{\mathrm{t}} & \rh_{2, 2} \end{array} \right),
\]
with $\rh_{1, 1}$ and $\rh_{2, 2}$
nondegenerate skew symmetric real $r \times r$ and
$(d - r) \times (d - r)$ matrices, and where all the entries of
$\rh_{1, 2}$ are in $\Z$.
\end{lem}

\begin{proof}
Let $H \S \Z^d$ be a subgroup of rank $r$ such that
$\te |_H$ is nondegenerate.
Since $\te$ is not rational, clearly $r \geq 2$.
Let $(v_1, v_2, \dots, v_r)$ be a basis for $H$ over $\Z$.
Choose $v_{r + 1}, \, \dots, \, v_d \in \Z^d$ such that
$(v_1, v_2, \dots, v_d)$ is a basis for $\Q^d$ over $\Q$.
For $r + 1 \leq k \leq d$, by hypothesis
$\te |_{H + \Z v_k}$ is degenerate.
By Lemma~\ref{NondegCond},
there exists $x_k \in \spn_{\Q} (H \cup \{ v_k \} ) \SM \{ 0 \}$
such that $\langle x_k, \, \te y \rangle \in \Q$
for all $y \in \spn_{\Q} (H \cup \{ v_k \} )$.
Since $\te |_H$ is nondegenerate, we have $x_k \not\in \spn_{\Q} (H)$.
Therefore $v_k \in \spn_{\Q} (H \cup \{ x_k \} )$.
It follows that
\[
v_{r + 1}, \, \dots, \, v_d \in \spn_{\Q}
 (v_1, \, v_2, \, \dots, \, v_r, \, x_{r + 1}, \, \dots, \, x_d),
\]
so that
$(v_1, \, v_2, \, \dots, \, v_r, \, x_{r + 1}, \, \dots, \, x_d)$
is a basis for $\Q^d$.
By construction, we have
$\langle x_k, \te v_l \rangle \in \Q$ for $1 \leq l \leq r$ and
$r + 1 \leq k \leq d$.
Choose $N \in \Z \SM \{ 0 \}$ such that
$N \langle x_k, \te v_l \rangle \in \Z$ for $1 \leq l \leq r$ and
$r + 1 \leq k \leq d$.

Let $B \in {\mathrm{GL}}_d (\Q)$ be the matrix whose
action on the standard basis vectors is
\[
B \dt_k = \left\{ \begin{array}{ll}
     v_k   & \hspace{3em}  1 \leq k \leq r  \\
    N x_k  & \hspace{3em}  r + 1 \leq k \leq d
    \end{array} \right..
\]
Then for $1 \leq l \leq r$ and $r + 1 \leq k \leq d$, we have
\[
\langle \dt_k, B^{\mathrm{t}} \te B \dt_l \rangle
   = N \langle x_k, \te v_l \rangle \in \Z.
\]
Since $B^{\mathrm{t}} \te B$ is skew symmetric, this shows that it
has a block decomposition of the required form.
It is immediate to check that the two diagonal blocks must be
nondegenerate, since otherwise $B^{\mathrm{t}} \te B$ would be
degenerate, contradicting Lemma~\ref{ConjByRat}.
\end{proof}

We will use below, and also on several later occasions, the
following consequence of H.\  Lin's classification theorem~\cite{Ln15}.
An AH algebra is a direct limit of finite direct sums of corners
of homogeneous \ca s whose primitive ideal spaces are
finite CW~complexes.
See, for example, the statement of Theorem~4.6 of~\cite{EG2},
except that we omit the restrictions there on the type of
CW~complexes which may appear;
or see~2.5 of~\cite{LnTAF}.
Also,
see Definition~\ref{DL1} and Remark~\ref{DL2} below for a careful
statement of what it means to satisfy the \uct.

\begin{lem}\label{ConseqOfClass}
Let $A$ be a simple infinite dimensional
separable unital nuclear \ca\  %
with tracial rank zero and which \suct.
Then $A$ is a simple AH algebra with real rank zero
and no dimension growth.
If $K_* (A)$ is torsion free, then $A$ is an AT algebra.
If, in addition, $K_1 (A) = 0$, then $A$ is an AF algebra.
\end{lem}

\begin{proof}
Theorems~6.11 and~6.13 of~\cite{LnTTR} show that
$K_0 (A)$ is weakly unperforated and is a Riesz group.
We now apply Theorem~4.20 of~\cite{EG2} to find a
simple unital AH algebra $B$ with real rank zero and no dimension
growth whose ordered scaled K-theory is the same as that of $A$.
If $K_* (A)$ is torsion free,
we claim that there is a simple AT algebra $B$ with real rank zero
whose ordered scaled K-theory is the same as that of $A$.
To prove this, note that
$K_0 (A)$ can't be $\Z$ because $A$ has real rank zero;
then we apply the proof of Theorem~8.3 of~\cite{Ell2}.
(As noted in the introduction to~\cite{Ell2},
the part of the order involving $K_1$
is irrelevant in the simple case.)
We can certainly take the groups in the direct limit decomposition
to be torsion free,
so that the proof shows that all the algebras in the direct system
constructed there may be taken to have primitive ideal space
the circle or a point.
Then Theorem~4.3 of~\cite{Ell2} shows they may all be taken
to have primitive ideal space the circle.
This gives the required AT algebra $B$.
Finally, if in addition $K_1 (A) = 0$,
following~\cite{Ef} we may find a simple AF algebra $B$
whose ordered scaled K-theory is the same as that of $A$.

Proposition~2.6 of~\cite{LnTAF}
(with ${\mathcal{C}}$ as defined in~2.5 of~\cite{LnTAF})
implies that simple AH algebras with real rank zero
and no dimension growth have tracial rank zero.
In particular, $B$ has tracial rank zero.
So the classification theorem for \ca s with tracial rank zero,
Theorem~5.2 of~\cite{Ln15},
implies that $A \cong B$.
\end{proof}

\begin{thm}\label{NCTIsAT}
Let $\te$ be a nondegenerate skew symmetric real $d \times d$ matrix,
with $d \geq 2$.
Then $A_{\te}$ is a simple AT algebra with real rank zero,
and in particular has tracial rank zero.
\end{thm}

\begin{proof}
We prove this by induction on $d$.
The first part of the conclusion
is true for $d = 2$ by the Elliott-Evans Theorem~\cite{EE},
and tracial rank zero follows from Proposition~2.6 of~\cite{LnTAF}
(with ${\mathcal{C}}$ as defined in~2.5 of~\cite{LnTAF}).
Suppose $d$ is given, and the theorem is known for all
skew symmetric real $k \times k$ matrices with $k < d$.
Let $\te$ be a nondegenerate skew symmetric real $d \times d$ matrix.
There are two cases.

First, suppose that there is a subgroup $H_0$ of $\Z^d$ of rank $d - 1$
such that $\te |_{H_0}$
(in the sense of Remark~\ref{Restrict}) is nondegenerate.
Set
\[
H = \{ x \in \Z^d \colon
 {\mbox{There is $n \in \Z$ such that $n x \in H_0$}} \}.
\]
Then $H$ is also a subgroup of $\Z^d$ of rank $d - 1$, and
$\te |_{H_0}$ is also nondegenerate.
Moreover, $\Z^d / H$ is torsion free and therefore isomorphic to $\Z$,
from which it follows that the quotient map splits.
Thus there is an isomorphism $\Z^d \to \Z^d$ which sends
$H$ isomorphically onto $\Z^{d - 1} \oplus \{ 0 \} \S \Z^d$.
Accordingly, we may assume that $H = \Z^{d - 1} \oplus \{ 0 \}$.
By the induction hypothesis,
$A_{\te |_H}$ is a simple AT algebra with real rank zero.
So Proposition~\ref{SingleStep} implies that
$A_{\te}$ is a simple AT algebra with real rank zero.

Now assume there is no such subgroup $H_0$ of rank $d - 1$.
Let $B$ be as in Lemma~\ref{StrOfNondeg}, with
\[
B^{\mathrm{t}} \te B
   = \left( \begin{array}{cc} \rh_{1, 1} & \rh_{1, 2} \\
             (\rh_{1, 2})^{\mathrm{t}} & \rh_{2, 2} \end{array} \right),
\]
and where in particular all the entries of
$\rh_{1, 2}$ are in $\Z$.
Then $A_{B^{\mathrm{t}} \te B} \cong A_{\rh}$ for
\[
\rh = \left( \begin{array}{cc} \rh_{1, 1} & 0 \\
                                   0   & \rh_{2, 2} \end{array} \right).
\]
Since $\rh$, $\rh_{1, 1}$, and $\rh_{2, 2}$ are all
nondegenerate, it easily follows that
$A_{\rh} \cong A_{ \rh_{1, 1} } \otimes A_{ \rh_{2, 2} }$.
By the induction hypothesis,
both $A_{ \rh_{1, 1} }$ and $A_{ \rh_{2, 2} }$
are simple AT algebras with real rank zero.
Therefore $A_{ \rh_{1, 1} } \otimes A_{ \rh_{2, 2} }$ is a simple
direct limit, with no dimension growth, of homogeneous \ca s.
Since it has a unique tracial state, Theorems~1 and~2 of~\cite{BDR}
imply that $A_{\rh}$ has stable rank one and real rank zero.
This algebra has weakly unperforated
K-theory by Theorem~6.1 of~\cite{Rf1}.
(Actually, this is true for any direct limit of the type at hand.)
It now follows from Theorem~4.6 of~\cite{Ln5} that
$A_{ \rh_{1, 1} } \otimes A_{ \rh_{2, 2} }$ has tracial rank zero.
So Corollary~\ref{TAFConjByRat}
shows that $A_{\te}$ has tracial rank zero.
Using Theorem~1.17 of~\cite{RSUCT} (see the preceding discussion
for the definition of ${\mathcal{N}}$),
it follows from Lemma~\ref{ItCrPrd}
that $A_{\te}$ satisfies the Universal Coefficient Theorem.
Clearly $A_{\te}$ is separable and nuclear.
Since
\[
K_0 (A_{\bt}) \cong K_1 (A_{\bt}) \cong \Z^{2^{n - 1}}
\]
by Lemma~\ref{ItCrPrd} and by repeated application of the
Pimsner-Voiculescu exact sequence~\cite{PV},
Lemma~\ref{ConseqOfClass} implies that $A_{\te}$ is an AT algebra.
\end{proof}

We note that one could use the earlier Theorem~3.11 of~\cite{EG}
to show that $A_{ \rh_{1, 1} } \otimes A_{ \rh_{2, 2} }$ is an
AT algebra with real rank zero, from which it follows that
this algebra has tracial rank zero.
The use of H.\  Lin's classification theorem,
Theorem~5.2 of~\cite{Ln15},
remains essential, because we can only relate tracial rank zero
to crossed products and fixed point algebras of actions by
finite cyclic groups, not the property of being an AT algebra
or even an AH algebra.

\begin{rmk}\label{NoHg}
Since the paper~\cite{Hg} remains unpublished, it is worth pointing
out that the proof of Theorem~\ref{NCTIsAT} does not actually
depend on this paper.
In the proof of Lemma~\ref{2D}, we need to know that the
ordinary irrational rotation algebras have real rank zero,
but this follows from Remark~6 of~\cite{EE}.
In the proof of Proposition~\ref{IrratAPR}, we need to
know that traces determine order on \pj s in $A_{\te}$ whenever
$A_{\te}$ is simple.
The proof of this in~\cite{BKR} does not rely on~\cite{Hg},
and in any case an independent proof (valid whenever
$\te$ is not purely rational) is contained in~\cite{Rf1}.
And in the application of Theorem~6.4 of~\cite{Ks4} in
the proof of Proposition~\ref{SingleStep}, we use the fact
that $A_{\af}$  has a unique tracial state, rather than
real rank zero, to show that Kishimoto's conditions hold.
\end{rmk}

One might also hope to prove Lemma~\ref{IrratAPR}
using Theorem~\ref{ARPFromPosElts2}, as is done for the
noncommutative Fourier transform in
Lemma~\ref{NCFT1} and Proposition~\ref{FTHasARP}.
However, Theorem~\ref{ARPFromPosElts2} requires that one know
ahead of time that the algebra involved has tracial rank zero.

Recall that the opposite algebra $A^{\mathrm{op}}$ of a \ca\  $A$
is the algebra $A$ with the multiplication reversed
but all other operations, including the scalar multiplication, the same.

\begin{cor}\label{IsomOpp}
Let $\te$ be a nondegenerate skew symmetric real $d \times d$ matrix,
with $d \geq 2$.
Then $(A_{\te})^{\mathrm{op}} \cong A_{\te}$.
\end{cor}

\begin{proof}
Every simple AT algebra $A$ with real rank zero is isomorphic to
its opposite algebra,
because the ordered K-theory of $A^{\mathrm{op}}$
is the same as the ordered K-theory of $A$.
\end{proof}

As far as we know,
it is unknown whether $(A_{\te})^{\mathrm{op}} \cong A_{\te}$
for general degenerate $\te$. 

\begin{rmk}\label{Morita}
In~\cite{RS}, certain Morita equivalences between
higher dimensional noncommutative toruses were exhibited.
Combining H.\  Lin's classification theorem for simple nuclear
\ca s with tracial rank zero and the computation of the range
of the trace on $K_0 (A_{\te})$ in~\cite{El0},
one should be able to completely determine the Morita equivalence
classes of simple higher dimensional noncommutative toruses.
We do not carry this out here.
\end{rmk}

\section{The noncommutative Fourier transform}\label{Sec:NCFT}

\indent
For $\te \in \R$ let $A_{\te}$ be the ordinary rotation algebra, the
universal \ca\  generated by unitaries $u$ and $v$
satisfying $v u = \exp (2 \pi i \te) u v$.
The group ${\mathrm{SL}}_2 (\Z)$ acts on $A_{\te}$ by sending the
matrix
\[
n = \left( \begin{array}{cc} n_{1, 1} & n_{1, 2} \\
           n_{2, 1} & n_{2, 2} \end{array} \right)
\]
to the automorphism determined by
\[
\af_n (u)
 = \exp (\pi i n_{1, 1} n_{2, 1} \te) u^{n_{1, 1}} v^{n_{2, 1}}
\andeqn
\af_n (v)
 = \exp (\pi i n_{1, 2} n_{2, 2} \te) u^{n_{1, 2}} v^{n_{2, 2}}.
\]
For fixed $\te$,
the actions of finite cyclic subgroups of ${\mathrm{SL}}_2 (\Z)$
are classified
and their fixed point algebras are found in~\cite{FW1};
the fixed point algebras are described in terms of generators
and relations in~\cite{BEEK} and~\cite{FW3}.
There are essentially only four actions:
an action of $\Zqt$ generated by the flip automorphism
\[
u \mapsto u^* \andeqn v \mapsto v^*,
\]
an action of $\Zqh$ generated by the automorphism
\[
u \mapsto e^{- \pi i \te} u^* v \andeqn v \mapsto u^*,
\]
an action of $\Zqf$
generated by the noncommutative Fourier transform
\[
u \mapsto v \andeqn v \mapsto u^*,
\]
and an action of $\Zqs$ generated by the automorphism
\[
u \mapsto v \andeqn v \mapsto e^{- \pi i \te} u^* v.
\]
One can define a few other automorphisms of the same order by changing
the scalar factors in the formulas above, but one does not get
anything essentially new.
See Proposition~21 of~\cite{FW4}.

Here, we will primarily be concerned with the crossed products,
although, as we will see,
for $\te$ irrational the crossed products are all Morita equivalent
to the corresponding fixed point algebras.
The most is known about the crossed product by the flip:
its (unordered) K-theory has been computed in~\cite{Kj},
and the crossed product has been proved to be an AF algebra~\cite{BK}.
The next best understood case is that of the
noncommutative Fourier transform, which has been intensively
studied in a series of papers culminating in~\cite{Wl2}.
It is proved there that for ``most'' irrational $\te$,
the crossed product of $A_{\te}$ by the noncommutative Fourier transform
has tracial rank zero in the sense of~\cite{LnTAF},
and that for ``most'' of those values of $\te$, it is in
addition an AF algebra.
There are three independent parts of the proof that this crossed
product if AF for ``most'' $\te$:
the proof that it satisfies the Universal Coefficient Theorem,
the proof that it has tracial rank zero,
and the computation of the K-theory.
In~\cite{Wl1} the K-theory computation is done for ``most'' $\te$,
and in~\cite{Wl2} the crossed product is shown to have tracial
rank zero for ``most'' $\te$, and to satisfy the
Universal Coefficient Theorem for all $\te$.
We prove here, using completely different methods,
that the crossed product and the fixed point algebra
have tracial rank zero for all irrational $\te$,
not just ``most'' $\te$.
Instead of the heavy use of theta functions in~\cite{Wl2},
we prove that the action of $\Zqf$ has the \aRp.
We conclude that the crossed product is a simple AH algebra
with real rank zero for all irrational $\te$.
To show that the crossed product is AF for all irrational $\te$
requires in addition an improvement on the K-theory calculation.
This will be done elsewhere,
in joint work with Wolfgang L\"{u}ck and Sam Walters.

The same methods show that the other
three finite cyclic group actions also have the \aRp,
and that their crossed products also satisfy the
Universal Coefficient Theorem.
These crossed products are therefore also simple AH algebras
with real rank zero, for all irrational $\te$.
Combining this with the K-theory computation in~\cite{Kj},
one obtains a new proof that for irrational $\te$ the crossed
product by the flip is AF.
This computation is subsumed by the work of the next section,
in which we show that the crossed product of any simple
higher dimensional noncommutative torus by the analog of the
flip is AF.
For the actions of $\Zqh$ and of $\Zqs$, we don't
know enough about the K-theory to conclude that any of the
crossed products are AF.

In this section, we do the case of the
noncommutative Fourier transform in detail.
In the next section, we give a brief description for the actions of
$\Zqh$ and $\Zqs$.

We start with a general method for proving
that an action has the \aRp\  %
when the algebra has tracial rank zero and a unique tracial state.
Most of the work is contained in the first lemma.
In Theorem~\ref{ARPFromPosElts2},
we give a useful further weakening of the hypotheses.

As Masaki Izumi has pointed out,
we could prove outerness of the actions we are interested in
by using the fact that the corresponding actions on the
factor obtained from the trace representation are outer.
We have decided to keep the original proof because we
hope it points the way to generalizations of the methods,
and perhaps to the right version of the \aRp\  for actions
on \ca s with few \pj s.

\begin{lem}\label{ARPFromPosElts}
Let $A$ be an infinite dimensional \suca\  with tracial rank zero and
with a unique tracial state $\ta$.
Let $\af \in \Aut (A)$ satisfy $\af^n = \id_A$.
Suppose that for every finite set $F \S A$ and every $\ep > 0$
there are positive elements
$a_0, a_1, \dots, a_{n - 1} \in A$ with $0 \leq a_j \leq 1$ such that:
\begin{itemize}
\item[(1)]
$a_j a_k = 0$ for $j \neq k$.
\item[(2)]
$\| \af (a_j) - a_{j + 1} \| < \ep$ for $0 \leq j \leq n - 2$.
\item[(3)]
$\| a_j c - c a_j \| < \ep$ for $0 \leq j \leq n - 1$ and all $c \in F$.
\item[(4)]
$\ta \left( 1 - \sum_{j = 0}^{n - 1} a_j \right) < \ep$.
\end{itemize}
Then the action of $\Zqn$ generated by $\af$ has the
tracial Rokhlin property.
\end{lem}

\begin{proof}
We verify the conditions of Definition~\ref{ARPDfn}.
Note that $\ta \circ \af = \ta$.
Let $F \S A$ be a finite set, let $\ep > 0$, let $N \in \N$,
and let $x \in A$ be a nonzero positive element.
\Wolog\  $F$ is $\af$-invariant and $\| a \| \leq 1$ for all $a \in F$;
also, $\ep < 1$.
Since $A$ has real rank zero (Theorem~\ref{TAFProp}),
there is a nonzero \pj\  $q \in {\overline{x A x}}$.
Since also $A$ is infinite dimensional, simple, and unital,
Theorem~1.1(i) of \cite{Zh7} provides a nonzero \pj\  $q_0 \in A$
such that
\[
\ta (q_0) < \min \left( \frac{\ta (q)}{2 n}, \,
    \frac{1}{8 n}, \, \frac{1}{4 N n^2} \right).
\]
By Lemma~2.3 of~\cite{BKR}, applied to $B = \C^{n + 1}$
(see Definition~2.2 of~\cite{BKR}), there is $\dt > 0$
such that whenever
$D$ is a \ca\  and $p_1, p_2, \dots, p_n \in D$ are \pj s such that
$\| p_j p_k \| < \dt$ for $1 \leq j, \, k \leq n$ with $j \neq k$,
then there are \mops\  $q_1, q_2, \dots, q_n \in D$ such that
$\| q_k - p_k \| < \frac{1}{10} \ep$ for $1 \leq k \leq n$.
Set
\[
\ep_0 = \min \left( \frac{\dt}{5}, \,
    \frac{1}{8 n (n + 4)}, \,
    \frac{1}{4 N n^2 (n + 4)}, \, \frac{\ta (q)}{2 n (n + 4)} \right).
\]

By Proposition~\ref{TAFCond},
there is a \pj\  $p \in A$ and a finite
dimensional unital subalgebra $E \S p A p$ such that:
\bei
\item[(1)]
$\| p a - a p \| < \ts{ \frac{1}{10} } \ep$ for all $a \in F$.
\item[(2)]
For every $a \in F$ there exists $b \in E$ such that
$\| p a p - b \| < \ts{ \frac{1}{10} } \ep$.
\item[(3)]
$1 - p$ is \mvnt\  to a \pj\  in $q_0 A q_0$.
\eei
Let $p^{(1)}, \, p^{(2)}, \dots, \, p^{(s)}$ be the minimal
central \pj s of $E$, write $p^{(l)} E p^{(l)} = M_{r (l)}$
with $r (l) \in \N$, and let
$\rule{0em}{2.3ex}
   \left\{ \rsz{ p^{(l)}_{j, k} \colon 1 \leq j, \, k \leq r (l) }
        \right\}$
be a system of matrix units for $p^{(l)} E p^{(l)}$.
Define $P \colon p A p \to p A p$ by
\[
P (a) = \sum_{l = 1}^s \sum_{k = 1}^{r (l)}
     p^{(l)}_{k, 1} a p^{(l)}_{1, k}.
\]
Then $P$ is a (slightly nonstandard)
conditional expectation from $p A p$ onto the
relative commutant $E' \cap p A p$ of $E$ in $p A p$.

Define
\[
S = \ts{ \left\{  \rsz{p^{(l)}_{j, k}} \colon {\mbox{$1 \leq l \leq s$
    and $1 \leq j, \, k \leq r (l)$}} \right\} },
\]
which is a complete system of matrix units for $E$,
with cardinality $\card (S)$.
We claim that if $a \in p A p$ then
\[
\| a - P (a) \| < \card (S) \max_{v \in S} \| v a - a v \|.
\]
To see this, first observe that
\[
\ts{ \left\| \rsz{p^{(l)}_{k, k} a p^{(l)}_{k, k}
                    - p^{(l)}_{k, 1} a p^{(l)}_{1, k} } \right\|
       = \left\|  \rsz{p^{(l)}_{k, k}
             \left[  \rsz{a p^{(l)}_{k, 1} - p^{(l)}_{k, 1} a} \right]
                               p^{(l)}_{1, k} } \right\|
       \leq \left\|  \rsz{a p^{(l)}_{k, 1} - p^{(l)}_{k, 1} a} \right\|. }
\]
Using this, we estimate:
\begin{align*}
\| a - P (a) \|
  & \leq \sum_{l = 1}^s \sum_{j, k = 1}^{r (l)}
     \left\|  \rsz{p^{(l)}_{j, j} [a - P (a) ] p^{(l)}_{k, k} }
                                           \right\|   \\
  & \leq \sum_{l = 1}^s \left( \ssum{j \neq k}{}
               \left\| \rsz{p^{(l)}_{j, j} a p^{(l)}_{k, k}} \right\|
     + \ssum{k = 2}{r (l)}
          \left\|  \rsz{p^{(l)}_{k, k} a p^{(l)}_{k, k}
             - p^{(l)}_{k, 1} a p^{(l)}_{1, k}} \right\| \right)  \\
  & \leq \sum_{l = 1}^s \left( \ssum{j \neq k}{}
           \left\|  \rsz{p^{(l)}_{j, j} a - a p^{(l)}_{j, j}} \right\|
       + \ssum{k = 2}{r (l)}
           \left\|  \rsz{p^{(l)}_{k, 1} a - a p^{(l)}_{k, 1}} \right\| \right) \\
  & \leq \card (S) \max_{v \in S} \| v a - a v \|.
\end{align*}
This proves the claim.

Define \cfn s $g_1, \, g_2, \, g_3 \colon \R \to [0, 1]$ by
\[
g_j (t) = \left\{ \begin{array}{ll}
     0     & \hspace{3em}  t \leq \frac{1}{3} (j - 1) \\
     3 t   & \hspace{3em}
              \frac{1}{3} (j - 1) \leq t \leq \frac{1}{3} j  \\
     1     & \hspace{3em}  \frac{1}{3} j \leq t
    \end{array} \right.
\]
for $j = 1, \, 2, \, 3$.
Then $g_1 g_2 = g_2$ and $g_2 g_3 = g_3$.
Using polynomial approximations to $g_1$, $g_2$, and $g_3$,
choose $\ep_1 > 0$ so small that
whenever $D$ is a \ca\  and $a, \, b \in D$ satisfy
$0 \leq a, \, b \leq 1$ and
\[
\| a - b \| < \max (n, \, 2 + \card (S) ) \ep_1,
\]
then
\[
\| g_1 (a) - g_1 (b) \| < \ep_0, \,\,\,\,\,\,
\| g_2 (a) - g_2 (b) \| < \ep_0,
\andeqn \| g_3 (a) - g_3 (b) \| < \ep_0.
\]

Apply the hypothesis with $S \cup \{ p \}$ in place of $F$ and
with $\min \left(\ep_1, \, \frac{1}{16} \ep_0^2 \right)$
in place of $\ep$, obtaining $a_0, a_1, \dots, a_{n - 1}$.

We have
$\rule{0em}{2.3ex}\left\| \rsz{(p a_0 p) p^{(l)}_{j, k}
         - p^{(l)}_{j, k} (p a_0 p)} \right\| < \ep_1$
for all $l$, $j$, and $k$,
whence
\[
\| P (p a_0 p) - p a_0 p \| < \card (S) \ep_1.
\]
Also,
\begin{align*}
\| a_0 - [(1 - p) a_0 (1 - p) + p a_0 p] \|
&  \leq \| p a_0 (1 - p) \| + \| (1 - p) a_0 p \|  \\
&  \leq 2 \| a_0 p - p a_0 \| < 2 \ep_1.
\end{align*}
Set $c = (1 - p) a_0 (1 - p) + P (p a_0 p)$.
It follows that
\[
\| a_0 - c \| < [2 + \card (S)] \ep_1.
\]
Evaluating functional calculus in the appropriate corners
on the right, we have
\[
g_2 (c) = g_2 ( (1 - p) a_0 (1 - p) ) + g_2 (P (p a_0 p)).
\]
Therefore, by the choice of $\ep_1$, we have
\[
\| g_2 (a_0) - [ g_2 ( (1 - p) a_0 (1 - p) ) + g_2 (P (p a_0 p)) ] \|
  < \ep_0.
\]
With $b = g_2 (P (p a_0 p))$,
which is equal to $g_2 (c) p \in E' \cap p A p$, we then get
\[
\| g_1 (a_0) b - b \|
  = \| [g_1 (a_0) g_2 (c) - g_2 (c) ] p \|
  \leq 2 \| g_2 (c) - g_2 (a_0) \| < 2 \ep_0.
\]
Similarly, we have $\| b g_1 (a_0) - b \| < 2 \ep_0$.

An inductive argument gives, for $0 \leq j \leq n - 1$,
the first inequality in the estimate
\[
\| \af^j (a_0) - a_j \| < j \ep_1 \leq n \ep_1.
\]
Therefore, again by the choice of $\ep_1$, we have
\[
\| \af^j ( g_1 (a_0)) - g_1 (a_j) \| < \ep_0
\]
for $0 \leq j \leq n - 1$.
Since $g_1 (a_0) g_1 (a_j) = 0$ for $1 \leq j \leq n - 1$,
it follows that
\[
\| \af^j ( g_1 (a_0)) g_1 (a_0) \| < \ep_0
\]
for those $j$.
Consequently, using $\| g_1 (a_0) b - b \| < \ep_0$
and $\| b \cdot g_1 (a_0) - b \| < \ep_0$, for $1 \leq j \leq n - 1$ we get
\[
\| \af^j (b) b \|
  < 4 \ep_0 + \| \af^j (b) \af^j ( g_1 (a_0)) g_1 (a_0) b \|
  < 5 \ep_0.
\]

We note that
\[
E' \cap p A p \cong \bigoplus_{l = 1}^s p^{(l)}_{1, 1} A p^{(l)}_{1, 1},
\]
which is a direct sum of \hsa s in $A$ and therefore has real rank zero.
So there exists a \pj\  $f$
in the \hsa\  of $E' \cap p A p$ generated by $g_3 (P (p a_0 p))$
such that $\| f g_3 (P (p a_0 p)) - g_3 (P (p a_0 p)) \| < \ep_0$.
Since $b \cdot g_3 (P (p a_0 p)) = g_3 (P (p a_0 p))$,
we get $b f = f b = f$.
Therefore $\| \af^j (f) f \| < 5 \ep_0$.
Since $\af^n = \id_A$, we get $\| \af^j (f) \af^k (f) \| < 5 \ep_0$
for $0 \leq j, \, k \leq n - 1$ with $j \neq k$.
By the choice $\ep_0 \leq \frac{1}{5} \dt$ at the beginning of the proof,
there are \mops\  $e_0, e_1, \dots, e_{n - 1} \in A$ such that
$\| e_j - \af^j (f) \| < \frac{1}{10} \ep$ for $0 \leq j \leq n - 1$.
We immediately get
\[
\| \af (e_j) - e_{j + 1} \|
  \leq \| \af (e_j) - \af^{j + 1} (f) \|
        + \| \af^{j + 1} (f) - e_{j + 1} \|
  < \ts{ \frac{2}{10} } \ep < \ep
\]
for $0 \leq j \leq n - 2$.
This is Condition~(1) of Definition~\ref{ARPDfn}.

We now prove Condition~(2).
Let $a \in F$ and let $0 \leq j \leq n - 1$.
Because we assumed $F$ is $\af$-invariant, we have $\af^{-j} (a) \in F$.
Choose $d_0 \in E$ such that
$\| d_0 - p \af^{-j} (a) p \| < \ts{ \frac{1}{10} } \ep$.
Set $d = d_0 + (1 - p) \af^{-j} (a) (1 - p)$.
Since $\| a \| \leq 1$ we get $\| d \| < 1 + \frac{1}{10} \ep$.
Also
\begin{align*}
\| \af^{-j} (a) - d \|
& \leq \| \af^{-j} (a)
         - [(1 - p) \af^{-j} (a) (1 - p) + p \af^{-j} (a) p] \|
      + \| p \af^{-j} (a) p - d_0 \|          \\
& \leq 2 \| \af^{-j} (a) p - p \af^{-j} (a) \|
      + \| p \af^{-j} (a) p - d_0 \|
  < \ts{ \frac{2}{10} } \ep + \ts{ \frac{1}{10} } \ep
  = \ts{ \frac{3}{10} } \ep.
\end{align*}
By construction, $f \in E' \cap p A p$, so $f$ commutes with $d$.
{}From $\| f - \af^{-j} (e_j) \| < \ts{ \frac{1}{10} } \ep$, we get
\[
\| \af^{-j} (e_j) d - d \af^{-j} (e_j) \|
   \leq 2 \| d \| \cdot \| f - \af^{-j} (e_j) \|
   < 2 \left( 1 + \ts{ \frac{1}{10} } \ep \right) \ts{ \frac{1}{10} } \ep.
\]
Therefore, since $\ep < 1$,
\begin{align*}
\| e_j a - a e_j \|
  &  \leq 2 \| \af^{-j} (a) - d \|
         + \| \af^{-j} (e_j) d - d \af^{-j} (e_j) \|  \\
  &  < \ts{ \frac{6}{10} } \ep
          + 2 (1 + \ts{ \frac{1}{10} } \ep) \ts{ \frac{1}{10} } \ep
     \leq \ts{ \frac{6}{10} } \ep + \ts{ \frac{4}{10} } \ep = \ep.
\end{align*}
This is Condition~(2) of Definition~\ref{ARPDfn}.

It remains to verify Conditions~(3) and~(4) of Definition~\ref{ARPDfn}.
This requires some work.
Let $X$ be the maximal ideal space of the unital \ca\  $C$ generated by
$a_0, \, a_1, \, \dots, \, a_{n - 1}$,
and let $h_0, \, h_1, \, \dots, \, h_{n - 1} \colon X \to [0, 1]$
be the elements of $C (X)$ corresponding to
$a_0, \, a_1, \, \dots, \, a_{n - 1}$.
Let $\mu$ be the probability measure on $X$ such that if $h \in C (X)$
corresponds to an element $a \in C$, then
$\int_X h \, d \mu = \ta (a)$.
Then the functions $h_j$ satisfy $h_j h_k = 0$ for $j \neq k$, and
\[
0 \leq \sum_{j = 0}^{n - 1} h_j \leq 1
\andeqn \sum_{j = 0}^{n - 1} \int_X h_j \, d \mu
   > 1 - \min \left(\ep_1, \, \ts{ \frac{1}{16} } \ep_0^2 \right).
\]

Set
\[
T_j = \left\{ x \in X \colon
   h_j (x) \geq 1 - \ts{ \frac{1}{4} } \ep_0 \right\} \S X.
\]
For $x \not\in \bigcup_{j = 0}^{n - 1} T_j$ we have
$\sum_{j = 0}^{n - 1} h_j (x) < 1 - \ts{ \frac{1}{4} } \ep_0$, whence
\begin{align*}
1 - \ts{ \frac{1}{16} } \ep_0^2
& \leq 1 - \min \left(\ep_1, \, \ts{ \frac{1}{16} } \ep_0^2 \right)
  < \int_X \left( \ssum{j = 0}{n - 1} h_j \right) \, d \mu   \\
& \leq \sum_{j = 0}^{n - 1} \mu (T_j)
      + \left( 1 - \ts{ \frac{1}{4} } \ep_0 \right)
           \left( 1 - \ssum{j = 0}{n - 1} \mu (T_j) \right)   \\
& = 1 - \frac{\ep_0}{4}
        \left( 1 - \ssum{j = 0}{n - 1} \mu (T_j) \right).
\end{align*}
It follows that
\[
\sum_{j = 0}^{n - 1} \mu (T_j)
  > 1 - \left( \frac{4}{\ep_0} \right) \left( \frac{\ep_0^2}{16} \right)
  = 1 - \ts{ \frac{1}{4} } \ep_0.
\]
Since $g_3 \circ h_j \geq 1 - \ts{ \frac{3}{4} } \ep_0$ on $T_j$,
we get
\[
\sum_{j = 0}^{n - 1} \ta (g_3 (a_j))
  = \sum_{j = 0}^{n - 1} \int_X (g_3 \circ h_j) \, d \mu
  \geq \left( 1 - \ts{ \frac{1}{4} } \ep_0 \right)
       \left( 1 - \ts{ \frac{3}{4} } \ep_0 \right)
  > 1 - \ep_0.
\]

By the choice of $\ep_1$, we have
$\| g_3 (a_j) - \af (g_3 (a_{j + 1})) \| < \ep_0$
for $0 \leq j \leq n - 2$.
Since $\af$ preserves the trace, we get
$\ta (g_3 (a_{j + 1})) < \ta (g_3 (a_j)) + \ep_0$, so inductively
$\ta (g_3 (a_j)) < \ta (g_3 (a_0)) + j \ep_0$.
Therefore
\begin{align*}
1 - \ep_0
& < \sum_{j = 0}^{n - 1} \ta (g_3 (a_j))
  < \sum_{j = 0}^{n - 1} [\ta (g_3 (a_0)) + j \ep_0 ]  \\
& < n \ta (g_3 (a_0)) + \ts{ \frac{1}{2} } n (n - 1) \ep_0
  \leq n \ta (g_3 (a_0)) + (n^2 - 1) \ep_0,
\end{align*}
that is,
\[
\ta (g_3 (a_0)) > \ts{ \frac{1}{n} } - n \ep_0.
\]

For the same reason as in a similar argument with $g_2$ earlier
in the proof, we get
\[
\| g_3 (a_0)
   - [ g_3 ( (1 - p) a_0 (1 - p) ) + g_3 (P (p a_0 p)) ] \|
 < \ep_0.
\]
So, using $1 - p \precsim q_0$ at the last step,
\begin{align*}
\ta ( g_3 (P (p a_0 p)) )
& > \ta (g_3 (a_0)) - \ta (g_3 ( (1 - p) a_0 (1 - p) ) ) - \ep_0  \\
& > \ts{ \frac{1}{n} } - n \ep_0 - \ta (1 - p) - \ep_0
  > \ts{ \frac{1}{n} } - (n + 1) \ep_0 - \ta (q_0).
\end{align*}
Recalling that
$\| f g_3 (P (p a_0 p)) - g_3 (P (p a_0 p)) \| < \ep_0$, we get
\[
\| f g_3 (P (p a_0 p)) f - g_3 (P (p a_0 p)) \| < 2 \ep_0,
\]
so that
\[
\ta (f) \geq \ta ( f g_3 (P (p a_0 p)) f )
 > \ta ( g_3 (P (p a_0 p)) ) - 2 \ep_0
 > \ts{ \frac{1}{n} } - (n + 3) \ep_0 - \ta (q_0).
\]

Because $\ep < 1$,
for each $j$ the \pj\  $e_j$ is unitarily equivalent to $\af^j (f)$,
whence $\ta (e_j) = \ta (\af^j (f)) = \ta (f)$.
Combining this with
\[
\ta (q_0) < \frac{\ta (q)}{2 n}
 \andeqn \ep_0 < \frac{\ta (q)}{2 n (n + 4)},
\]
we get
\[
\sum_{j = 0}^{n - 1} \ta (e_j) = n \ta (f)
  > 1 - n (n + 3) \ep_0 - n \ta (q_0)
  > 1 - \ta (q).
\]
Therefore $e = \sum_{j = 0}^{n - 1} e_j$ satisfies
$\ta (1 - e) < \ta (q)$.
Since in simple \ca s with tracial rank zero, the traces determine the
order on \pj s (Theorem~\ref{TAFProp}), it follows that
$1 - e \precsim q$.
Since $q \in {\overline{x A x}}$, this proves Condition~(3)
of Definition~\ref{ARPDfn}.

If we use instead
\[
\ta (q_0) < \frac{1}{4 N n^2}
 \andeqn \ep_0 < \frac{1}{4 N n^2 (n + 4)},
\]
we get instead $\ta (e) > 1 - \frac{1}{2 N n}$,
or $\ta (1 - e) < \frac{1}{2 N n}$.
Combining our estimate on $\ta (f)$ with
\[
\ta (q_0) < \frac{1}{8 n} \andeqn \ep_0 < \frac{1}{8 n (n + 4)},
\]
we get $\ta (f) > \frac{1}{2 n}$.
For every $j$, recalling from above that $\ta (e_j) = \ta (f)$,
we have
\[
\ta (e_j) = \ta (f) > \ts{ \frac{1}{2 n} } > N \ta (1 - e).
\]
Because the traces determine the order on \pj s,
this implies that there are $N$ \mops\  in $f A f$,
each \mvnt\  to $1 - e$.
We have proved Condition~(4) of Definition~\ref{ARPDfn}.
\end{proof}

\begin{thm}\label{ARPFromPosElts2}
Let $A$ be an infinite dimensional \suca\  with tracial rank zero and
with a unique tracial state $\ta$.
Let $\af \in \Aut (A)$ satisfy $\af^n = \id_A$.
Suppose that for every finite set $F \S A$ and every $\ep > 0$
there are positive elements
$a_0, a_1, \dots, a_{n - 1} \in A$ with $0 \leq a_j \leq 1$ such that:
\begin{itemize}
\item[(1)]
$\| a_j a_k \| < \ep$ for $j \neq k$.
\item[(2)]
$\| \af (a_j) - a_{j + 1} \| < \ep$ for $0 \leq j \leq n - 2$.
\item[(3)]
$\| a_j c - c a_j \| < \ep$ for $0 \leq j \leq n - 1$ and all $c \in F$.
\item[(4)]
$\left| \ta \left( 1 - \sum_{j = 0}^{n - 1} a_j \right) \right| < \ep$.
\end{itemize}
Then the action of $\Zqn$ generated by $\af$ has the
tracial Rokhlin property.
\end{thm}

\begin{proof}
We prove that the hypotheses imply those of Lemma~\ref{ARPFromPosElts}.
Let $F \S A$ be finite and let $\ep > 0$.
\Wolog\  $\| c \| \leq 1$ for all $c \in F$.
Let $B$ be the universal (nonunital) \ca\  generated by
selfadjoint elements $d_0, \, d_1, \, \dots, \, d_{n - 1}$ subject to the
relations $0 \leq d_j \leq 1$ and $d_j d_k = 0$ for $j \neq k$.
Then $B$ is isomorphic to the cone over $\C^n$, and is
hence a projective \ca.
(See Lemmas 8.1.3 and 10.1.5 and Theorem 10.1.11 of~\cite{Lr}.)
Therefore it is semiprojective (Definition 14.1.3 of~\cite{Lr}).
Thus, using Theorem 14.1.4 of~\cite{Lr}
(see Definition 14.1.1 of~\cite{Lr}, and take $B$ there to be $\{ 0 \}$),
there is $\dt > 0$ such that whenever $D$ is a \ca\  and
$d_0, \, d_1, \, \dots, \, d_{n - 1} \in D$
are positive elements with $0 \leq d_j \leq 1$ and such that
$\| d_j d_k \| < \dt$ for $j \neq k$, then there are
positive elements
$a_0, a_1, \dots, a_{n - 1} \in D$ with $0 \leq a_j \leq 1$ such that
$a_j a_k = 0$ for $j \neq k$ and
\[
\| a_j - d_j \|
  < \min \left( \ts{ \frac{1}{4} } \ep, \,
                 \ts{ \frac{1}{2} } n^{-1} \ep \right)
\]
for all $j$.
Apply our hypotheses with $F$ as given and with
$\min \left( \ts{ \frac{1}{2} } \ep, \, \dt \right)$ in place of $\ep$,
obtaining $d_0, \, d_1, \, \dots, \, d_{n - 1} \in A$,
and let $a_0, a_1, \dots, a_{n - 1} \in A$ be as above.
The relation $a_j a_k = 0$ for $j \neq k$
is Condition~(1) of Lemma~\ref{ARPFromPosElts}.
For~(2), estimate
\[
\| \af (a_j) - a_{j + 1} \|
  \leq \| \af (d_j) - d_{j + 1} \| + \| a_j - d_j \|
            + \| a_{j + 1} - d_{j + 1} \|
  < \ts{ \frac{1}{2} } \ep
            + \ts{ \frac{1}{4} } \ep + \ts{ \frac{1}{4} } \ep
  = \ep.
\]
For~(3), use $\| c \| \leq 1$ for $c \in F$ to estimate
\[
\| a_j c - c a_j \| \leq \| d_j c - c d_j \| + 2 \| a_j - d_j \|
     < \ts{ \frac{1}{2} } \ep + 2 \left( \ts{ \frac{1}{4} } \ep \right)
     = \ep.
\]
For~(4), estimate
\[
\ta \left( 1 - \ssum{j = 0}{n - 1} d_j \right)
  \leq \left| \ta \left( 1 - \ssum{j = 0}{n - 1} a_j \right) \right|
            + \sum_{j = 0}^{n - 1} \| a_j - d_j \|
     < \ts{ \frac{1}{2} } \ep
            + n \left( \ts{ \frac{1}{2} } n^{-1} \ep \right)
     = \ep.
\]
This completes the proof.
\end{proof}

In the rest of this section, we will use the \ct\  field of
rotation algebras $A_{\te}$, with section algebra
$A$ equal to the \ca\  of the discrete Heisenberg group,
as described in Theorem~\ref{FieldOfRotAlg},
but we revert to tradition and use $u_{\te}$ and $v_{\te}$
for the unitary generators of $A_{\te}$
and $u$, $v$, and $z$ for the unitary generators of $A$.
We let $\ta_{\te}$ be the standard trace on $A_{\te}$,
as in Lemma~\ref{ContOfTrace}.
Since this \ct\  field will be the only one used in this section,
we continue to follow Notation~\ref{SectionNot},
letting $\Gm (E)$ be the set of \ct\  sections of this field
over the set $E$, and writing $a (\te)$ for $a (\exp (2 \pi i \te))$
when $a$ is a section.

\begin{dfn}\label{NcFTDfn}
Let $A$ and $A_{\te}$ be as in Theorem~\ref{FieldOfRotAlg},
with notation as above.
The {\emph{noncommutative Fourier transform}} on $A_{\te}$
is the unique automorphism $\sm_{\te}$ of order~$4$ satisfying
$\sm_{\te} (u_{\te}) = v_{\te}$
and $\sm_{\te} (v_{\te}) = u_{\te}^*$.
The {\emph{noncommutative Fourier transform}} on $A$
is the unique automorphism $\sm$ of order~$4$ satisfying
$\sm (u) = v$, $\sm (v) = u^*$, and $\sm (z) = z$.
\end{dfn}

\begin{lem}\label{NCFT1}
Let the notation be as before and in Definition~\ref{NcFTDfn}.
Let $\ep > 0$.
Then there is $\dt > 0$ such that,
with $I = \exp (2 \pi i (- \dt, \, \dt) )$,
there are \ct\  sections $b_0, b_1, b_2, b_3 \in \Gm (I)$ satisfying:
\begin{itemize}
\item[(1)]
$0 \leq b_j (\te) \leq 1$ for $0 \leq j \leq 3$ and $| \te | < \dt$.
\item[(2)]
$\| b_j (\te) b_k (\te) \| < \ep$ for $j \neq k$ and $| \te | < \dt$.
\item[(3)]
Setting $b_4 = b_0$,
we have $\| \sm_{\te} (b_j (\te)) - b_{j + 1} (\te) \| < \ep$
for $0 \leq j \leq 3$ and $| \te | < \dt$.
\item[(4)]
$| \ta_{\te} (1 - b_0 (\te) - b_1 (\te) - b_2 (\te) - b_3 (\te) ) |
      < \ep$
for $| \te | < \dt$.
\end{itemize}
\end{lem}

\begin{proof}
Set $\ep_0 = \frac{1}{16} \ep$.
Choose a \cfn\  $h \colon S^1 \to [0, 1]$ such that
$h (\zt) = 0$ for ${\mathrm{Im}} (\zt) \leq 0$ and
$h (\zt) = 1$ when $\zt = \exp (2 \pi i t)$ with
$t \in \left[ \ep_0, \, \frac{1}{2} - \ep_0 \right]$.
Identifying $A_0$ with $C (S^1 \times S^1)$ in the obvious way,
define elements $g_0, g_1, g_2, g_3 \in A_0$ by
\[
g_0 (\zt_1, \zt_2) = h (\zt_1) h (\zt_2), \,\,\,\,\,\,
g_1 (\zt_1, \zt_2) = h ({\overline{\zt}}_1) h (\zt_2),
\]
\[
g_2 (\zt_1, \zt_2) = h ({\overline{\zt}}_1) h ({\overline{\zt}}_2),
\andeqn
g_3 (\zt_1, \zt_2) = h (\zt_1) h ({\overline{\zt}}_2).
\]
Clearly $0 \leq g_j \leq 1$, $\sm_{0} (g_j) = g_{j + 1}$, and
$g_j g_k = 0$ for $j \neq k$.
Moreover,
\[
\ta_0 (g_0) \geq \left( \ts{ \frac{1}{2} } - 2 \ep_0 \right)^2
  > \ts{ \frac{1}{4} } - 2 \ep_0,
\]
and $\ta (g_j) = \ta (g_0)$ for all $j$, so
\[
\ta_0 ( 1 - g_0 - g_1 - g_2 - g_3) < 8 \ep_0
 = \ts{ \frac{1}{2} } \ep.
\]

Since $\ev_0 \colon A \to A_0$ is surjective, there exist
selfadjoint elements $a_0, a_1, a_2, a_3 \in A$
such that $\ev_0 (a_j) = g_j$ and $0 \leq a_j \leq 1$
for $0 \leq j \leq 3$.
Since $\ev_0 (a_j) \ev_0 (a_k) = 0$ for $j \neq k$
and $\ev_0 (\sm (a_j)) - \ev_0 (a_{j + 1}) = 0$
for $0 \leq j \leq 2$,
Theorem~\ref{FieldOfRotAlg} provides $\dt_1 > 0$ such that
whenever $| \te | < \dt_1$, we have
$\| \ev_{\te} (a_j) \ev_{\te} (a_k) \| < \ep$ for $j \neq k$ and
$\| \ev_{\te} (\sm (a_j)) - \ev_{\te} (a_{j + 1}) \| < \ep$
for $0 \leq j \leq 3$.
Since
\[
| \ta_0 (1 - \ev_0 (a_0) - \ev_0 (a_1)
         - \ev_0 (a_2) - \ev_0 (a_3)) |
 < \ts{ \frac{1}{2} } \ep,
\]
Lemma~\ref{ContOfTrace} provides $\dt_2 > 0$ such that
\[
| \ta_{\te} (1 - \ev_{\te} (a_0) - \ev_{\te} (a_1)
         - \ev_{\te} (a_2) - \ev_{\te} (a_3)) |
 < \ts{ \frac{1}{2} } \ep
\]
for $| \te | < \dt_2$.
The lemma is now proved by taking $\dt = \min ( \dt_1, \dt_2)$
and taking $b_j$ to be the restriction of $a_j$,
regarded as a section over $S^1$ of the \ct\  field of
Theorem~\ref{FieldOfRotAlg},
to $I = \exp (2 \pi i (- \dt, \, \dt) )$.
\end{proof}

\begin{lem}\label{Dioph}
Let $\te \in \R \SM \Q$.
Then for every $\ep > 0$ there is $n \in \N$ with $n > 0$ such that
\[
\dist (n \te, \, \Z) < \ep \andeqn \dist (n^2 \te, \, \Z) < \ep.
\]
\end{lem}

\begin{proof}
Define $h \colon S^1 \times S^1 \to S^1 \times S^1$ by
\[
h (\zt_1, \zt_2)
 = ( \exp ( 2 \pi i \te) \zt_1, \, \exp ( 2 \pi i \te) \zt_1^2 \zt_2).
\]
As in the discussion preceding Proposition~1.5 of~\cite{Fr},
the map $h$ is a minimal homeomorphism of $S^1 \times S^1$,
and the forward orbit of $(1, 1)$ consists of all points
$(\exp ( 2 \pi i n \te), \, \exp ( 2 \pi i n^2 \te) )$ for
$n \in \N$.
The result is immediate from the density of the forward orbit.
\end{proof}

\begin{prp}\label{FTHasARP}
Let $\te \in \R \SM \Q$.
Then the action of $\Zqf$ on $A_{\te}$ generated by the
noncommutative Fourier transform $\sm_{\te}$ of
Definition~\ref{NcFTDfn} has the \aRp.
\end{prp}

\begin{proof}
Let the notation be as before and in Definition~\ref{NcFTDfn}.
It suffices to verify the conditions of Theorem~\ref{ARPFromPosElts2}.
Moreover, we may take the finite set $F$ to be $F = \{ u, v \}$.
We know from Theorem~\ref{NCTIsAT} that $A_{\te}$
has tracial rank zero, and we know that
$A_{\te}$ has a unique tracial state $\ta_{\te}$.

Accordingly, let $\ep > 0$.
Apply Lemma~\ref{NCFT1},
obtaining $\dt_0 > 0$, an arc $I = \exp (2 \pi i (- \dt_0, \, \dt_0) )$,
and \ct\  sections $b_0, b_1, b_2, b_3 \in \Gm (I)$
satisfying Properties (1)--(4) there.
Choose polynomials $f_0, f_1, f_2, f_3$ in four
noncommuting variables such that in $A_0 = C (S^1 \times S^1)$ we have
\[
\| f_j (u_0, u_0^*, v_0, v_0^*) - b_j (0) \| < \ts{ \frac{1}{5}} \ep.
\]
for $0 \leq j \leq 3$.
Because we are dealing with a \ct\  field
(by Theorem~\ref{FieldOfRotAlg}), there is
$\dt_1 > 0$ with $\dt_1 \leq \dt_0$ such that
for $| \et | < \dt_1$ and $0 \leq j \leq 3$ we have
\[
\| f_j (u_{\et}, u_{\et}^*, v_{\et}, v_{\et}^*) - b_j (\et) \|
  < \ts{ \frac{2}{5} } \ep.
\]

Choose (see Proposition~4.3 of~\cite{BKR}) $\dt_2 > 0$ such that
whenever $D$ is a \ca\  and $x_1, x_2, x_3, x_4, y \in D$
are elements of norm~$1$ which satisfy $\| [x_k, y] \| < \dt_2$
for all $k$, then
\[
\| [f_j (x_1, x_2, x_3, x_4), \, y ] \| < \ts{ \frac{1}{5} } \ep
\]
for $0 \leq j \leq 3$.
Now use Lemma~\ref{Dioph} to choose $n \in \N$ with $n > 0$
and $m \in \Z$ such that
\[
| \exp ( 2 \pi i n \te) - 1 | < \dt_2 \andeqn | n^2 \te - m | < \dt_1.
\]

Let $\et = n^2 \te - m$,
and let $\ps \colon A_{\et} \to A_{\te}$ be the \hm\  determined by
$\ps (u_{\et}) = u_{\te}^n$ and $\ps (v_{\et}) = v_{\te}^n$.
Set $a_j = \ps (b_j (\et))$ for $0 \leq j \leq 3$.
We verify the hypotheses (1)--(4) of Theorem~\ref{ARPFromPosElts2}.
We have $0 \leq a_j \leq 1$ because $0 \leq b_j (\et) \leq 1$.
For~(1), we have $\| a_j a_k \| < \ep$ for $j \neq k$ because
$| \et | < \dt_1 \leq \dt_0$ implies
$\| b_j (\et) b_k (\et) \| < \ep$ for $j \neq k$.
For~(2), we observe that
$\sm_{\te} \circ \ps = \ps \circ \sm_{\et}$.
Since $| \et | < \dt_0$ we have
$\| \sm_{\et} (b_j (\et)) - b_{j + 1} (\et) \| < \ep$,
and this implies $\| \sm_{\te} (a_j) - a_{j + 1} \| < \ep$.
For~(4), we note that $\et$ and $\te$ are both irrational, so
uniqueness of the traces on $A_{\te}$ and $A_{\et}$ implies
$\ta_{\te} \circ \ps = \ta_{\et}$.
Therefore $| \et | < \dt_0$ implies
\[
| \ta_{\te} (1 - a_0 - a_1 - a_2 - a_3 ) |
 = | \ta_{\et} (1 - b_0 (\et) - b_1 (\et) - b_2 (\et) - b_3 (\et) ) |
 < \ep.
\]

Finally, we prove~(3).
We have
$u_{\te}^n v_{\te} = \exp (- 2 \pi i n \te) v_{\te} u_{\te}^n$, so
\[
\| u_{\te}^n v_{\te} - v_{\te} u_{\te}^n \|
  = | \exp (- 2 \pi i n \te) - 1 | < \dt_2.
\]
For either similar reasons or trivially, we get
$\| x y - y x \| < \dt_2$ for all
\[
x \in \{ \ps (u_{\et}), \, \ps (u_{\et})^*, \, \ps (v_{\et}), \,
      \ps (v_{\et})^* \}
\andeqn y \in \{ u_{\te}, v_{\te} \}.
\]
Therefore, by the choice of $\dt_2$,
\[
\| [ \ps (f_j (u_{\et}, u_{\et}^*, v_{\et}, v_{\et}^*) ), \, y ] \|
  < \ts{ \frac{1}{5} } \ep
\]
for $y \in \{ u_{\te}, v_{\te} \}$.
Furthermore, from $| \et | < \dt_1$ we get
\[
\| \ps (f_j (u_{\et}, u_{\et}^*, v_{\et}, v_{\et}^*) ) - a_j \|
  < \ts{ \frac{2}{5} } \ep.
\]
Therefore
\begin{align*}
\| [a_j, u_{\te}] \|
&  \leq \| [ \ps (f_j (u_{\et}, u_{\et}^*, v_{\et}, v_{\et}^*) ),
                        \, u_{\te}] \|
          + 2 \| \ps (f_j (u_{\et}, u_{\et}^*, v_{\et}, v_{\et}^*) )
                         - a_j \|      \\
&  < \ts{ \frac{1}{5} } \ep + \ts{ \frac{4}{5} } \ep = \ep,
\end{align*}
and similarly $\| [a_j, v_{\te}] \| < \ep$.
This shows that $\| [a_j, y] \| < \ep$ for $y \in F$, and completes the
proof of~(3).

Having verified the hypotheses of Theorem~\ref{ARPFromPosElts2},
we conclude from that theorem that the action of $\Zqf$ generated
by $\sm_{\te}$ has the \aRp.
\end{proof}

Next, we show that the crossed product by the noncommutative
Fourier transform satisfies the Universal Coefficient Theorem.
We formulate this condition in a more convenient way than the
usual definition.
This reformulation will be more important
in Section~\ref{Sec:CrPrdFlip} than here.

\begin{dfn}\label{DL1}
Let $A$ be a separable nuclear \ca.
We say that $A$ {\emph{satisfies the \uct}}
if for every separable \ca\  $B$ such that $K_* (B)$ is an injective
abelian group, the natural map
\[
\gm = \gm_{A, B} \colon KK^* (A, B) \to \Hom ( K_* (A), K_* (B))
\]
is an isomorphism.
\end{dfn}

\begin{rmk}\label{DL2}
The definition in \cite{RSUCT} is that $A$ satisfies the \uct\  if
for every separable \ca\  $B$ (actually, every \ca\  $B$
with a countable approximate identity; see Remark~\ref{DL3} below),
there is a natural short exact sequence
\[
0 \longrightarrow \Ext_1^{\Z} ( K_* (A), K_* (B))
\longrightarrow KK^* (A, B)
\longrightarrow \Hom ( K_* (A), K_* (B))
\longrightarrow 0,
\]
in which:
\bei
\item[(1)]
The second map is the map $\gm$ of Definition~\ref{DL1}.
\item[(2)]
The first map has degree one and is the inverse of the map
\[
{\mathrm{Ker}} (\gm) \to \Ext_1^{\Z} ( K_* (A), K_* (B))
\]
which sends a class in $KK^1 (A, B)$ represented by an extension
\[
0 \longrightarrow A \longrightarrow E \longrightarrow B
\longrightarrow 0
\]
to the $\Ext$ class of the short exact sequence
\[
0 \longrightarrow K_* (A) \longrightarrow K_* (E)
\longrightarrow K_* (B) \longrightarrow 0,
\]
and is the suspension of this on $KK^0 (A, B)$.
\eei

This statement is unwieldy (both in the definition and verification)
if the maps are defined, and imprecise if they are not.
Moreover, the condition in Definition~\ref{DL1} is what one in practice
verifies.
The proof of Theorem~4.1 of \cite{RSUCT} shows that if $A$ is nuclear
and satisfies the condition of Definition~\ref{DL1},
then for every separable \ca\  $B$ the
map ${\mathrm{Ker}} (\gm) \to \Ext_1^{\Z} ( K_* (A), K_* (B))$
of~(2) above is in fact invertible, and the resulting sequence
\[
0 \longrightarrow \Ext_1^{\Z} ( K_* (A), K_* (B))
\longrightarrow KK^* (A, B)
\longrightarrow \Hom ( K_* (A), K_* (B))
\longrightarrow 0,
\]
is in fact exact.
\end{rmk}

\begin{rmk}\label{DL3}
The blanket assumption in \cite{RSUCT} on the second variable $B$
is that it has a countable approximate identity.
There seems to be a gap in the proof of Theorem~4.1 of \cite{RSUCT}
in this generality, because it is not clear that the algebra
$D_0$ occurring in the proof has a countable approximate identity.
\end{rmk}

Much of the following proof is due to Walters.
See Theorem~9.4 of~\cite{Wl2}.

\begin{lem}\label{CPbyFThasUCT}
Let $\te \in \R \SM \Q$.
Then the crossed product $\Cs{4}{A_{\te}}{\sm_{\te}}$,
with $\sm_{\te}$  as in Definition~\ref{NcFTDfn}, satisfies the
Universal Coefficient Theorem.
\end{lem}

\begin{proof}
Let $A$, $A_{\te}$, and $\ev_{\te} \colon A \to A_{\te}$ be
as in Theorem~\ref{FieldOfRotAlg}, and let $\sm \in \Aut (A)$ be as
in Definition~\ref{NcFTDfn}.
Then the crossed product $\Cs{4}{A}{\sm}$
is the \ca\  of the semidirect product of $\Zqf$ and the
discrete Heisenberg group $H$, formed from the obvious
corresponding order four automorphism of $H$.
This group is amenable.
Using Lemma~3.5 of~\cite{Tu},
we may apply Proposition~10.7 of~\cite{Tu}, specialized to groups,
to $H$, and conclude that
$\Cs{4}{A}{\sm}$ \suct.
If $\te \in \Q$, then $A_{\te}$ is type I, so
$\Cs{4}{A_{\te}}{\sm_{\te}}$ is type I by
Theorem~4.1 of~\cite{RfFG}.
Therefore it \suct.

Now let $\te_1, \, \te_2 \in [0, 1)$ satisfy $\te_1 < \te_2$.
Let $I_{\te_1, \te_2}^{(0)}$ be the set of sections of the
\ct\  field of Theorem~\ref{FieldOfRotAlg} which vanish on the
closed arc from $\exp (2 \pi i \te_1)$ to $\exp (2 \pi i \te_2)$,
and let $J_{\te_1, \te_2}^{(0)}$
be the set of sections which vanish on the other closed arc
with the same endpoints.
Let $I_{\te_1, \te_2}$ and $J_{\te_1, \te_2}$
be the crossed products of these
by the action of $\Zqf$ generated by $\sm$.
Using the maps $\ev_{\te_1}$ and $\ev_{\te_2}$,
we obtain a $\Zqf$-equivariant short exact sequence of \ca s
\[
0 \longrightarrow I_{\te_1, \te_2}^{(0)} \oplus J_{\te_1, \te_2}^{(0)}
 \longrightarrow A \longrightarrow A_{\te_1} \oplus A_{\te_2}
 \longrightarrow 0.
\]
Since crossed products preserve exact sequences
(see Lemma~2.8.2 of~\cite{Ph1}),
we obtain the short exact sequence
\begin{align*}
0 \longrightarrow I_{\te_1, \te_2} \oplus J_{\te_1, \te_2}
& \longrightarrow \Cs{4}{A}{\sm}  \\
& \longrightarrow \Csw{4}{A_{\te_1}}{\sm_{\te_1}}
           \oplus \Csw{4}{A_{\te_2}}{\sm_{\te_2}}
 \longrightarrow 0.
\end{align*}
For $\te_1, \, \te_2 \in \Q$
we have seen that the second and third terms of this sequence
satisfy the \uct,
so $I_{\te_1, \te_2} \oplus J_{\te_1, \te_2}$ \suct\  by
Proposition~2.3(a) of~\cite{RSUCT}.
It is easy to check from Definition~\ref{DL1} that
$I_{\te_1, \te_2}$ and $J_{\te_1, \te_2}$ therefore each separately
satisfy the \uct.

Now let $\te \in (0, 1)$ be arbitrary.
Choose
\[
\bt_1, \, \bt_2, \dots, \gm_1, \, \gm_2, \dots \in (0, 1) \cap \Q
\]
such that
\[
\bt_1 < \bt_2 < \cdots < \te < \cdots < \gm_2 < \gm_1
\andeqn
\limi{n} \bt_n = \limi{n} \gm_n = \te.
\]
Let $L$ be the kernel of the map of crossed products
\[
\Cs{4}{A}{\sm} \to \Cs{4}{A_{\te}}{\sm_{\te}}
\]
Then, since crossed products preserve direct limits
and exact sequences,
\[
J_{\bt_1, \gm_1} \S J_{\bt_1, \gm_1} \S \cdots
\andeqn
L
 = {\overline{ {\ts{ {\ds{\bigcup}}_{n = 1}^{\I} }} J_{\bt_n, \gm_n} }}.
\]
So $L$ \suct\  by Proposition~2.3(b) of~\cite{RSUCT},
whence $\Cs{4}{A_{\te}}{\sm_{\te}}$
\suct\  by Proposition 2.3(a) of~\cite{RSUCT}.
\end{proof}

\begin{thm}\label{CPbyFTisAH}
Let $\te \in \R \SM \Q$.
Then the crossed product $\Cs{4}{A_{\te}}{\sm_{\te}}$,
with $\sm_{\te}$  as in Definition~\ref{NcFTDfn},
is a simple AH algebra with slow dimension growth and real rank zero.
\end{thm}

\begin{proof}
We know from Theorem~\ref{NCTIsAT} that $A_{\te}$
has tracial rank zero.
Combining Proposition~\ref{FTHasARP} and Theorem~\ref{RokhTAF},
we find that the crossed product
$\Cs{4}{A_{\te}}{\sm_{\te}}$ has tracial
rank zero.
By Lemma~\ref{CPbyFThasUCT}
it satisfies the Universal Coefficient Theorem,
by Corollary~\ref{CrPrIsSimple} it is simple,
and it is clearly separable and nuclear.
The result therefore follows from Lemma~\ref{ConseqOfClass}.
\end{proof}

The following corollary is Theorem~9.3 of~\cite{Wl2}.

\begin{cor}\label{CPbyFTisAF}
There exists a dense $G_{\dt}$-set $E \S \R \setminus \Q$
such that for every $\te \in E$,
the crossed product $\Cs{4}{A_{\te}}{\sm_{\te}}$
is a simple AF~algebra.
\end{cor}

\begin{proof}
It is shown in~\cite{Wl1} that there is
a dense $G_{\dt}$-set $E \S \R \setminus \Q$
such that for every $\te \in E$,
the crossed product $\Cs{4}{A_{\te}}{\sm_{\te}}$
has trivial $K_1$ and torsion-free $K_0$.
The result therefore follows from Theorem~\ref{CPbyFTisAH}
and Lemma~\ref{ConseqOfClass}.
\end{proof}

\begin{cor}\label{FPofFTisAF}
There exists a dense $G_{\dt}$-set $E \S \R \setminus \Q$
such that for every $\te \in E$,
the fixed point algebra $A_{\te}^{\sm_{\te}}$ of the automorphism
$\sm_{\te} \in \Aut (A_{\te})$ is a simple AF algebra.
\end{cor}

\begin{proof}
For every $\te$ in the dense $G_{\dt}$-set $E \S \R \setminus \Q$ of
Corollary~\ref{CPbyFTisAF}, Proposition~\ref{FPIsSimple}
shows that $A_{\te}^{\sm_{\te}}$ is strongly Morita
equivalent to a simple AF algebra.
For these $\te$, it follows that $A_{\te}^{\sm_{\te}}$ is
a simple AF algebra.
\end{proof}

\section{Other finite cyclic group actions on the
             irrational rotation algebra}\label{Sec:OtherActions}

\indent
In this section, we consider the crossed products of the irrational
rotation algebra by the automorphisms of orders~$3$ and~$6$ coming
from the action of ${\mathrm{SL}}_2 (\Z)$.
There are four differences from the previous section.
First, the actions don't extend to an action on the C*-algebra of
the discrete Heisenberg group, and we must use in its place a slightly
larger group and corresponding larger continuous field.
Second, the analysis of the case $\te = 0$, in the proof of
the analog of Lemma~\ref{NCFT1}, is messier.
Third, rather than doing everything twice we prove a lemma
which reduces the \aRp\  for the action of $\Zqh$ to the
case of $\Zqs$.
Finally, lacking information analogous to that of~\cite{Wl1}
on the K-theory, we do not prove that any of the crossed products
are AF.

As in Section~\ref{Sec:NCFT},
for $\te \in \R$ let $A_{\te}$ be the ordinary rotation algebra, the
universal \ca\  generated by unitaries $u_{\te}$ and $v_{\te}$
satisfying $v_{\te} u_{\te} = \exp (2 \pi i \te) u_{\te} v_{\te}$.

\begin{dfn}\label{Ord6Auto}
Define an automorphism $\ph_{\te} \colon A_{\te} \to A_{\te}$
of order $6$ by
\[
\ph_{\te} (u_{\te}) = v_{\te} \andeqn
\ph_{\te} (v_{\te}) = e^{- \pi i \te} u_{\te}^* v_{\te}.
\]
\end{dfn}

One readily checks that $\ph_{\te}$ really does define an automorphism
of order $6$.
This is the automorphism of order $6$ discussed in the introduction
to Section~\ref{Sec:NCFT}.

\begin{rmk}\label{Ord3Auto}
A computation easily shows that $\ph_{\te}^2$ is the
automorphism of $A_{\te}$ order $3$ discussed in the introduction
to Section~\ref{Sec:NCFT}.
\end{rmk}

The following result is the analog of Theorem~\ref{FieldOfRotAlg}.
We need the slight modification here because $\ph_{\te}$ does not
come from an automorphism of the discrete Heisenberg group,
as a result of the factor $\exp (- \pi i \te)$ in the definition of
$\ph (v_{\te})$.

\begin{lem}\label{FieldOfRotAlg2}
Let $H$ be the group on generators $x$, $y$, and $z$,
subject to the relations
\[
y x = z^2 x y, \,\,\,\,\,\, z x = x z, \andeqn z y = y z.
\]
Let $B = C^* (H)$.
Then there is a \ct\  field of \ca s over $S^1$ whose fiber
over $\exp (2 \pi i \te)$ is $A_{2 \te}$, whose \ca\  of \ct\  sections
is $B$, and such that the
evaluation map $\ev_{\te} \colon B \to A_{2 \te}$ of sections
at $\exp (2 \pi i \te)$ is determined by
\[
\ev_{\te} (x) = u_{2 \te}, \,\,\,\,\,\, \ev_{\te} (y) = v_{2 \te},
\andeqn \ev_{\te} (z) = \exp (4 \pi i \te) \cdot 1.
\]
Moreover, there is a unique automorphism $\ph$ of $B$ of order~$6$
such that
\[
\ph (x) = y, \,\,\,\,\,\, \ph (y) = z^* x^* y, \andeqn \ph (z) = z,
\]
and $\ev_{\te} \circ \ph = \ph_{2 \te} \circ \ev_{\te}$ for all $\te$.
\end{lem}

\begin{proof}
Apply Corollary~3.6 of~\cite{RfF}, with the group being $\Z$,
to the continuous field over $S^1$
whose fiber over $\exp (2 \pi i \te)$ is $C (S^1)$ for all $\te$,
with standard generator $u_{2 \te}$,
and with the action on that fiber being generated by the
automorphism $\af_{\te} (u_{2 \te} ) = \exp (4 \pi i \te) u_{2 \te}$.
The section algebra of this continuous field is $C (S^1 \times S^1)$,
generated by two sections: the constant section $x$
whose value at $\exp (2 \pi i \te)$ is the standard generator
of $C (S^1)$ (called $u_{2 \te}$ above),
and the section $z (\exp (2 \pi i \te)) = \exp (2 \pi i \te) \cdot 1$
for all $\te$.
It is the universal \ca\  on unitary generators
$x$ and $z$ satisfying $z x = x z$.
The action of $\Z$ on the sections is generated by the automorphism
$\af (z) = z$ and $\af (x) = z^2 x$.
Therefore the crossed product is $C^* (H)$, with $y$ being the
implementing unitary for the action.
Corollary~3.6 of~\cite{RfF} thus shows that there is a Hilbert-\ct\  %
field with the desired properties, and the remark in
Definition~3.3 of~\cite{RfF} shows that there is a \ct\  field.

One checks that the formula for $\ph$ defines a unique \hm\  by
checking the commutation relations for the images of the generators.
A computation shows that $\ph^6 = \id_B$, so $\ph$ is an automorphism.
It is immediate to check that
$\ev_{\te} \circ \ph = \ph_{2 \te} \circ \ev_{\te}$.
\end{proof}

\begin{rmk}\label{CompOfContSec}
By comparing the actions, one immediately sees that a section $b$ of the
continuous field of Lemma~\ref{FieldOfRotAlg2} over a small
\nbhd\  of $0$ is \ct\  \ifo\  the section $\zt \mapsto b (\zt^2)$
is a \ct\  section of the
continuous field of Theorem~\ref{FieldOfRotAlg}.
(In fact, the continuous field of Lemma~\ref{FieldOfRotAlg2}
is really just
the pullback of that of Theorem~\ref{FieldOfRotAlg} by the map
$\zt \mapsto \zt^2$.)
\end{rmk}

\begin{lem}\label{Ord6Lem}
Let $\ep > 0$.
Then there is $\dt > 0$ such that
there are \ct\  sections $b_0, b_1, \dots, b_5$ of the
continuous field of Lemma~\ref{FieldOfRotAlg2} defined on the
arc $I = \exp (2 \pi i (- \dt, \, \dt) )$ and satisfying:
\begin{itemize}
\item[(1)]
$0 \leq b_j (\te) \leq 1$ for $0 \leq j \leq 5$ and $| \te | < \dt$.
\item[(2)]
$\| b_j (\te) b_k (\te) \| < \ep$ for $j \neq k$ and $| \te | < \dt$.
\item[(3)]
Setting $b_6 = b_0$,
we have $\| \ph_{\te} (b_j (\te)) - b_{j + 1} (\te) \| < \ep$
for $0 \leq j \leq 5$ and $| \te | < \dt$.
\item[(4)]
With $\ta_{\te}$ being the standard trace on $A_{\te}$ as in
Lemma~\ref{ContOfTrace}, we have the estimate
$\ta_{2 \te} \left( 1 - \sum_{j = 0}^5 b_j (\te) \right) < \ep$
for $| \te | < \dt$.
\end{itemize}
\end{lem}

\begin{proof}
Set $\ep_0 = \frac{1}{12} \ep$.

Define $h_0 \colon \R^2 \to \R^2$
by $h_0 (r_1, r_2) = (r_1 - r_2, \, r_1)$.
Then $h_0$ is an invertible linear map whose matrix is in
${\mathrm{GL}}_2 (\Z)$, so $h_0 (\Z^2) = \Z^2$ and
$h_0$ induces a \hme\  $h \colon \R^2 / \Z^2 \to \R^2 / \Z^2$.
With the obvious identification of $A_0$ with $C (\R^2 / \Z^2)$,
one checks that $\ph_0 (f) = f \circ h^{-1}$ for all
$f \in C (\R^2 / \Z^2)$.

Define open sets $E_0, E_1, \dots, E_5 \S \R^2$ as follows.
Let ${\mathrm{Conv}} (S)$
denote the convex hull of a subset $S \S \R^2$.
Then set
\begin{align*}
E_0 & = \sint \left( {\mathrm{Conv}} \left( \left\{ (0, 0), \, (1, 0),
    \, \left( \ts{\frac{2}{3}}, \ts{\frac{1}{3}} \right) \right\}
        \right) \right), \\
E_1 & = \sint \left( {\mathrm{Conv}} \left( \left\{ (0, 0), \, (1, 1),
    \, \left( \ts{\frac{1}{3}}, \ts{\frac{2}{3}} \right) \right\}
        \right) \right), \\
E_2 & = \sint \left( {\mathrm{Conv}} \left( \left\{ (1, 0), \, (1, 1),
    \, \left( \ts{\frac{2}{3}}, \ts{\frac{1}{3}} \right) \right\}
        \right) \right), \\
E_3 & = \sint \left( {\mathrm{Conv}} \left( \left\{ (0, 1), \, (1, 1),
    \, \left( \ts{\frac{1}{3}}, \ts{\frac{2}{3}} \right) \right\}
        \right) \right), \\
E_4 & = \sint \left( {\mathrm{Conv}} \left( \left\{ (0, 0), \, (1, 1),
    \, \left( \ts{\frac{2}{3}}, \ts{\frac{1}{3}} \right) \right\}
        \right) \right), \\
E_5 & = \sint \left( {\mathrm{Conv}} \left( \left\{ (0, 0), \, (0, 1),
    \, \left( \ts{\frac{1}{3}}, \ts{\frac{1}{2}} \right) \right\}
        \right) \right).
\end{align*}
These sets can be described as follows.
Divide $[0, 1]^2$ in half along the main diagonal $r_1 = r_2$.
Divide each of the two resulting triangles using line segments
from the point $\left( \ts{\frac{1}{3}}, \ts{\frac{2}{3}} \right)$
or $\left( \ts{\frac{2}{3}}, \ts{\frac{1}{3}} \right)$ as appropriate
to the three vertices.
The $E_j$ for even $j$ are in the lower right triangle.

These $6$ sets are disjoint open subsets of $[0, 1]^2$ whose union is
dense in $[0, 1]^2$ and each of which has measure $\frac{1}{6}$.
Calculations show that
\[
h_0 (E_0) = E_1, \,\,\,\,\,\, h_0 (E_1) = E_2 + (-1, \, 0),
\,\,\,\,\,\, h_0 (E_2) = E_3,
\]
\[
h_0 (E_3) = E_4 + (-1, \, 0),
\,\,\,\,\,\, h_0 (E_4) = E_5, \andeqn h_0 (E_5) = E_0 + (-1, \, 0).
\]
Since $(-1, \, 0) \in \Z^2$, the images $U_j$ of the $E_j$ in
$\R^2 / \Z^2$ are disjoint open sets with measure $\frac{1}{6}$
which are cyclically permuted by $h$.

Choose a closed set $T_0 \S U_0$ with measure greater than
$\frac{1}{6} - \ep_0$, and choose $f_0 \in C (\R^2 / \Z^2)$
with $0 \leq g \leq 1$, with $\supp (g) \S U_0$, and such that
$g_0 (x) = 1$ for all $x \in T_0$.
Define $g_j = g_0 \circ h^{-j}$ for $1 \leq j \leq 5$.
Since the sets $U_j$ are pairwise disjoint, we have
$g_j g_k = 0$ for $j \neq k$.
Moreover,
\[
\ta_0 \left( 1 - \ssum{j = 0}{5} g_j \right) < 1 - 6 \ep_0.
\]
The proof is now finished as in the last paragraph of the proof
of Lemma~\ref{NCFT1}.
We use Lemma~\ref{FieldOfRotAlg2}
in place of Theorem~\ref{FieldOfRotAlg}.
Lemma~\ref{ContOfTrace} still implies that traces of
\ct\  sections are \ct, by Remark~\ref{CompOfContSec}.
\end{proof}

\begin{prp}\label{Ord6HasARP}
Let $\te \in \R \SM \Q$.
Then the action of $\Zqs$ on $A_{\te}$ generated by the
automorphism $\ph_{\te}$ of Definition~\ref{Ord6Auto} has the \aRp.
\end{prp}

\begin{proof}
Except for one point (see below),
the proof is essentially the same as for Proposition~\ref{FTHasARP}.
We use Lemma~\ref{FieldOfRotAlg2}
in place of Theorem~\ref{FieldOfRotAlg}.
We use Lemma~\ref{Ord6Lem} in place of Lemma~\ref{NCFT1},
adjusting for the fact that $B_{\te} \cong A_{2 \te}$
rather than $A_{\te}$.

The difference is in the proof that
$\ps \circ \ph_{\et} = \ph_{\te} \circ \ps$.
With $\dt_1$ and $\dt_2$ chosen to satisfy conditions
analogous to those in the proof of Proposition~\ref{FTHasARP},
we require $n \in \N$ with $n > 0$ and $m \in \Z$ such that
\[
| \exp (2 \pi i n \te) - 1 | < \dt_2
\andeqn | n^2 \te - m | < \dt_1,
\]
and in addition such that $m$ is even.
To get this, apply Lemma~\ref{Dioph} with $\frac{1}{2} \te$ in place
of $\te$,
obtaining $n \in \N$ with $n > 0$ and $m_0 \in \Z$ such that
\[
\left| \exp (\pi i n \te) - 1 \right| < \ts{ \frac{1}{2} } \dt_2
\andeqn
\left| \ts{ \frac{1}{2} } n^2 \te - m_0 \right|
           < \ts{ \frac{1}{2} } \dt_1.
\]
Then take $m = 2 m_0$.

This done, as in the proof of Proposition~\ref{FTHasARP},
set $\et = n^2 \te - m$
and let $\ps \colon A_{\et} \to A_{\te}$ be the \hm\  such that
$\ps (u_{\et}) = u_{\te}^n$ and $\ps (v_{\et}) = v_{\te}^n$.
We must check that $\ps \circ \ph_{\et} = \ph_{\te} \circ \ps$.
First, use the relation
$v_{\te} u_{\te}^* = \exp (- 2 \pi i \te) u_{\te}^* v_{\te}$
to get
$u_{\te}^{-n} v_{\te}^n
   = \exp (\pi i n (n - 1) \te) (u_{\te}^* v_{\te})^n$.
Further, since $m$ is even, we get
$\exp (- \pi i \et) = \exp (- \pi i n^2 \te)$.
Now
\[
\ps \circ \ph_{\et} (u_{\et}) = \ps ( v_{\et} )
  = v_{\te}^n
  = \ph_{\te} (u_{\te}^n)
  = \ph_{\te} \circ \ps (u_{\et})
\]
and
\begin{align*}
\ps \circ \ph_{\et} (v_{\et})
 & = \ps \left( e^{- \pi i \et} u_{\et}^* v_{\et} \right)
   = e^{- \pi i \et} u_{\te}^{-n} v_{\te}^n
   = e^{- \pi i n^2 \te} e^{ \pi i n (n - 1) \te}
            (u_{\te}^* v_{\te})^n   \\
 & = \left( e^{- \pi i \te} u_{\te}^* v_{\te} \right)^n
   = \ph_{\te} (v_{\te})^n 
   = \ph_{\te} \circ \ps (v_{\et}).
\end{align*}
This completes the proof.
\end{proof}

Rather than repeating everything for the action of $\Zqh$,
we prove the following lemma.

\begin{lem}\label{TRPForSubgp}
Let $A$ be a stably finite \suca, and let $\af \in \Aut (A)$ generate
an action of $\Zqn$ which has the \aRp.
Let $n = l m$ be a nontrivial factorization in positive integers.
Then $\af^l$ generates an action of $\Zq{m}$ which has the \aRp.
\end{lem}

\begin{proof}
Let $F \S A$ be a finite set,
let $\ep > 0$, let $N \in \N$,
and let $x \in A$ be a nonzero positive element.
Apply Lemma~\ref{StARPDfn} with $\ep / l^2$ in place of $\ep$,
and with the given values of $F$, $N$, and $x$.
Call the resulting \pj s $f_0, f_1, \dots, f_{n - 1}$.
For $0 \leq k \leq m - 1$, define
\[
e_k = \sum_{j = k l}^{k l + l - 1} f_j.
\]
Then
\begin{align*}
\| \af^l (e_k) - e_{k + 1} \|
  & \leq \sum_{j = k l}^{k l + l - 1} \| \af^l (f_j) - f_{j + l} \|  \\
  & \leq \sum_{j = k l}^{k l + l - 1} \sum_{r = 0}^{l - 1}
                   \| \af (f_{j + r}) - f_{j + r + 1} \|
    < l^2 (\ep / l^2) = \ep.
\end{align*}
For $a \in F$, we have
\[
\| a e_k - e_k a \|
  \leq \sum_{j = k l}^{k l + l - 1} \| a f_j - f_j a \|
  < l (\ep / l^2) < \ep.
\]
This verifies Conditions~(1) and~(2) of Definition~\ref{ARPDfn}.
Conditions~(3) and~(4) are immediate since
$\sum_{k = 0}^{l - 1} e_k = \sum_{j = 0}^{n - 1} f_j$, and
since every $e_k$ dominates some $f_j$.
\end{proof}

\begin{cor}\label{Ord3HasARP}
Let $\te \in \R \SM \Q$, and let
$\ph_{\te}$ be as in Definition~\ref{Ord6Auto}.
Then the action of $\Zqh$ on $A_{\te}$ generated by
$\ph_{\te}^2$ has the \aRp.
\end{cor}

\begin{proof}
Combine Proposition~\ref{Ord6HasARP} and Lemma~\ref{TRPForSubgp}.
\end{proof}

\begin{lem}\label{CPbyOrd6hasUCT}
Let $\te \in \R \SM \Q$, and let $\ph_{\te}$ be as in
Definition~\ref{Ord6Auto}.
Then the crossed products $C^* (\Zqs, \, A_{\te}, \, \ph_{\te} )$
and $C^* (\Zqh, \, A_{\te}, \, \ph_{\te}^2 )$ satisfy the
Universal Coefficient Theorem.
\end{lem}

\begin{proof}
The proof is the same as for Lemma~\ref{CPbyFThasUCT}, once we
know that the group $H$ in Lemma~\ref{FieldOfRotAlg2} is amenable.
This is immediate from the fact that the subgroup generated by
$x$, $y$, and $z^2$ is an amenable subgroup (in fact, the discrete
Heisenberg group) of index~$2$.
\end{proof}

\begin{thm}\label{CPbyOrd6isAH}
Let $\te \in \R \SM \Q$, and let $\ph_{\te}$ be as in
Definition~\ref{Ord6Auto}.
Then the crossed products $C^* (\Zqs, \, A_{\te}, \, \ph_{\te} )$
and $C^* (\Zqh, \, A_{\te}, \, \ph_{\te}^2 )$
are simple AH algebras with slow dimension growth and real rank zero.
\end{thm}

\begin{proof}
The proof is the same as for Theorem~\ref{CPbyFTisAH}.
\end{proof}

We don't know whether any of these algebras are AF, because
we don't know the K-theory.
However, the K-theory computations for $\te$ rational,
done in~\cite{FW5} and~\cite{FW6}, suggest that all should be AF.
These computations are for the fixed point algebras rather than
the crossed products, but one expects the K-theory for the
irrational fixed point algebras to be the same as for the generic
rational ones, and in the irrational case
Proposition~\ref{Ord6HasARP}, Corollary~\ref{Ord3HasARP},
Corollary~\ref{CrPrIsSimple}, and Proposition~\ref{FPIsSimple}
show that the fixed point algebra
is always Morita equivalent to the crossed product.

\section{The crossed product of a simple noncommutative
            torus by the flip}\label{Sec:CrPrdFlip}

\indent
In this section, we prove that the crossed product of any
higher dimensional simple noncommutative torus by the flip
is a simple AF algebra.
This generalizes Theorem~3.1 of~\cite{Bc}, and completely answers
a question raised in the introduction to~\cite{FW2}.
As in Section~\ref{Sec:NCFT}, there are three parts:
the proof that the action has the \aRp,
the proof that the crossed product \suct,
and the computation of the K-theory of the crossed product.

\begin{dfn}\label{FlipDfn}
Let $\te$ be a skew symmetric real $d \times d$ matrix, and let
$A_{\te}$ be the corresponding noncommutative torus with
unitary generators
\[
u_1 (\te), \, u_2 (\te), \, \dots, \, u_d (\te)
\]
being the unitaries $u_1, u_2, \dots, u_d$ of Notation~\ref{StdNtn}.
When no confusion is likely to result,
we will suppress $\te$ in the notation.
Then the {\emph{flip automorphism}}
$\sm_{\te} \colon A_{\te} \to A_{\te}$
is the automorphism of order two determined by
$\sm_{\te} (u_k) = u_k^*$ for $1 \leq k \leq n$.
We further set, in this section,
$B_{\te} = C^* (\Zqt, \, A_{\te}, \, \sm_{\te})$,
and call it the crossed product by the flip.
\end{dfn}

\begin{ntn}\label{RestField}
For $d \geq 1$ let $R_d$ be the set of all
skew symmetric real $d \times d$ matrices $\te$ such that
$\te_{j, k} \in [-1, \, 1]$ for all $j$ and $k$.
For any $d$, we let $r \colon R_{d + 1} \to R_d$ be the restriction map
$r ( \te) = \te |_{\Z^{d} \times \{ 0 \} }$,
which deletes the last row and column of $\te$.
(See Remark~\ref{Restrict}.
The intended value of $d$ will always be clear.)
For any continuous field $C$ over $R_{d - 1}$,
we let $r^* (C)$ denote the continuous field over $R_d$
obtained as the pullback of $C$.
(See~\cite{Dx}, especially Sections~10.1 and~10.3,
for information on \ct\  fields of \ca s.
See Lemma~1.3 of~\cite{KP2} for pullbacks of continuous fields.)
In a variation of Notation~\ref{SectionNot},
we denote the set of all \ct\  sections of $C$ by $\Gm (C)$.
\end{ntn}

\begin{lem}\label{RC}
Let the notation be as in Notation~\ref{RestField}.
\begin{itemize}
\item[(1)]
The space $R_d$ is homeomorphic to $[-1, \, 1]^{d (d - 1)/2}$.
In particular, $R_1 = \{ 0 \}$.
\item[(2)]
For any continuous field $C$ over $R_{d}$,
there is an isomorphism
$\ph \colon \Gm (r^* (C)) \to C ([0, 1]^d, \, \Gm (C))$
such that every evaluation map
$\ev_{\te} \colon \Gm (r^* (C)) \to C_{r (\te)} = r^* (C)_{\te}$
factors as
\[
\Gm (r^* (C)) \stackrel{\ph}{\longrightarrow}
 C ([0, 1]^d, \, \Gm (C)) \stackrel{\ev_t}{\longrightarrow}
 \Gm (C) \stackrel{\ev_{r (\te)}}{\longrightarrow}
 C_{r (\te)}
\]
for some $t \in [0, 1]^d$.
\item[(3)]
For $\te \in R_d$ and a continuous field $C$ over $R_{d}$,
the map $\ev_{\te} \colon \Gm (r^* (C)) \to C_{r (\te)}$
is the composition of $\ev_{r (\te)} \colon \Gm (C) \to C_{r (\te)}$
and a homotopy equivalence.
\end{itemize}
\end{lem}

\begin{proof}
Part~(1) is immediate.
For Part~(2), we use the homeomorphism
from $R_{d + 1}$ to $[0, 1]^d \times R_d$ given by
\[
\rh \mapsto ( \rh_{d + 1, \, 1}, \, \rh_{d + 1, \, 2},
               \, \dots, \, \rh_{d + 1, \, d}, \, r (\rh))
\]
to construct $\ph$.
Part~(3) then follows from the fact that the \hm\  $\ev_t$
in Part~(2) is a homotopy equivalence.
\end{proof}

\begin{lem}\label{HDContField}
Let the notation be as in
Definition~\ref{FlipDfn} and Notation~\ref{RestField}.
\begin{itemize}
\item[(1)]
For every skew symmetric real $d \times d$ matrix $\te$ there is
$\te' \in R_d$ such that $A_{\te'} \cong A_{\te}$
and $B_{\te'} \cong B_{\te}$.
\item[(2)]
There exists a \ct\  field $A^{(d)}$ of \ca s over $R_d$,
whose fiber over
$\te \in R_d$ is $A_{\te}$, and such that the sections
\[
\te \mapsto
  f (\te) u_1 (\te)^{n_1} u_2 (\te)^{n_2} \cdots u_d (\te)^{n_d},
\]
for $f \in C (R_d)$ and $n_1, n_2, \dots, n_d \in \Z$,
span a dense subalgebra of the \ca\  of all \ct\  sections.
\item[(3)]
There is an automorphism $\sm$ of the \ca\  $\Gm (A^{(d)})$
of \ct\  sections
of the \ct\  field of Part~(3) such that for any section $a$, we
have $\sm (a) (\te) = \sm_{\te} (a (\te))$ for all $\te$.
Moreover, there exists a \ct\  field $B^{(d)}$ of \ca s over $R_d$,
whose fiber over
$\te \in R_d$ is $B_{\te}$, and whose \ca\  of \ct\  sections
is $C^* ( \Zqt, \, \Gm (A^{(d)}), \, \sm)$.
\end{itemize}
\end{lem}

\begin{proof}
Part~(1) is immediate from the fact that $A_{\te}$ is unchanged
if some entry $\te_{j, k}$ is replaced by
$\te_{j, k} + n$ for some $n \in \Z$.
Part~(2) is contained in the discussion following
Corollary~2.9 of~\cite{RfF}, but we will need to reprove it
anyway in the course of proving Part~(3).

To prove Part~(3), we prove by induction on $d$ that both
$A^{(d)}$ and $B^{(d)}$ are Hilbert continuous fields in the
sense of Definition~3.3 of~\cite{RfF}.
The comment in the last part of that definition shows that
such a field is automatically continuous in the usual sense,
so this does in fact prove Part~(3) (as well as Part~(2)).

If $d = 1$ then $R_d$ consists of just one point,
and the statement is trivial.
Suppose now that the result is known for $d$.
Since $A^{(d)}$ is a Hilbert continuous field,
it is immediate from Lemma~\ref{RC}(2)
that $r^* (A^{(d)})$ is as well.
Let $\Z$ act on the fiber $r^* (A^{(d)})_{\te} = A_{r (\te)}$
according to the
automorphism of Lemma~\ref{ItCrPrd}, so that the crossed
product is $A_{\te}$.
{}From the formula given in the proof, it is clear that
we get a continuous field of actions of $\Z$ in the sense of
Definition~3.1 of~\cite{RfF}.
So Theorem~3.5 of~\cite{RfF} shows that
$A^{(d + 1)}$ is Hilbert \ct\  and has the right algebra of
sections.
Moreover, the actions $\sm_{\te}$ define a
continuous field of actions of $\Zqt$,
so the same theorem now yields Hilbert continuity of $B^{(d + 1)}$.
\end{proof}

We now prove that the flip has the \aRp.
We use the method of proof of Proposition~\ref{IrratAPR}.
There is one difficulty:
the restriction of the flip to the irrational rotation subalgebra
$A_{\et} \S A_{\te}$ that we construct
need not be the flip on that algebra, but rather has the form
\[
v_{\et} \mapsto \ld v_{\et}^* \andeqn w_{\et} \mapsto \zt w_{\et}^*
\]
for some a priori unknown $\ld, \, \zt \in S^1$.
Fortunately, all such automorphisms are conjugate to the flip,
and the set of possibilities is compact.

\begin{prp}\label{FlipHasARP}
Let $\te \in R_d$ be nondegenerate (Definition~\ref{NondegDfn}).
Then the action of $\Zqt$ on $A_{\te}$ generated by the
automorphism $\sm_{\te}$ of Definition~\ref{FlipDfn} has the \aRp.
\end{prp}

\begin{proof}
We verify the conditions of Theorem~\ref{ARPFromPosElts2}.
Theorem~\ref{NCTIsAT} implies that $A_{\te}$ has tracial rank zero,
and Theorem~\ref{Simplicity} implies that $A_{\te}$ has
a unique tracial state.

For the rest of the proof,
let the notation for rotation algebras be as in
Theorem~\ref{FieldOfRotAlg}, Notation~\ref{SectionNot},
and Lemma~\ref{ContOfTrace}.
In particular, $A$ is the \ca\  of the discrete Heisenberg group.

Let $\ep > 0$ and let $F \S A_{\te}$ be a finite set.
Choose and fix $\et_0 \in \R \SM \Q$.

For $\et \in \R$ and $\ld, \, \zt \in S^1$
let $\sm_{\et}^{(\ld, \zt)} \in \Aut (A_{\et})$
be the automorphism of order two determined by
\[
\sm_{\et}^{(\ld, \zt)} ( v_{\et} ) = \ld v_{\et}^*
\andeqn
\sm_{\et}^{(\ld, \zt)} ( w_{\et} ) = \zt w_{\et}^*.
\]
Thus $\sm_{\et}^{(1, 1)} = \sm_{\et}$.
Moreover,
we claim that $\sm_{\et}^{(\ld, \zt)}$ is conjugate to $\sm_{\et}$.
To see this, choose $\ld_0, \, \zt_0 \in S^1$ such that
$\ld_0^2 = \ld$ and $\zt_0^2 = \zt$,
let $\ph \in \Aut (A_{\et})$ be the automorphism determined by
\[
\ph ( v_{\et} ) = \ld_0 v_{\et}
\andeqn
\ph ( w_{\et} ) = \zt_0 w_{\et},
\]
and check that $\sm_{\et} \circ \ph = \ph \circ \sm_{\et}^{(\ld, \zt)}$.

It follows from Proposition~\ref{FTHasARP} and Lemma~\ref{TRPForSubgp}
that $\sm_{\et}$ has the \aRp\  for $\et \not\in \Q$,
so that each $\sm_{\et}^{(\ld, \zt)}$ also has the \aRp.
Taking $\et = \et_0$,
for each $\ld, \, \zt \in S^1$ choose orthogonal \pj s
$p_0^{(\ld, \zt)}, \, p_1^{(\ld, \zt)} \in A_{\et_0}$ such that
\[
\left\| \rsz{ \sm_{\et_0}^{(\ld, \zt)}
        \left( \rsz{ p_0^{(\ld, \zt)} } \right)
                           - p_1^{(\ld, \zt)} } \right\|
      < \ts{ \frac{1}{2} } \ep
\andeqn
\ta_{\et_0}
       \left( \rsz{1 - p_0^{(\ld, \zt)} - p_1^{(\ld, \zt)} } \right)
    < \ts{ \frac{1}{2} } \ep.
\]
Choose \ct\  selfadjoint sections
$b_0^{(\ld, \zt)}, \, b_1^{(\ld, \zt)} \in A$
with $0 \leq b_j^{(\ld, \zt)} \leq 1$ and such that
$b_j^{(\ld, \zt)} (\et_0) = p_j^{(\ld, \zt)}$.

We claim that for fixed $\ld_0, \, \zt_0 \in S^1$ the function
on $S^1 \times S^1 \times \R$ defined by
\[
(\ld, \zt, \et) \mapsto
  \left\| \rsz{ \sm_{\et}^{(\ld, \zt)}
                \left( \rsz{ b_0^{(\ld_0, \zt_0)} (\et) } \right)
         - b_1^{(\ld_0, \zt_0)} (\et) } \right\|
\]
is jointly \ct.
Identify the \cfn\  algebra $C (S^1 \times S^1, \, A)$
with the universal \ca\  generated by unitaries $v, w, x, y, z$
subject to the relations
that $w v = z v w$ and $x, y, z$ are all central,
by taking $A = C^* (v, w, z)$ and identifying $x$ and $y$ with the
functions $(\ld, \zt) \mapsto \ld$ and $(\ld, \zt) \mapsto \zt$.
Using Theorem~\ref{FieldOfRotAlg}, this algebra can be realized
in an obvious way as the section algebra of a \ct\  field over
$S^1 \times S^1 \times S^1$ whose fiber over
$(\ld, \, \zt, \, \exp (2 \pi i \et) )$ is $A_{\et}$.
Let $\sm \in \Aut (A)$ be determined by
\[
\sm (v) = x v^*, \,\,\,\,\,\,
\sm (w) = y w^*, \,\,\,\,\,\,
\sm (x) = x, \,\,\,\,\,\,
\sm (y) = y, \andeqn
\sm (z) = z.
\]
The evaluation map
$\ev_{\ld, \, \zt, \, \exp (2 \pi i \et)} \colon
       C (S^1 \times S^1, \, A) \to A_{\et}$
satisfies
\[
\sm_{\et}^{(\ld, \zt)} \circ \ev_{\ld, \, \zt, \, \exp (2 \pi i \et)}
  = \ev_{\ld, \, \zt, \, \exp (2 \pi i \et)} \circ \sm,
\]
and
since $\sm$ sends \ct\  sections to \ct\  sections the claim follows.

Using the claim, for every $\ld_0, \, \zt_0 \in S^1$ there are
open \nbhd s $U$ of $(\ld_0, \zt_0)$ and $V$ of $\et_0$ such that
\[
\left\| \rsz{ \sm_{\et}^{(\ld, \zt)}
    \left( \rsz{ b_0^{(\ld_0, \zt_0)} (\et) } \right)
         - b_1^{(\ld_0, \zt_0)} (\et) } \right\|
   < \ep
\]
for $(\ld, \zt) \in U$ and $\et \in V$.
Using compactness of $S^1 \times S^1$ to cover it with finitely
many of the sets $U$, and intersecting the corresponding sets $V$,
we find a \nbhd\  $W_0$ of $\et_0$ and finitely many pairs of sections
chosen from the ones above,
which for simplicity we call
$b_0^{(m)}, \, b_1^{(m)}$ for $1 \leq m \leq r$,
such that for every $(\ld, \zt) \in S^1$ there is $m$ such that
for every $\et \in W_0$ we have
\[
\left\| \rsz{ \sm_{\et}^{(\ld, \zt)}
               \left( \rsz{ b_0^{(m)} (\et) } \right)
         - b_1^{(m)} (\et) } \right\| < \ep.
\]
Since $b_0^{(m)} (\et_0) b_1^{(m)} (\et_0) = 0$ and
$\ta_{\et_0}
       \left( \rsz{ 1 - b_0^{(m)} (\et_0) - b_1^{(m)} (\et_0) } \right)
     < \ts{ \frac{1}{2} } \ep$,
we can use the continuous field structure and Lemma~\ref{ContOfTrace}
to find an open set $W$ with $\et_0 \in W \S W_0$ such that
for every $(\ld, \zt) \in S^1 \times S^1$ there is $m$ such that
for every $\et \in W$ we have:
\begin{itemize}
\item[(1)]
$0 \leq b_0^{(m)} (\et), \, b_1^{(m)} (\et) \leq 1$.
\item[(2)]
$\left\| \rsz{ b_0^{(m)} (\et) b_1^{(m)} (\et) } \right\| < \ep
    \rule{0em}{2.5ex}$.
\item[(3)]
$\left\| \rsz{ \sm_{\et}^{(\ld, \zt)}
            \left( \rsz{ b_0^{(m)} (\et) } \right)
         - b_1^{(m)} (\et) } \right\| < \ep
    \rule{0em}{2.5ex}$.
\item[(4)]
$\left| \rsz{ \ta_{\et} \left( \rsz{ 1 - b_0^{(m)} (\et)
            - b_1^{(m)} (\et) } \right) } \right| < \ep
    \rule{0em}{2.5ex}$.
\end{itemize}

Apply Lemma~\ref{MakeSubalgStep3}
with $\te$, $\et_0$, $F$, and $\ep$ as given,
with $k = n = 1$, with $U = \exp (2 \pi i W)$,
and with
$S = \left\{ \rsz{ b_0^{(m)} |_U, \, b_1^{(m)} |_U
               \colon 1 \leq m \leq r } \right\}$.
We obtain $\et \in \R \SM \Q$ and
\[
l = (l_1, l_2, \dots, l_d) \in \Z^d \andeqn
m = (m_1, m_2, \dots, m_d) \in \Z^d.
\]
Set
\[
x = u_1^{l_1} u_2^{l_2} \cdots u_d^{l_d} \andeqn
y = u_1^{m_1} u_2^{m_2} \cdots u_d^{m_d},
\]
so that $y x = \exp ( 2 \pi i \et) x y$.
Let $\ph \colon A_{\et} \to A_{\te}$ be the \hm\  such that
$\ph (v_{\et}) = x$ and $\ph (w_{\et}) = y$.
We get
\[
\left\| \left[ \rsz{ \ph \left( \rsz{ b_j^{(m)} (\et) } \right),
                           \, c } \right] \right\| < \ep
\]
for all $c \in F$, for $1 \leq m \leq r$, and for $j = 0, \, 1$.

Using the commutation relations in Notation~\ref{StdNtn},
there are $\ld, \zt \in S^1$ such that
\[
x^* = {\overline{\ld}} (u_1^*)^{l_1} (u_2^*)^{l_2} \cdots (u_d^*)^{l_d}
 \andeqn
y^* = {\overline{\zt}} (u_1^*)^{m_1} (u_2^*)^{m_2} \cdots (u_d^*)^{m_d}.
\]
It follows that
$\sm_{\te} \circ \ph = \ph \circ \sm_{\et}^{(\ld, \zt)}$.
By the above, there exists $m$ such that conditions (1)--(4) above
hold with this $\ld, \zt, \et$.
Set $a_j = \ph \left( \rsz{ b_j^{(m)} (\et) } \right)$
for $j = 0, \, 1$.
We have $0 \leq a_0, \, a_1 \leq 1$,
and we claim that $a_0$ and $a_1$ satisfy the hypotheses (1)--(4)
of Theorem~\ref{ARPFromPosElts2}.
Parts~(1) and~(4) follow from Conditions~(2) and~(4) above,
Part~(2) follows from Condition~(3) above and the relation
$\sm_{\te} \circ \ph = \ph \circ \sm_{\et}^{(\ld, \zt)}$,
and Part~(3) follows from the commutator estimate at the
end of the previous paragraph.

It now follows from Theorem~\ref{ARPFromPosElts2}
that the action of $\Zqt$ on $A_{\te}$ generated by $\sm_{\te}$
has the \aRp.
\end{proof}

Now we turn to K-theory.
In~\cite{Kj} and~\cite{FW2}, the K-theory calculations used the
exact sequence of~\cite{Nt}.
However, for the Universal Coefficient Theorem we need information
about KK-theory, so we use~\cite{Pm}.

\begin{ntn}\label{Dih}
Throughout this section, we let $G$ be the infinite dihedral group
$G = \Z \rtimes \Zqt$.
We let $g_0$ be the generator of $\Z$ corresponding to $1$,
and we let $h_0$ be the nontrivial element of $\Zqt$,
so that the defining relations are
$h_0^2 = 1$ and $h_0 g_0 h_0^{-1} = g_0^{-1}$.
We further let $H_0 = \{ 1, h_0 \}$ be the obvious copy of $\Zqt$.
Note that also $h_1 = g_0 h_0$ satisfies $h_1^2 = 1$,
and set $H_1 = \{ 1, h_1 \}$, another subgroup isomorphic to $\Zqt$.
It is known that $G$ is the free product of $H_0$ and $H_1$;
see for example the discussion of the structure of this group in the
preliminaries in~\cite{Kj}.
\end{ntn}

\begin{lem}\label{StrOfAction}
With notation as in Notation~\ref{Dih} and Notation~\ref{RestField},
there is a continuous field $\af$
of actions of $G$ on $r^* (A^{(d - 1)})$
(in the sense of Definition~3.1 of~\cite{Rf1}),
determined by actions $\af^{(\te)}$ of $G$ on $r^* (A^{(d - 1)})_{\te}$,
such that:
\begin{itemize}
\item[(1)]
$\af^{(\te)}_{g_0}$ acts on the generators of
$r^* (A^{(d - 1)})_{\te} = A_{ r (\te) }$
by the formula
$\af^{(\te)}_{g_0} (u_j) = \exp (2 \pi i \te_{j, d}) u_j$
for $1 \leq j \leq d - 1$,
and the continuous field of crossed products
$\te
 \mapsto C^* (\Z, \, r^* (A^{(d - 1)})_{\te}, \, \af^{(\te)}_{g_0} )$
is $A^{(d)}$.
\item[(2)]
$\af^{(\te)}_{h_0}$ acts on the generators of
$r^* (A^{(d - 1)})_{\te} = A_{r (\te)}$
by $\af^{(\te)}_{g_0} (u_j) = u_j^*$
for $1 \leq j \leq d - 1$,
and the continuous field
$\te \mapsto
    C^* (\Zqt, \, r^* (A^{(d - 1)})_{\te}, \, \af^{(\te)}_{h_0} )$
of crossed products is $r^* (B^{(d - 1)})$ with $B^{(d - 1)}$
as in Lemma~\ref{HDContField}(3).
\item[(3)]
$\af_{h_1}$ is conjugate to $\af_{h_0}$ via the continuous family of
automorphisms $\te \mapsto \ph^{(\te)}$ determined by
$\ph^{(\te)} (u_j) = \exp (- \pi i \af_{j, d}) u_j$
for $1 \leq j \leq d - 1$.
\item[(4)]
The continuous field of crossed products
$\te \mapsto C^* (G, \, r^* (A^{(d - 1)})_{\te}, \, \af^{(\te)} )$
is $B^{(d)}$.
\end{itemize}
\end{lem}

\begin{proof}
Lemma~\ref{ItCrPrd} shows
that there exist automorphisms $\bt^{(\te)}$ on the individual fibers
with the properties claimed for $\af^{(\te)}_{g_0}$ in Part~(1).
Let $S$ be the set of all \ct\  sections of $\Gm (r^* (A^{(d - 1)}))$
of the form
\[
\te \mapsto f (\te) u_1 (\te)^{n_1} u_2 (\te)^{n_2}
                          \cdots u_{d - 1} (\te)^{n_{d - 1}}
\]
with $f \in C (R_d)$.
It is easy to check that if $s \in S$ then
$\te \mapsto \bt^{(\te)} (s (\te))$ is again in $S$, hence \ct.
Since ${\mathrm{Span}} (S)$ is dense in $\Gm (r^* (A^{(d - 1)}))$,
this implies that $\te \mapsto \bt^{(\te)}$ defines an
automorphism of $r^* (A^{(d - 1)})$.

Applying the pullback via $r$ to the automorphism of
Lemma~\ref{HDContField}(3),
we obtain
an automorphism $\te \mapsto \gm^{(\te)}$ of $r^* (A^{(d - 1)})$
with the properties claimed in Part~(2).
A computation shows that
\[
{\ts{ \left( \gm^{(\te)} \right)^2 }} = \id_{r^* (A^{(d - 1)})_{\te} }
\andeqn
\gm^{(\te)} \circ \bt^{(\te)}
            \circ {\ts{ \left( \gm^{(\te)} \right)^{-1} }}
  = {\ts{ \left( \bt^{(\te)} \right)^{-1} }}
\]
for all $\te$, from which it follows that there is an action of $G$
satisfying~(1) and~(2).

To get Part~(4), we rewrite
\begin{align*}
C^* (G, \, r^* (A^{(d - 1)})_{\te}, \, \af^{(\te)})
&  = C^* (\Z \rtimes \Zqt, \, r^* (A^{(d - 1)})_{\te}, \, \af^{(\te)})
                     \\
&  = C^* (\Zqt, \,
       C^* (\Z, \, r^* (A^{(d - 1)})_{\te}, \, \af^{(\te)}_{g_0}),
                     \, \ov{\gm}^{(\te)}),
\end{align*}
where $\ov{\gm}^{(\te)}$ acts on $r^* (A^{(d - 1)})_{\te}$ via
$\af^{(\te)}_{h_0}$ and on the canonical unitary $u$ in the
crossed product by $\Z$ like conjugation by $h_0$ on $g_0$, that is,
$\ov{\gm}^{(\te)} (u) = u^{-1}$.
Identifying $u_j (r (\te))$ with $u_j (\te)$ and $u$ with $u_d (\te)$,
we identify $\ov{\gm}^{(\te)}$ with
$\sm_{\te}$ as in Lemma~\ref{HDContField}(3).

It remains to proves Part~(3).
A computation shows that
\[
\ph^{(\te)} \circ \af^{(\te)}_{g_0 h_0}
   = \af^{(\te)}_{h_0} \circ \ph^{(\te)}
\]
for all $\te$.
\end{proof}

\begin{lem}\label{KKExactSeqs}
Let $d \in \N$ and let $\te \in R_d$.
Then for every separable \ca\  $D$ there is a commutative
diagram with exact columns:
\[
\begin{CD}
@VVV   @VVV\\
KK^{1 - j} (D, \, \Gm (B^{(d)} ))
                  @>{(\ev_{\te})_*}>> KK^{1 - j} (D, \, B_{\te} )   \\
@VVV   @VVV\\
KK^{j} (D, \, \Gm (r^* (A^{(d - 1)} )))
                  @>{(\ev_{\te})_*}>> KK^{j} (D, \, A_{r (\te)} )   \\
@VVV   @VVV\\
\begin{array}{l} KK^{j} (D, \, \Gm (r^* (B^{(d - 1)} ))) \\
   \hspace*{2em} \mbox{} \oplus KK^{j} (D, \, \Gm (r^* (B^{(d - 1)} )))
  \end{array}
     @>{(\ev_{\te})_*}>>
        \begin{array}{l} KK^{j} (D, \, B_{r (\te)} ) \\
             \hspace*{2em} \mbox{} \oplus KK^{j} (D, \, B_{r (\te)} )
                              \end{array} \\
@VVV   @VVV\\
KK^{j} (D, \, \Gm (B^{(d)}) )
                  @>{(\ev_{\te})_*}>> KK^{j} (D, \, B_{\te} )   \\
@VVV   @VVV\\
KK^{1 - j} (D, \, \Gm (r^* (A^{(d - 1)} )))
              @>{(\ev_{\te})_*}>> KK^{1 - j} (D, \, A_{r (\te)} )   \\
@VVV   @VVV
\end{CD}
\]
Moreover, for every separable nuclear \ca\  $D$ there is a commutative
diagram with exact columns:
\[
\begin{CD}
@VVV   @VVV\\
KK^{1 - j} (A_{r (\te)}, \, D )
    @>{(\ev_{\te})^*}>> KK^{1 - j} (\Gm (r^* (A^{(d - 1)} )), \, D) \\
@VVV   @VVV\\
KK^{j} (B_{\te}, \, D )
                @>{(\ev_{\te})^*}>> KK^{j} (\Gm (B^{(d)}), \, D )   \\
@VVV   @VVV\\
\begin{array}{l} KK^{j} (B_{r (\te)}, \, D) \\
     \hspace*{2em} \mbox{} \oplus KK^{j} (B_{r (\te)}, \, D)
  \end{array}
      @>{(\ev_{\te})^*}>>
         \begin{array}{l} KK^{j} (\Gm (r^* (B^{(d - 1)})), \, D)  \\
              \hspace*{2em} \mbox{}
                    \oplus KK^{j} (\Gm (r^* (B^{(d - 1)})), \, D)
                              \end{array} \\
@VVV   @VVV\\
KK^{j} (A_{r (\te)}, \, D)
       @>{(\ev_{\te})^*}>> KK^{j} (\Gm (r^* (A^{(d - 1)})), \, D)   \\
@VVV   @VVV\\
KK^{1 - j} (B_{\te}, \, D)
             @>{(\ev_{\te})^*}>> KK^{1 - j} (\Gm (B^{(d)}), \, D)   \\
@VVV   @VVV
\end{CD}
\]
\end{lem}

\begin{proof}
We derive these exact sequences from~\cite{Pm},
so we begin by constructing a suitable action of $G$
on a suitable tree $X$.
We take the vertices to be $X^0 = \Z$,
and the oriented edges to be
\[
X^1 = \{ (2 n - 1, \, 2 n), \, (2 n + 1, \, 2 n) \colon n \in \Z \}.
\]
We define the actions of $g_0$ and $h_0$ to be
\[
g_0 \cdot n = n + 2 \andeqn h_0 \cdot n = - n
\]
on $X^0$, and
\[
g_0 \cdot (2 n - 1, \, 2 n) = (2 n + 1, \, 2 n + 2), \,\,\,\,\,\,
g_0 \cdot (2 n + 1, \, 2 n) = (2 n + 3, \, 2 n + 2)
\]
and
\[
h_0 \cdot (2 n - 1, \, 2 n) = (- 2 n + 1, \, - 2 n), \,\,\,\,\,\,
h_0 \cdot (2 n + 1, \, 2 n) = (- 2 n - 1, \, - 2 n)
\]
on $X^1$.
One proves that these define an action by checking that the maps on
$X^0$ and $X^1$ are bijective and satisfy the two defining
relations of Notation~\ref{Dih}.
Following the discussion after Lemma~3 of~\cite{Pm},
we choose subsets $\Sm^0 \S X^0$ and $\Sm^1 \S X^1$ which are
representatives of the quotients of $X^0$ and $X^1$ by $G$.
We choose $\Sm^0 = \{ 0, 1 \}$ and $\Sm^1 = \{ (1, 0) \}$.
The stabilizers are then
\[
G_0 = \{ 1, h_0 \}, \,\,\,\,\,\, G_1 = \{ 1, \, g_0 h_0 \},
\andeqn G_{(1, 0)} = \{ 1 \}.
\]
We now apply the exact sequence of Theorem~16 of~\cite{Pm}
for the full crossed products, to the actions of $G$ on both
$\Gm (r^* (A^{(d - 1)} ) )$ and on $A_{r (\te)}$.
By naturality we obtain a commutative diagram with exact columns,
in which for simplicity
we suppress the notation for the actions in the crossed products.
\[
\begin{CD}
@VVV   @VVV\\
KK^{1 - j} (D, \, C^* (G, \, \Gm (r^* (A^{(d - 1)} ) ) ) )
      @>{(\ev_{\te})_*}>>
            KK^{1 - j} (D, \, C^* (G, \, A_{r (\te)} ) )  \\
@VVV   @VVV\\
KK^{j} (D, \, \Gm (r^* (A^{(d - 1)} )))
                  @>{(\ev_{\te})_*}>> KK^{j} (D, \, A_{r (\te)} )   \\
@VVV   @VVV\\
\begin{array}{l} KK^{j} (D, \, C^* (G_0, \, \Gm (r^* (A^{(d - 1)} ) ) ) ) \\
     \hspace*{1.3em} \mbox{} \oplus
            KK^{j} (D, \, C^* (G_1, \, \Gm (r^* (A^{(d - 1)} ) ) ) )
           \end{array}
     @>{(\ev_{\te})_*}>>
        \begin{array}{l} KK^{j} (D, \, C^* (G_0, \, A_{r (\te)})) \\
             \hspace*{1.3em} \mbox{} \oplus
                        KK^{j} (D, \, C^* (G_1, \, A_{r (\te)}))
           \end{array} \\
@VVV   @VVV\\
KK^{j} (D, \, C^* (G, \, \Gm (r^* (A^{(d - 1)} ) ) ) )
        @>{(\ev_{\te})_*}>> KK^{j} (D, \, C^* (G, \, A_{r (\te)} ) )   \\
@VVV   @VVV\\
KK^{1 - j} (D, \, \Gm (r^* (A^{(d - 1)} )))
              @>{(\ev_{\te})_*}>> KK^{1 - j} (D, \, A_{r (\te)} )   \\
@VVV   @VVV
\end{CD}
\]
Then we use Lemma~\ref{StrOfAction} to identify
$C^* (G, \, \Gm (r^* (A^{(d - 1)} ) ) ) \cong \Gm (B^{(d)} )$,
etc.
Note that Lemma~\ref{StrOfAction}(3) ensures that the crossed
products by $G_1$ are all isomorphic to the corresponding crossed
products by $G_0$.
This gives the first diagram of the present lemma.

To get the second one, we use Theorem~17 of~\cite{Pm}
in place of Theorem~16 of~\cite{Pm}, and proceed the same way.
\end{proof}

\begin{prp}\label{CptKTh}
For every skew symmetric real $d \times d$ matrix $\te$,
we have $K_0 (B_{\te}) \cong \Z^{3 \cdot 2^{d - 1}}$
and $K_1 (B_{\te}) = 0$.
\end{prp}

\begin{proof}
By Lemma~\ref{HDContField}(1), it suffices to consider $\te \in R_d$.
We prove by induction on $d$ that the maps
\[
(\ev_{\te})_* \colon K_* (\Gm (A^{(d)} )) \to K_* (A_{\te})
\andeqn
(\ev_{\te})_* \colon K_* (\Gm (B^{(d)} )) \to K_* (B_{\te})
\]
are isomorphisms for all $\te$.
The proof of Theorem~7 of~\cite{FW2} shows that $K_* (B_{\te})$ is as
claimed in the statement for at least some values of $\te \in R_d$,
so it will follow that this is correct for all $\te \in R_d$.

When $d = 1$, the space $R_d$ consists of a single point.
So the only possible $\ev_{\te}$ is an isomorphism of \ca s.
Suppose the result is known for all $\te \in R_d$ for some $d$.
Let $\te \in R_{d + 1}$.
In the first diagram in Lemma~\ref{KKExactSeqs},
replace $d$ by $d + 1$, and use this $\te$.
Also take $D = \C$, giving a diagram in ordinary K-theory.
The induction hypothesis and Lemma~\ref{RC}(3) imply that
\[
\ev_{\te} \colon \Gm (r^* (A^{(d)} )) \to A_{r (\te)}
\andeqn
\ev_{\te} \colon \Gm (r^* (B^{(d)} )) \to  B_{r (\te)}
\]
are isomorphisms on K-theory.
Now
$(\ev_{\te})_* \colon K_* (\Gm (B^{(d + 1)} )) \to K_* (B_{\te})$
is an isomorphism by the Five Lemma.
That $(\ev_{\te})_* \colon K_* (\Gm (A^{(d + 1)} )) \to K_* (A_{\te})$
is an isomorphism follows in the same way
from a similar diagram in which the columns
are obtained from the Pimsner-Voiculescu exact sequence for the
K-theory of crossed products by $\Z$~\cite{PV}.
This completes the induction, and the proof.
\end{proof}

\begin{prp}\label{UCTForFlip}
For every skew symmetric real $d \times d$ matrix $\te$,
the algebra $B_{\te}$ satisfies the \uct.
(See Definition~\ref{DL1} and Remark~\ref{DL2}.)
\end{prp}

\begin{proof}
By Lemma~\ref{HDContField}(1), it suffices to consider $\te \in R_d$.
We prove this by induction on $d$.
When $d = 1$, the only possible algebra is
$B_0 = \Csw{2}{C (S^1)}{\sm_0}$,
which satisfies the \uct\  because it is type~I.

Suppose the result is known for all $\te \in R_d$ for some $d$.
Let $\te \in R_{d + 1}$.
Let $D$ be any separable nuclear \ca\  %
such that $K_* (D)$ is an injective abelian group.
Using the diagrams of Lemma~\ref{KKExactSeqs},
we can construct the commutative diagram:
\[
\begin{CD}
@VVV   @VVV\\
KK^{1 - j} (A_{r (\te)}, \, D )
   @>{(\ev_{\te})^*}>>  \Hom (K_* (A_{r (\te)}), \, K_* (D) )_{1 - j} \\
@VVV   @VVV\\
KK^{j} (B_{\te}, \, D )
        @>{(\ev_{\te})^*}>> \Hom (K_* (B_{\te}), \, K_* (D) )_j   \\
@VVV   @VVV\\
\begin{array}{l} KK^{j} (B_{r (\te)}, \, D) \\
     \hspace*{2em} \mbox{} \oplus KK^{j} (B_{r (\te)}, \, D)
  \end{array}
      @>{(\ev_{\te})^*}>>
         \begin{array}{l} \Hom (K_* (B_{r (\te)}), \, K_* (D) )_j  \\
              \hspace*{2em} \mbox{}
                    \oplus \Hom (K_* (B_{r (\te)}), \, K_* (D) )_j
                              \end{array} \\
@VVV   @VVV\\
KK^{j} (A_{r (\te)}, \, D)
       @>{(\ev_{\te})^*}>> \Hom (K_* (A_{r (\te)}), \, K_* (D) )_j   \\
@VVV   @VVV\\
KK^{1 - j} (B_{\te}, \, D)
    @>{(\ev_{\te})^*}>>  \Hom (K_* (B_{\te}), \, K_* (D) )_{1 - j}   \\
@VVV   @VVV
\end{CD}
\]
The left column is exact by Lemma~\ref{KKExactSeqs},
and the right column is exact because $K_* (D)$ is injective
and by Lemma~\ref{KKExactSeqs}.
Using Theorem~1.17 of~\cite{RSUCT} (see the preceding discussion
for the definition of ${\mathcal{N}}$),
it follows from Theorem~\ref{NCTIsAT}
that $A_{r (\te)}$ satisfies the Universal Coefficient Theorem.
Also, $B_{r (\te)}$ satisfies the Universal Coefficient Theorem
by the induction hypothesis.
So $B_{\te}$ satisfies the Universal Coefficient Theorem
by the Five Lemma.
\end{proof}

\begin{thm}\label{FlipAF}
Let $\te$ be a nondegenerate skew symmetric real $d \times d$ matrix,
with $d \geq 2$.
Then $\Cs{2}{A_{\te}}{\sm_{\te}}$ is a simple AF~algebra.
\end{thm}

\begin{proof}
Combining Proposition~\ref{FlipHasARP} and Theorem~\ref{RokhTAF},
we find that the crossed product
$B_{\te} = \Cs{2}{A_{\te}}{\sm_{\te}}$ has tracial rank zero.
By Proposition~\ref{UCTForFlip}
it satisfies the Universal Coefficient Theorem,
and it is clearly separable and nuclear.
Proposition~\ref{CptKTh} shows that $K_0 (B_{\te})$ is torsion free
and $K_1 (B_{\te}) = 0$.
The result therefore follows from Lemma~\ref{ConseqOfClass}.
\end{proof}

\section{More on the Rokhlin property and tracial approximate
   innerness}\label{Sec:More}

\indent
In this section we prove, in analogy with Theorem~\ref{RokhTAF},
that the crossed product of an AF~algebra by an action with the
\sRp\  is again an AF~algebra.
We also prove that tracially approximately inner automorphisms
act trivially on the trace space and on $K_0$ mod infinitesimals.
These results are relevant to the examples we construct in
the next section.

It is useful to start the first proof with the following analog of
Lemma~\ref{PreserveFD}.

\begin{lem}\label{PreserveFD2}
Let $A$ be a unital \ca,
and let $\af \in \Aut (A)$ be an automorphism which
satisfies $\af^n = \id_A$ and such that the action of $\Zqn$
generated by $\af$ has the strict Rokhlin property
(Definition~\ref{ERPDfn}).
Then for every finite set
$F \S A$, every \fd\  subalgebra $E \S A$, and every $\ep > 0$,
there are \mops\  $e_0, e_1, \dots, e_{n - 1} \in A$
and a unitary $v \in A$ such that:
\begin{itemize}
\item[(1)]
$\| \af (e_j) - e_{j + 1} \| < \ep$ for $0 \leq j \leq n - 1$,
where, following Convention~\ref{ARPConv},  we take $e_n = e_0$.
\item[(2)]
$\| e_j a - a e_j \| < \ep$ for $0 \leq j \leq n - 1$ and all $a \in F$.
\item[(3)]
$\sum_{j = 0}^{n - 1} e_j = 1$.
\item[(4)]
$\| v - 1 \| < \ep$, and $e_j$ commutes with all elements of
$v E v^*$ for $0 \leq j \leq n - 1$.
\end{itemize}
\end{lem}

\begin{proof}
The proof is a slightly simpler version of the proof of
Lemma~\ref{PreserveFD}.
\end{proof}

\begin{thm}\label{SRokhAF}
Let $A$ be a unital AF~algebra,
and let $\af \in \Aut (A)$ be an automorphism which
satisfies $\af^n = \id_A$ and such that the action of $\Zqn$
generated by $\af$ has the strict Rokhlin property.
Then $\CZnAa$ is an AF algebra.
\end{thm}

\begin{proof}
We prove that for every finite set $S \S \CZnAa$
and every $\ep > 0$,
there is a \fd\  C*-subalgebra $D \S \CZnAa$
such that every element of $S$ is within $\ep$ of an element of $D$.
Theorem~2.2 of~\cite{Brt} will then imply that
$\CZnAa$ is AF.
Let $u \in \CZnAa$
be the canonical unitary implementing the automorphism $\af$.
\Wolog\  we may assume that $S = F \cup \{ u \}$ for
a finite subset $F$ of the unit ball of $A$.

Set
\[
\ep_0 = \frac{\ep}{12 (n + 1)^5}.
\]
Choose $\dt > 0$ with $\dt < \ep_0$,
and so small that whenever $e$ and $f$ are \pj s in a \ca\  $C$
such that $\| e - f \| < \dt$, then there is a partial isometry
$s \in C$ such that
\[
s s^* = e, \,\,\,\,\,\, s^* s = f, \andeqn
\| s - e \| < \ep_0.
\]

Since $A$ is AF, there is a \fd\  C*-subalgebra $E_0 \S A$
such that for every $a$ in the finite set
\[
S_0 = F \cup \af (F) \cup \cdots \cup \af^{n - 1} (F),
\]
there exists $b \in E_0$ such that $\| a - b \| < \dt$.
Apply Lemma~\ref{PreserveFD2} with $S_0$ in place of $F$,
with $E_0$ in place of the \fd\  subalgebra $E$,
and with $\dt$ in place of $\ep$.
We obtain a unitary
$y \in A$ and \mops\  $e_0, \, e_1, \, \dots, \, e_{n - 1} \in A$ which
commute with all elements of $y E_0 y^*$,
such that $\| e_j a - a e_j \| < \dt$ for all $a \in S_0$,
such that $\sum_{j = 0}^{n - 1} e_j = 1$,
such that $\| \af (e_j) - e_{j + 1} \| < \dt$,
and such that $\| y - 1 \| < \dt$.

According to the choice of $\dt$, for $1 \leq j \leq n - 1$
there are partial isometries $w_j \in A$ such that
\[
w_j w_j^* = e_j, \,\,\,\,\,\, w_j^* w_j = \af ( e_{j - 1} ), \andeqn
\| w_j - e_j \| < \ep_0.
\]
Apply Lemma~\ref{P} to the $e_j$ and $w_j$,
with $\ep_0$ in place of $\ep$.
We obtain a unitary $w$ as there such that $\| w - 1 \| < 2 n^2 \ep_0$
and satisfying the condition of Lemma~\ref{PertOfCrPrd},
such that moreover the automorphism $\bt = \Ad (w) \circ \af$
satisfies $\bt^n = \id_A$,
\[
\| \bt^k (a) - \af^k (a) \| \leq 4 k n^2 \ep_0 \| a \|
\]
for $a \in A$, and $\bt (e_j) = e_{j + 1}$ for all $j$.

Define
\[
E = \bigoplus_{k = 0}^{n - 1} \bt^k (e_0 y E_0 y^* e_0)
  = \bigoplus_{k = 0}^{n - 1} e_k \bt^k (y E_0 y^*) e_k.
\]
Since the $e_k$ are orthogonal,
$\sum_{k = 0}^{n - 1} e_k = 1$,
and $e_0$ commutes with every element of $y E_0 y^*$,
it follows that $E$ is a $\bt$-invariant \fd\  subalgebra of $A$
such that $1_A \in E$.

Let $a \in F$.
We estimate the distance from $a$ to $E$.
We begin by estimating,
using $\af^{n - k} (a) \in S_0$ at the fourth step,
\begin{align*}
\| [ e_k, a] \|
 & = \| \bt^{n - k} ([ e_k, a]) \|
   = \| [ e_0, \, \bt^{n - k} (a) ] \|  \\
 & \leq 2 \| \bt^{n - k} (a) - \af^{n - k} (a) \|
        + \| [ e_0, \, \af^{n - k} (a) ] \|  \\
 & < 8 (n - k) n^2 \ep_0 + \dt
   < (8 n^3 + 1) \ep_0.
\end{align*}
It follows that if $k \neq l$ then
$\| e_k a e_l \| = \| [e_k, a] e_l \| < (8 n^3 + 1) \ep_0$.
Since there are fewer than $n^2$ terms in the sum
in the second expression, we can estimate
\[
\left\| a - \ssum{k = 0}{n - 1} e_k a e_k \right\|
  \leq \sum_{k = 0}^{n - 1} \sum_{l \neq k} \| e_k a e_l \|
  < n^2 (8 n^3 + 1) \ep_0.
\]
Moreover, by construction there exists $b_k \in E_0$ such that
$\| b_k - \af^{n - k} (a) \| < \dt$.
Then
\[
\| b_k - \bt^{n - k} (a) \|
  < \dt + 4 (n - k) n^2 \ep_0 < (4 n^3 + 1) \ep_0,
\]
whence
\[
\| \bt^k (e_0 b_k e_0) - e_k a e_k \|
 = \| e_0 b_k e_0 - e_0 \bt^{n - k} (a) e_0 \| < (4 n^3 + 1) \ep_0.
\]
It follows that $c = \sum_{k = 0}^{n - 1} \bt^k ( e_0 b_k e_0) \in E$
and satisfies
\begin{align*}
\| a - c \|
& \leq \left\| a - \ssum{k = 0}{n - 1} e_k a e_k \right\|
       + \sum_{k = 0}^{n - 1} \| e_k a e_k - \bt^k (e_0 b_k e_0) \|  \\
& < n^2 (8 n^3 + 1) \ep_0 + n (4 n^3 + 1) \ep_0
   \leq 12 (n + 1)^5 \ep_0 < \ep.
\end{align*}

Let
\[
D_0 = C^* (\Zqn, \, E, \, \bt |_E),
\]
which is a \fd\  subalgebra of $\Cs{n}{A}{\bt}$.
Let
\[
\ph \colon
 \Cs{n}{A}{\bt} \to \CZnAa
\]
be the isomorphism of Lemma~\ref{PertOfCrPrd}.
We take the \fd\  subalgebra $D$ to be $D = \ph (D_0)$.
Because $\ph (a) = a$ for $a \in A$,
we have shown that for every $a \in F$
there is $c \in E \S D_0$ such that $\| c - a \| < \ep$.
Let $v \in \Cs{n}{A}{\bt}$ be the
canonical unitary implementing the automorphism $\bt$.
Then $\ph (v) \in D$ and Lemma~\ref{PertOfCrPrd} gives
$\| \ph (v) - u \| = \| w - 1 \| < 2 n^2 \ep_0 < \ep$.
\end{proof}

The following two lemmas give essentially the only restrictions we know
on the behavior of tracially approximately inner automorphisms:
they act as the identity on $K_0$ mod infinitesimals
and on the tracial states.
It is also easy to see that if
$\af \in \Aut (A)$ is tracially approximately inner
but not approximately inner, then there must be
arbitrarily small positive elements $\et \in K_0 (A)$
such that $\af_* (\et) = \et$.
The dual automorphisms of the actions in Examples~\ref{CAR2}
and~\ref{CAR4} are tracially approximately inner,
but the one in Example~\ref{CAR2} is not the identity on $K_0$
and the one in Example~\ref{CAR4} is not the identity on $K_1$.

Recall that an element $\et$ of a partially ordered group $G$
with order unit $u \in G_+ \SM \{ 0 \}$ is
{\emph{infinitesimal}} if $- m u \leq n \et \leq m u$ for all
$m, \, n \in \N$ with $m > 0$.
See Definition~1.10 of~\cite{GPS}, where this definition is given for
simple dimension groups.
Clearly we need only consider $m = 1$.
By Proposition~4.7 of~\cite{Gd0}, an equivalent condition is that all
states on $(G, u)$ vanish on $\et$.

\begin{lem}\label{TAIAndInf}
Let $A$ be a stably finite \suca,
and let $\af \in \Aut (A)$ be tracially approximately inner.
Then $\af_* (\et) - \et$ is infinitesimal for every $\et \in K_0 (A)$.
\end{lem}

\begin{proof}
We prove that for every $\et \in K_0 (A)$
we have $- [1_A] \leq \af_* (\et) - \et \leq [1_A]$.
This implies the result, because replacing $\et$ by $n \et$
gives $- [1_A] \leq n [ \af_* (\et) - \et ] \leq [1_A]$.

Accordingly, let $\et \in K_0 (A)$,
and choose $n \in \N$ and \pj s $p, \, r \in M_n (A)$
such that $\et = [p] - [r]$.
Let $\ch \colon \R \setminus \left\{ \frac{1}{2} \right\} \to \R$
be the characteristic function of $\left( \frac{1}{2}, \I \right)$.
Choose $\ep > 0$ so small that $(n^2 + 2) \ep < \frac{1}{6}$, and
also so small that whenever $C$ is a unital \ca\  and
$f, \, p \in C$ are \pj s such that $\| f q - q f \| < (n^2 + 2) \ep$,
then $\frac{1}{2}$ is not in the spectrum of either
$f q f$ or $(1 - f) q (1 - f)$, and moreover
the \pj s $q_0 = \ch (f q f)$ and $q_1 = \ch ((1 - f) q (1 - f))$
satisfy
$\| q_0 + q_1 - q \| < \frac{1}{6}$.

Apply Definition~\ref{TAInnDfn} with
$F = \{ p_{j, k}, \, r_{j, k} \colon 1 \leq j, \, k \leq n \}$,
the set of all matrix entries of $p$ and $r$,
with $\ep$ as just chosen, with $N = 2 n$,
and with $x = 1$.
Let $e \in A$ and $v \in e A e$ be the resulting \pj\  and unitary.

We have
\[
\| (1 \otimes e) p - p (1 \otimes e) \|
  \leq \sum_{j, \, k = 1}^n \| e p_{j, k} - p_{j, k} e \| < n^2 \ep.
\]
By the choice of $\ep$ the \pj s
\[
p_0 = \ch ((1 \otimes e) p (1 \otimes e) ) \in M_n (e A e)
\]
and
\[
p_1 = \ch ((1 - 1 \otimes e) p (1 - 1 \otimes e) )
  \in M_n ((1 - e) A (1 - e))
\]
are defined and satisfy $\| p_0 + p_1 - p \| < \frac{1}{6}$.
To get a similar result for $(\id \otimes \af) (p)$,
we begin by observing that
\[
\| (\id \otimes \af)^{-1} (1 \otimes e) - 1 \otimes e \|
  = \| e - \af (e) \| < \ep.
\]
So
\begin{align*}
& \| (1 \otimes e) [(\id \otimes \af) (p)]
         - [(\id \otimes \af) (p)] (1 \otimes e) \|   \\
& \hspace*{4em}
  = \| [(\id \otimes \af)^{-1} (1 \otimes e)] p
         - p [(\id \otimes \af)^{-1} (1 \otimes e)] \|  \\
& \hspace*{4em}
  \leq 2 \| (\id \otimes \af)^{-1} (1 \otimes e) - 1 \otimes e \|
      + \| (1 \otimes e) p - p (1 \otimes e) \|    \\
& \hspace*{4em}
  < 2 \ep + n^2 \ep = (n^2 + 2) \ep.
\end{align*}
Therefore we get \pj s $q_0  \in M_n (e A e)$
and $q_1 \in M_n ((1 - e) A (1 - e))$
such that $\| q_0 + q_1 - (\id \otimes \af) (p) \| < \frac{1}{6}$.
Note in particular that
\[
[p] = [p_0] + [p_1] \andeqn \af_* ([p]) = [q_0] + [q_1]
\]
in $K_0 (A)$.

We now claim that $[p_0] = [q_0]$ in $K_0 (A)$.
First,
\[
\| (1 \otimes e) p (1 \otimes e) - p_0 \|
   = \| (1 \otimes e) [p - (p_0 + p_1)] (1 \otimes e) \|
   < {\textstyle{ \frac{1}{6} }}.
\]
Next,
\begin{align*}
& \| (\id \otimes \af) ((1 \otimes e) p (1 \otimes e)) - q_0 \| \\
& \hspace*{4em}
   \leq 2 \| (\id \otimes \af) (1 \otimes e) - 1 \otimes e \|
     + \| (1 \otimes e) (\id \otimes \af) (p) (1 \otimes e) - q_0 \|  \\
& \hspace*{4em}
   < 2 \ep + \| (\id \otimes \af) (p) - (q_0 + q_1) \|
   < 2 \ep + {\textstyle{ \frac{1}{6} }}.
\end{align*}
Finally,
\begin{align*}
& \| (1 \otimes v) (1 \otimes e) p (1 \otimes e) (1 \otimes v)^*
           - (\id \otimes \af) ((1 \otimes e) p (1 \otimes e)) \|   \\
& \hspace*{10em}
  \leq \sum_{j, \, k = 1}^n \|v e p_{j, k} e v^* - \af (e p_{j, k} e)\|
  < n^2 \ep.
\end{align*}
Putting these estimates together gives
\[
\| (1 \otimes v) p_0 (1 \otimes v)^* - q_0 \|
     < {\textstyle{ \frac{1}{6} }} + n^2 \ep + 2 \ep
         + {\textstyle{ \frac{1}{6} }}
     < {\textstyle{ \frac{1}{2} }}.
\]
The claim follows.

Repeating the argument of the last two paragraphs with $r$ in place
of $p$, we find \pj s
\[
r_0, \, s_0 \in M_n (e A e) \andeqn
r_1, \, s_1 \in M_n ((1 - e) A (1 - e))
\]
such that
\[
[r] = [r_0] + [r_1], \,\,\,\,\,\, \af_* ([r]) = [s_0] + [s_1],
  \andeqn [r_0] = [s_0]
\]
in $K_0 (A)$.

Now recall from Condition~(5) of Definition~\ref{TAInnDfn} that
there are $2 n$ \mops\  $f_1, f_2, \dots, f_{2 n} \leq e$,
each of which is \mvnt\  to $1 - e$.
Since $p_1, \, q_1, \, r_1, \, s_1 \in M_n ((1 - e) A (1 - e))$,
this implies that
\[
[p_1] + [s_1] \leq [1_A] \andeqn [q_1] + [r_1] \leq [1_A]
\]
in $K_0 (A)$.
Therefore
\[
\af_* (\et) - \et
  = \af_* ([p]) - [p] - \af_* ([r]) + [r]
  = [q_1] - [p_1] - [s_1] + [r_1],
\]
and
\[
- [1_A] \leq - [p_1] - [s_1]
        \leq \af_* (\et) - \et
        \leq [q_1] + [r_1]
        \leq [1_A].
\]
This completes the proof.
\end{proof}

\begin{lem}\label{TAIAndTraces}
Let $A$ be a \suca, and let $\af \in \Aut (A)$ be
tracially approximately inner.
Then $\ta \circ \af = \ta$ for every tracial state $\ta$ on $A$.
\end{lem}

\begin{proof}
Let $a \in A$ and let $\ep > 0$.
We show that $| \ta ( \af (a)) - \ta (a) | < \ep$.
\Wolog\  $a \geq 0$ and $\| a \| \leq 1$.
Choose $N \in \N$ such that $\frac{1}{N} < \frac{1}{3} \ep$.
Apply Definition~\ref{TAInnDfn} with $F = \{ a \}$,
with $\frac{1}{3} \ep$ in place of $\ep$, with $N$ as chosen here,
and with $x = 1$.
Let $e \in A$ and $v \in e A e$ be the resulting \pj\  and unitary.
Condition~(5) of Definition~\ref{TAInnDfn} implies that
$\ta (1 - e) \leq \frac{1}{N} < \frac{1}{3} \ep$.
We have
\[
\ta (a) = \ta (e a e) + \ta ( (1 - e) a (1 - e) ),
\]
with
\[
0 \leq \ta ( (1 - e) a (1 - e) ) \leq \ta (1 - e)
 < {\textstyle{ \frac{1}{3} }} \ep,
\]
so $| \ta (a) - \ta (e a e) | < {\textstyle{ \frac{1}{3} }} \ep$.
Similarly $\ta ( \af (1 - e) ) \leq \frac{1}{N} < \frac{1}{3} \ep$
and $| \ta ( \af (a)) - \ta ( \af (e a e) ) | < \frac{1}{3} \ep$.
Finally, $\| v e a e v^* - \af (e a e) \| < \frac{1}{3} \ep$
and $\ta (v e a e v^*) = \ta (e a e)$.
Putting these together gives
\begin{align*}
| \ta ( \af (a)) - \ta (a) |
 & \leq | \ta ( \af (a)) - \ta ( \af (e a e) ) |
       + | \ta (\af (e a e) ) - \ta (e a e) |
       + | \ta (e a e) - \ta (a) |   \\
 & < {\textstyle{ \frac{1}{3} }} \ep + {\textstyle{ \frac{1}{3} }} \ep
       + {\textstyle{ \frac{1}{3} }} \ep
   = \ep,
\end{align*}
as desired.
\end{proof}

\section{Examples}\label{Sec:Exs}

\indent
We look at several examples of actions on the $2^{\infty}$ UHF algebra.
They demonstrate the following:
\begin{itemize}
\item
Even on a UHF algebra, an action with the \aRp\  need not have
the \sRp, and in fact the crossed product by such an action need
not be AF.
See Example~\ref{CAR4}.
\item
If an automorphism $\af$ of a \suca\  $A$ with tracial rank zero
is approximately inner
and generates an action of $\Zqn$  with the \aRp,
it does not follow that the automorphisms of the dual action are
approximately inner---even if $A$ is UHF and $\af$ is the pointwise
limit of inner automorphisms $\Ad (u_k)$ with $u_k^n = 1$.
See Example~\ref{CAR2}.
\item
If an automorphism $\af$ of a \suca\  $A$ with tracial rank zero
is tracially approximately inner
and generates an action of $\Zqn$  with the \sRp,
it does not follow that the automorphisms of the dual action are
approximately inner---even if $A$ is AF.
Use the dual of the action in Example~\ref{CAR4}.
\item
A tracially approximately inner automorphism of a \suca\  $A$
with tracial rank zero
need not be trivial on $K_0 (A)$---even if $A$ is AF and
$\af$ generates an action of $\Zqn$  with the \sRp.
Use the dual of the action in Example~\ref{CAR2}.
\item
A tracially approximately inner automorphism of a \suca\  $A$
with tracial rank zero
need not be trivial on $K_1 (A)$---even if $A$ is AT and
$\af$ generates an action of $\Zqn$  with the \sRp.
Use the dual of the action in Example~\ref{CAR4}.
\item
There is an automorphism $\af$ of a simple AF algebra such that
$\af^n = \id_A$, such that $\CZnAa$ is
again a simple AF algebra, but such that this action does not
have the \aRp.
This can happen even when the dual action has the \sRp\  and
$\af$ is approximately inner.
See Example~\ref{CAR3}.
\item
The dual of an action with the \aRp\  need not have the \aRp,
even when both the original algebra and the crossed product are
simple AF algebras,
and even when the original action also has the strict Rokhlin property.
Use the dual of the action in Example~\ref{CAR3}.
\end{itemize}

The first three examples are product type actions
of $\Zqt$ on the $2^{\infty}$ UHF algebra $D$.
In each case we represent $D$ as an infinite tensor product,
and the automorphism generating the action as an infinite
tensor product of inner automorphisms.
It is useful to give a general lemma.

\begin{lem}\label{Z2Struct}
Let $D$ be an infinite tensor product \ca\  and let $\af \in \Aut (D)$
be an automorphism of order two, of the form
\[
D = \bigotimes_{n = 1}^{\I} M_{k (n)} \andeqn
\af = \bigotimes_{n = 1}^{\I}
 \Ad ( e_n - f_n ),
\]
with $k (n) \in \N$ and where $e_n, \, f_n \in M_{k (n)}$
are \pj s with $e_n + f_n = 1$.
Let $D_n = \bigotimes_{m = 1}^{n} M_{k (m)}$,
so that $D = \dirlim D_n$,
and write $D_n = M_{t (n)}$, where $t (n) = \prod_{m = 1}^{n} k (m)$.
Then the direct system of crossed
products can be identified as
\[
C^* (\Zqt, D_n) \cong M_{t (n)} \oplus M_{t (n)},
\]
with the dual action given by the flip $\sm_n (a, b) = (b, a)$
for $a, \, b \in M_{t (n)}$,
and where the maps
\[
\ps_n \colon M_{t (n - 1)} \oplus M_{t (n - 1)}
         \to M_{t (n)} \oplus M_{t (n)}
\]
are given by
\[
\ps_n (a, b)
 = (a \otimes e_n + b \otimes f_n, \, b \otimes e_n + a \otimes f_n)
\]
for $a, \, b \in M_{t (n - 1)}$.
Assuming moreover that $e_n, \, f_n \neq 0$ for all $n$, we have:
\begin{itemize}
\item[(1)]
$\Cs{2}{D}{\af}$ is a simple unital AF~algebra.
\item[(2)]
The action of $\Zqt$ generated by $\af$
is approximately representable
in the sense of Definition~3.6(2) of~\cite{Iz}.
\item[(3)]
The dual action on $\Cs{2}{D}{\af}$ has
the strict Rokhlin property.
\item[(4)]
The action of $\Zqt$ generated by $\af$ has the \sRp\  \ifo\  %
the the dual action is approximately representable.
\item[(5)]
If the action of $\Zqt$ generated by $\af$ has the \aRp, then
the generating automorphism ${\widehat{\af}}$ of the dual action is
tracially approximately inner.
\end{itemize}
\end{lem}

\begin{proof}
For $n \in \N$ and a unitary $v \in M_n$ with $v^2 = 1$,
we use the isomorphism
$\Csw{2}{M_n}{\Ad (v)} \to M_n \oplus M_n$
which sends $a \in M_n$ to $(a, a)$ and the canonical unitary
of the crossed product to $(v, \, - v)$.
The identification of the direct system of crossed products is then a
calculation, which we omit.
Now assume that $e_n, \, f_n \neq 0$ for all $n$.
Simplicity of $\Cs{2}{D}{\af}$
follows from the fact that the partial embedding multiplicities
in the direct system of crossed products are all nonzero,
and the rest of~(1) is immediate.
Part~(2) is immediate.
Part~(3) now follows from Lemma~3.8(2) of~\cite{Iz}.
Part~(4) is Lemma~3.8(1) of~\cite{Iz}.
To prove Part~(5), and we use Theorem~\ref{DualHasRokh}.
The algebra $D$ has cancellation of \pj s because it is AF,
and $\Cs{2}{D}{\af}$ has
the weak divisibility property (Definition~\ref{WDivDfn})
because it is also AF.
Therefore Theorem~\ref{DualHasRokh} applies.
\end{proof}

Presumably the converse of Part~(5) is true as well,
but we have no need for it.

Our first example is the most regular possible,
and is given to contrast with the remaining ones.

\begin{exa}\label{CAR1}
Let $\af$ be the automorphism of order~$2$ given by
\[
D = \bigotimes_{n = 1}^{\I} M_2 \andeqn
\af = \bigotimes_{n = 1}^{\I}
 \Ad \left( \begin{array}{cc} 1 & 0 \\ 0 &  -1 \end{array} \right).
\]
Then the action of $\Zqt$ generated by $\af$ has the
strict Rokhlin property, the crossed product is again the
$2^{\infty}$ UHF algebra, and the dual action is just another copy
of the given action.
All this is easily proved using Lemma~\ref{Z2Struct},
and is also a special case of Example~3.2 of~\cite{Iz}.
\end{exa}

\begin{exa}\label{CAR2}
Let $\bt$ be the automorphism of order~$2$ given by
\[
D = \bigotimes_{n = 1}^{\I} M_{2^n} \andeqn
\af = \bigotimes_{n = 1}^{\I}
 \Ad ( 1_{2^{n - 1} + 1} \oplus (- 1_{2^{n - 1} - 1}) ).
\]
The automorphism in the $n$-th tensor factor is conjugation by a
diagonal unitary in which $2^{n - 1} + 1$ diagonal entries are equal
to $1$ and $2^{n - 1} - 1$ diagonal entries are equal to $- 1$.
The crossed product $\Cs{2}{D}{\bt}$ is a simple unital AF~algebra
by Lemma~\ref{Z2Struct}(1).

In this case, the action of $\Zqt$ generated by $\bt$ has the
\aRp, but does not have the strict Rokhlin property.
The dual action has the \sRp\  %
and its generator is tracially approximately inner,
but the generator is not approximately inner
and does not induce the identity
map on $K_0 (\Cs{2}{D}{\bt} )$.
We prove this in the next three propositions.
\end{exa}

\begin{prp}\label{CAR2_1}
The action of $\Zqt$ generated by the automorphism $\bt$
of Example~\ref{CAR2} has the \aRp.
\end{prp}

\begin{proof}
We verify the conditions of Theorem~\ref{ARPFromPosElts2}.
(We will not use the full strength of this theorem because the
elements we construct will in fact be \pj s.)
Let $D_n = \bigotimes_{k = 1}^{n} M_{2^k}$, so that
$D_{n} = D_{n - 1} \otimes M_{2^n}$ and
$D = \dirlim D_n$.
Let
\[
v_n = 1_{2^{n - 1} + 1} \oplus (- 1_{2^{n - 1} - 1}) \in M_{2^n}
\andeqn
u_n = \bigotimes_{k = 1}^{n} v_n,
\]
so that for $a \in D_n$ we have $\bt (a) = u_n a u_n^*$.

Let $F \S D$ be finite and let $\ep > 0$.
Choose $m$ and a finite set $S \S D_m$ such that
every element of $F$ is within $\frac{1}{2} \ep$ of an element of $S$.
Choose $n > m$ and so large that $2^{- m + 1} < \ep$.
Define $p = \diag (0, 1, 1, \dots, 1) \in M_{2^{m - 1}}$.
Take
\[
q_0 = \frac{1}{2} \left( \begin{array}{cc}
  p     & p       \\  p     & p
\end{array} \right)
\andeqn
q_1 = \frac{1}{2} \left( \begin{array}{cc}
  p     &  - p        \\   - p     &  p
\end{array} \right),
\]
and for $j = 0, \, 1$ define
\[
e_j = 1_{D_{m - 1}} \otimes q_j \in D_m \S D.
\]
Then each $e_j$ is a \pj, and $e_0 e_1 = 0$.
Moreover, $e_j$ commutes exactly with every element
of $S$ (indeed, with every element of $D_n$), so that
$\| e_j a - a e_j \| < \ep$ for all $a \in F$.

To compute $\bt (e_j)$, we write
\[
v_n = \left( \begin{array}{cc} 1 & 0 \\ 0 & 1 - 2 p \end{array} \right).
\]
Then calculations show that $v_n q_0 v_n^* = q_1$.
Consequently $\bt (e_0) = e_1$.
Finally, the unique tracial state $\ta$ on $D$ restricts to the
usual tracial state on $1_{D_{m - 1}} \otimes M_{2^m} \otimes 1 \S D$,
so
\[
\ta (1 - e_0 - e_1)
  = \frac{2 \cdot \rank (p)}{2^m} = \frac{1}{2^{m - 1}}
\]
satisfies $| \ta (1 - e_0 - e_1) | < \ep$.
This completes the verification of the hypotheses of
Theorem~\ref{ARPFromPosElts2},
so the action generated by $\bt$ has the \aRp.
\end{proof}

\begin{prp}\label{CAR2_2}
The action of $\Zqt$ generated by the automorphism $\bt$
of Example~\ref{CAR2} does not have the strict Rokhlin property.
\end{prp}

\begin{proof}
We show that there is no \pj\  $e \in D$ such that
$\| \bt (e) - (1 - e) \| < 1$.
Let $D_n$, $v_n$, and $u_n$ be as in
the proof of Proposition~\ref{CAR2_1}.
Write $D_n = M_{t (n)}$, where $t (n) = 2^{n (n + 1) / 2}$.
Then $u_n$ is a $t (n) \times t (n)$ diagonal matrix whose
diagonal entries are all either $1$ or $- 1$.
Let $r (n)$ be the number of entries equal to $1$, and
let $s (n)$ be the number of entries equal to $- 1$.
We prove by induction that $r (n) - s (n) = 2^n$.
This is certainly true for $n = 1$, when $t (n) = 2$ and $u_n = 1$.
If it is true for $n$, then
\[
r (n + 1) = (2^n + 1) r (n) + (2^n - 1) s (n)
\]
and
\[
s (n + 1) = (2^n - 1) r (n) + (2^n + 1) s (n),
\]
so
\[
r (n + 1) - s (n + 1) = 2 (r (n) - s (n)) = 2^{n + 1}.
\]
This completes the induction.

Now suppose there is a \pj\  $e \in D$ such that
$\| \bt (e) - (1 - e) \| < 1$.
Set $\ep = \frac{1}{2} ( 1 - \| \bt (e) - (1 - e) \|) > 0$.
Choose $n$ and a \pj\  $f \in D_n$ such that $\| e - f \| < \ep$.
Then $\| u_n f u_n^* - (1 - f) \| < 1$.
It follows that $f \sim 1 - f$.
Now recall that $D_n = M_{t (n)}$, and represent this algebra
on the Hilbert space $\C^{t (n)}$ in the usual way.
Then $\rank (f) = \frac{1}{2} t (n)$.
We can write $u_n = q_0 - q_1$ where $q_0$ and $q_1$ are orthogonal
\pj s of ranks $r (n)$ and $s (n)$.
{}From the previous paragraph, $r (n) > \frac{1}{2} t (n)$,
so $E = q_0 \C^{t (n)} \cap f \C^{t (n)}$ is nontrivial.
Choose $\xi \in E$ with $\| \xi \| = 1$.
Then
\[
(1 - f) \xi = 0, \,\,\,\,\, q_0 \xi = \xi, \andeqn q_1 \xi = 0,
\]
so
\[
[u_n f u_n^* - (1 - f)] \xi = \xi.
\]
It follows that $\| u_n f u_n^* - (1 - f) \| \geq 1$.
This contradiction shows that $e$ does not exist, and that
the action generated by $\bt$ does not have the strict Rokhlin
property.
\end{proof}

\begin{prp}\label{CAR2_3}
Let $\bt \in \Aut (D)$ be as in Example~\ref{CAR2}.
Let ${\widehat{\bt}}$ be the nontrivial automorphism of the
dual action on $\Cs{2}{D}{\bt}$.
Then ${\widehat{\bt}}$ is tracially approximately inner
and generates an action with the \sRp,
but ${\widehat{\bt}}$ is not approximately inner
and is nontrivial on $K_0 (\Cs{2}{D}{\bt} )$.
\end{prp}

\begin{proof}
It follows from Lemma~\ref{Z2Struct}(5) and Proposition~\ref{CAR2_1}
that ${\widehat{\bt}}$ is tracially approximately inner,
and from Lemma~\ref{Z2Struct}(3)
that ${\widehat{\bt}}$ generates an action with the \sRp.

We show that ${\widehat{\bt}}$
is nontrivial on $K_0 (\Cs{2}{D}{\bt} )$;
that ${\widehat{\bt}}$ is not approximately inner follows.
Using Lemma~\ref{Z2Struct}, we can identify
$K_0 (\Cs{2}{D}{\bt} )$ as the direct limit of
the system
\[
\Z^2 \stackrel{T_1}{\longrightarrow}
\Z^2 \stackrel{T_2}{\longrightarrow}
\Z^2 \stackrel{T_3}{\longrightarrow} \cdots,
\]
with
\[
T_n = \left( \begin{array}{cc}  2^{n - 1} + 1  & 2^{n - 1} - 1  \\
               2^{n - 1} - 1  & 2^{n - 1} + 1 \end{array} \right),
\]
and where ${\widehat{\bt}}_*$ is the direct limit of the
maps $(j, k) \mapsto (k, j)$ on $\Z^2$.
One proves by induction that
\[
T_n \circ T_{n - 1} \circ \cdots \circ T_1 (1, \, -1)
  = (2^n, \, - 2^n).
\]
This expression is nonzero for all $n$, so the image of $(1, \, -1)$
is a nonzero element $\et$ of $K_0 (\Cs{2}{D}{\bt} )$.
The element $\et$ is not torsion,
since $\Cs{2}{D}{\bt}$ is an AF algebra.
Evidently ${\widehat{\bt}}_* (\et) = - \et \neq \et$.
\end{proof}

The following example was suggested by Izumi.
Some of the properties given here are folklore,
but we have been unable to find a reference for the proofs.
This example is a special case of an example used for other purposes
in Example~3.14 of~\cite{Iz0}.

\begin{exa}\label{CAR3}
Let $\gm$ be the automorphism of order~$2$ given by
\[
D = \bigotimes_{n = 1}^{\I} M_{2^n} \andeqn
\gm = \bigotimes_{n = 1}^{\I}
 \Ad ( 1_{2^{n} - 1} \oplus (- 1) ).
\]
The automorphism in the $n$-th tensor factor is conjugation by a
diagonal unitary in which $2^{n} - 1$ diagonal entries are equal
to $1$ and one diagonal entry is equal to $- 1$.

We prove the following facts in the next three propositions.
The automorphism $\gm$ is approximately inner,
the action of $\Zqt$ it generates does not have the \aRp,
and the generator of the dual action is not
tracially approximately inner.
Nevertheless, the crossed product is a simple unital AF~algebra,
and the dual action has the strict Rokhlin property.
As an ``explanation'', in the factor representation of $D$
associated to the trace, $\gm$ becomes inner.
\end{exa}

\begin{prp}\label{CAR3_1}
Let $\gm \in \Aut (D)$ be as in Example~\ref{CAR3}.
Then the crossed product $C^* (\Zqt, \, D, \, \gm)$
is a simple AF~algebra which has exactly two extreme tracial states.
The dual action exchanges these two tracial states.
\end{prp}

\begin{proof}
Let $D_n = \bigotimes_{k = 1}^{n} M_{2^k}$, so that
$D_{n} = D_{n - 1} \otimes M_{2^n}$ and $D = \dirlim D_n$.
We identify the direct system of crossed products as in
Lemma~\ref{Z2Struct},
and we follow the notation there, with
\[
e_n = 1_{2^{n} - 1} \oplus 0 \in M_{2^n}
\andeqn f_n = 0_{2^{n} - 1} \oplus 1 \in M_{2^n}.
\]
and with $t (n) = 2^{n (n + 1) / 2}$.
The algebra $C^* (\Zqt, \, D, \, \gm )$ is a simple AF~algebra
by Lemma~\ref{Z2Struct}(1).

We want to identify all tracial states
on $C^* (\Zqt, \, D, \, \gm )$,
but we begin with some convenient notation.
For $\ld \in \R$ define the matrix
\[
T (\ld) = \frac{1}{2} \left( \begin{array}{cc}
           1 + \ld & 1 - \ld \\ 1 - \ld & 1 + \ld \end{array} \right)
        = \left( \begin{array}{cc}
               1  & 1  \\ 1  & -1 \end{array} \right)^{-1}
          \left( \begin{array}{cc}
               1  & 0  \\ 0  & \ld \end{array} \right)
          \left( \begin{array}{cc}
               1  & 1  \\ 1  & -1 \end{array} \right).
\]
{}From the second expression, we see that
$T (\ld \mu) = T (\ld) T (\mu)$ for $\ld, \, \mu \in \R$,
and that $T (1) = 1$.
{}From the first expression, we see that the matrix of partial
embedding multiplicities of $\ps_n$ is exactly
$2^n T \left( 1 - \frac{1}{2^{n - 1}} \right)$.
Moreover, for $\ld \in [0, 1]$, if $(r_0, s_0), \, (r, s) \in \R^2$
satisfy $T (\ld) (r_0, s_0) = (r, s)$,
and if $r_0, s_0 \in [0, 1]$ with $r_0 + s_0 = 1$,
then $r, s \in [0, 1]$ with $r + s = 1$.

Let ${\mathrm{tr}}_m$ denote the normalized trace on $M_m$.
Tracial states on $C^* (\Zqt, \, D, \, \gm )$ are in one to one
correspondence with sequences $(\ta_n)_{n \in \N}$ of tracial states
$\ta_n$ on $M_{t (n)} \oplus M_{t (n)}$ satisfying the
compatibility conditions $\ta_n \circ \ps_n = \ta_{n - 1}$
for all $n$.
The tracial state $\ta_n$ has the form
$\ta_n (a, b)
  = r_n {\mathrm{tr}}_{t (n)} (a) + s_n {\mathrm{tr}}_{t (n)} (b)$
for $r_n, s_n \in [0, 1]$ with $r_n + s_n = 1$,
and the compatibility condition is exactly
$T \left( 1 - \frac{1}{2^{n - 1}} \right) (r_n, s_n)
  = (r_{n - 1}, \, s_{n - 1})$
in $\R^2$.

Define
\[
\ld_n = \prod_{k = n + 1}^{\I} \left( 1 - \frac{1}{2^{k - 1}} \right)
 \in [0, 1).
\]
We claim that a sequence $(r_n, s_n)_{n \in \N}$ corresponds to a
tracial state on the crossed product
$C^* (\Zqt, \, D, \, \gm )$ \ifo\  there is
$r \in [0, 1]$ such that for all $n \in \N$
we have $(r_n, s_n) = T (\ld_n) (r, \, 1 - r)$.
One direction is easy: given $r$, the sequence $(r_n, s_n)_{n \in \N}$
defined by $T (\ld_n) (r, \, 1 - r) = (r_n, s_n)$ clearly satisfies
$r_n, s_n \in [0, 1]$ and $r_n + s_n = 1$,
and the compatibility condition follows from the relation
$\left( 1 - \frac{1}{2^{n - 1}} \right) \ld_n = \ld_{n - 1}$.
For the converse, one observes that $\log (1 - x) \geq - 2 x$
for $0 \leq x \leq \frac{1}{2}$, so that
\[
\log (\ld_n)
  = \sum_{k = n + 1}^{\I} \log \left( 1 - \frac{1}{2^{k - 1}} \right)
  \geq \sum_{k = n + 1}^{\I} 2 \left( 1 - \frac{1}{2^{k - 1}} \right)
  > - \I.
\]
Therefore $\ld_n > 0$, so that $T (\ld_n)$ is invertible.
Define $r, s \in \R$ by $(r, s) = T (\ld_n)^{-1} (r_n, s_n)$.
The compatibility condition guarantees that this definition does not
depend on $n$.
Moreover, $\limi{n} \ld_n = 1$, whence $\limi{n} T (\ld_n)^{-1} = 1$.
Therefore $\limi{n} r_n = r$ and $\limi{n} s_n = s$.
It follows that $r, s \in [0, 1]$ and $r + s = 1$.
This completes the proof of the claim.

Since $T (\ld_n)$ is invertible, we have an affine parametrization
of the tracial states on $\Cs{2}{D}{\gm}$ by $[0, 1]$,
from which it is clear that
there are exactly two extreme tracial states.
That the dual action exchanges them
is clear from the identification of the dual action with the flip
in Lemma~\ref{Z2Struct}.
\end{proof}

\begin{prp}\label{CAR3_2}
Let $\gm \in \Aut (D)$ be as in Example~\ref{CAR3}.
Then:
\begin{itemize}
\item[(1)]
The action of $\Zqt$ generated by $\gm$ does not have the \aRp.
\item[(2)]
The automorphism $\gm$ is approximately inner.
\item[(3)]
The dual action on $\Cs{2}{D}{\gm}$ has
the strict Rokhlin property.
\item[(4)]
The generating automorphism ${\widehat{\gm}}$ of the dual action is not
tracially approximately inner.
\end{itemize}
\end{prp}

\begin{proof}
Proposition~\ref{CAR3_1} shows that ${\widehat{\gm}}$
is not the identity on the tracial state space.
So Lemma~\ref{TAIAndTraces} implies that
${\widehat{\gm}}$ is not tracially approximately inner.
This is~(4).
The rest follows from Lemma~\ref{Z2Struct}.
\end{proof}

\begin{prp}\label{CAR3_3}
Let $\gm \in \Aut (D)$ be as in Example~\ref{CAR3}.
Let $\pi$ be the Gelfand-Naimark-Segal representation
associated with the unique tracial state $\ta$ on $D$.
Then the automorphism ${\overline{\gm}}$ of $\pi (D)''$
induced by $\gm$ is inner.
\end{prp}

\begin{proof}
In a slight modification of the notation used previously, set
\[
e_n^{(0)} = 1_{2^{n} - 1} \oplus 0 \in M_{2^n}
\andeqn f_n^{(0)} = 0_{2^{n} - 1} \oplus 1 \in M_{2^n}.
\]
Then define \pj s $e_n, \, f_n \in D = \bigotimes_{n = 1}^{\I} M_{2^n}$
by
\[
e_n = 1 \otimes \cdots 1 \otimes e_n^{(0)} \otimes 1 \otimes \cdots
\andeqn
f_n = 1 \otimes \cdots 1 \otimes f_n^{(0)} \otimes 1 \otimes \cdots,
\]
where the nontrivial tensor factors are in the $n$-th positions.
Define unitaries in $D$ by $v_n = e_n - f_n$ and
$u_n = \prod_{k = 1}^{n} v_k$.
Further note that $\ta (f_n) = \frac{1}{2^n}$.

We show that $\limi{n} \pi (u_n)$ exists in $\pi (D)''$
in the strong operator topology, and that the limit is a unitary
which implements $\gm$ on $\pi (D)$.
This is easily seen to imply the result.
Let $\xi$ be the standard cyclic vector for the representation $\pi$,
and let $H$ be the Hilbert space on which it acts.
For $a \in D$ we have
\begin{align*}
\| \pi (u_n) \pi (a) \xi - \pi (u_{n - 1}) \pi (a) \xi \|^2
& = \| \pi (v_n) \pi (a) \xi - \pi (a) \xi \|^2  \\
& = \langle \pi (a^* (v_n - 1)^* (v_n - 1) a) \xi, \, \xi \rangle  \\
& = \ta (a^* (v_n - 1)^* (v_n - 1) a)   \\
& = \ta ((v_n - 1) a a^* (v_n - 1)^* )   \\
& \leq \ta ((v_n - 1) (v_n - 1)^* ) \| a \|^2
  = 4 \ta (f_n) \| a \|^2.
\end{align*}
Since $\sum_{n = 1}^{\I} \ta (f_n)^{1/2} < \I$,
it follows that $(\pi (u_n) \pi (a) \xi)_{n \in \N}$
is a Cauchy sequence in $H$, and therefore converges.
The elements $\pi (a) \xi$ form a dense subset of $H$,
and $\sup_{n \in \N} \| \pi (u_n) \| < \I$,
so a standard argument shows that
$u \et = \limi{n} \pi (u_n) \et$ exists for all $\et \in H$.
Further, $u$ is clearly isometric, hence bounded.
Since $u_n$ is selfadjoint and the selfadjoint elements are
strong operator closed in $L (H)$, we get $u^* = u$.
Since multiplication is jointly strong operator continuous on
bounded sets, $u^2 = 1$.
Therefore $u$ is a selfadjoint unitary.

For any $a \in D$, we have
$\limi{n} u_n a u_n = \limi{n} u_n a u_n^* = \gm (a)$ in norm,
and, again using joint strong operator continuity of
multiplication on bounded sets,
$\limi{n} u_n a u_n = u a u = u a u^*$
in the strong operator topology.
Therefore $u a u^* = \gm (a)$, as desired.
\end{proof}

\begin{exa}\label{CAR4}
Let $A$ be the $2^{\infty}$ UHF algebra,
and let $\af$ be the automorphism
constructed in Section~5 of~\cite{Bl0}.
It follows from Corollary 5.3.2 of~\cite{Bl0}
and Takai duality that $\Cs{2}{A}{\af}$ is not an AF~algebra.
We prove in Proposition~\ref{BlAutHasARP} below
that $\af$ generates an action of $\Zqt$ with the \aRp,
and in Proposition~\ref{CAR4_4} below
that this action does not have the strict Rokhlin property,
and that the generator of the dual action
is tracially approximately inner
but induces a nontrivial automorphism of $K_1$.
By construction, the action of $\Zqt$ generated by $\af$
is approximately representable
in the sense of Definition~3.6(2) of~\cite{Iz}.
(The construction is recalled below.)
It follows from Lemma~3.8(2) of~\cite{Iz}
that the dual action on $\Cs{2}{A}{\af}$ has
the strict Rokhlin property.
\end{exa}

We remark that the methods used to prove the \aRp\  in this example
seem likely to be more
typical of proofs that actions on AH algebras have the \aRp\  than
the methods used for
Proposition~\ref{IrratAPR} and Proposition~\ref{FTHasARP}.

We begin with a convenient description
of the construction in~\cite{Bl0}, following Section~5 there.
We make the convention that in any block matrix decomposition,
all blocks are to be the same size,
and we write $1_{n}$ for the identity of $M_n$.
We take the identification of $M_m \otimes M_n$ with $M_{m n}$
to send $a \otimes e_{j, k}$
to the $n \times n$ block matrix with $m \times m$ blocks,
of which the $(j, k)$ block is $a$ and the rest are zero.
Thus, $a \otimes 1 = \diag (a, a, \dots, a)$.
Also, we identify the circle $S^1$ with $\R / \Z$,
and write elements of $C (S^1)$ as functions on $[0, 1]$
whose values at $0$ and $1$ are equal.

Following Definition 3.1.1 of~\cite{Bl0},
we choose a standard twice around embedding
$\ph^+ \colon C (S^1) \to C (S^1, M_2)$,
given by choosing a \ct\  unitary path $c \in C ([0, 1], \, M_2)$
with
\[
c (0) = 1 \andeqn
c (1)
 = \left( \begin{array}{cc} 0 & 1 \\ 1 & 0 \end{array} \right),
\]
and then setting
\[
\ph^+ (f) (t) = c (t)
  \left( \begin{array}{cc} f \left( \frac{1}{2} t \right) & 0
    \\ 0 & f \left( \frac{1}{2} (t + 1) \right) \end{array} \right)
   c (t)^*
\]
for $f \in C (S^1)$.
Further let $\ph^- \colon C (S^1) \to C (S^1, M_2)$
be the standard $- 2$ times around embedding
$\ph^- (f) (t) = \ph^+ (f) (1 - t)$.
Extend everything, using the same notation, to embeddings
of $C (S^1, M_m)$ in $C (S^1, M_{2 m})$,
by using $\ph^+ \otimes \id_{M_m}$, etc.

Following Section~5 of~\cite{Bl0},
set $A_n = C (S^1, \, M_{4^n})$
and, remembering our convention on block sizes,
define a unitary in $M_{4^n} \S A_n$ by $u_n = \diag (1, \, -1)$.
Further define $\ps_n \colon A_n \to A_{n + 1}$ by
\[
\ps_n \left( \begin{array}{cc} f_{1, 1} & f_{1, 2}
    \\ f_{2, 1} & f_{2, 2} \end{array} \right)
 = \left( \begin{array}{cccc}
\ph^+ (f_{1, 1}) &  0     & \ph^+ (f_{1, 2}) &  0      \\
  0     & \ph^- (f_{2, 2}) &  0     & \ph^- (f_{2, 1}) \\
\ph^+ (f_{2, 1}) &  0     & \ph^+ (f_{2, 2}) &  0      \\
  0     & \ph^- (f_{1, 2}) &  0     & \ph^- (f_{1, 1})
\end{array} \right).
\]
Theorem 4.1.1, Proposition 5.1.1, and Proposition 5.1.2 of~\cite{Bl0}
show that the direct limit of the $A_n$ using the maps $\ps_n$ is
the $2^{\infty}$ UHF algebra $A$,
and that the automorphisms $\af_n = \Ad (u_n)$ of $A_n$ define an
automorphism $\af$ of $A$ of order two.
It follows from Takai duality and
Proposition 5.2.2 and Proposition 5.4.1 of~\cite{Bl0}
that $\Cs{2}{A}{\af}$ is isomorphic to the tensor product
of $A$ and the $2^{\infty}$ Bunce-Deddens algebra.

Further let $\io_n \colon M_{4^n} \to A_n$ be the embedding of
matrices as constant functions,
and define $\sm_n \colon M_{4^n} \to M_{4^{n + 1}}$ by
\[
\io_n \left( \begin{array}{cc} a_{1, 1} & a_{1, 2}
    \\ a_{2, 1} & a_{2, 2} \end{array} \right)
 = \left( \begin{array}{cccc}
a_{1, 1} \otimes 1_2 &  0     & a_{1, 2} \otimes 1_2 &  0      \\
  0     & a_{2, 2} \otimes 1_2 &  0     & a_{2, 1} \otimes 1_2 \\
a_{2, 1} \otimes 1_2 &  0     & a_{2, 2} \otimes 1_2 &  0      \\
  0     & a_{1, 2} \otimes 1_2 &  0     & a_{1, 1} \otimes 1_2
\end{array} \right).
\]
It is a consequence of the next lemma that
$\io_{n + 1} \circ \sm_n = \ps_n \circ \io_n$.
Moreover, we get an automorphism $\mu_n$ of $M_{4^n}$
by defining
$\mu_n = \Ad (u_n)$, and $\io_n \circ \mu_n = \af_n \circ \io_n$.

\begin{lem}\label{PresConst}
Let the notation be as above.
Let $\ep \geq 0$, let $f \in A_n$, and let $a \in M_{4^n}$.
Suppose $f (t) = a$ for all $t \in [\ep, \, 1 - \ep]$.
Then $\ps_n (f)(t) = \sm_n (a)$
for all $t \in [2 \ep, \, 1 - 2 \ep]$.
In particular, $\io_{n + 1} \circ \sm_n = \ps_n \circ \io_n$.
\end{lem}

\begin{proof}
Write
\[
f = \left( \begin{array}{cc} f_{1, 1} & f_{1, 2}
    \\ f_{2, 1} & f_{2, 2} \end{array} \right)
\andeqn
a = \left( \begin{array}{cc} a_{1, 1} & a_{1, 2}
    \\ a_{2, 1} & a_{2, 2} \end{array} \right).
\]
Then for each $j$ and $k$,
we have $f_{j, k} (t) = a_{j, k}$ for all $t \in [\ep, \, 1 - \ep]$.
For $t \in [2 \ep, \, 1 - 2 \ep]$ we therefore get
\[
\left( \begin{array}{cc} f_{j, k} \left( \frac{1}{2} t \right) & 0
    \\ 0 & f_{j, k} \left( \frac{1}{2} (t + 1) \right) \end{array} \right)
= \left( \begin{array}{cc} a_{j, k} & 0
    \\ 0 & a_{j, k} \end{array} \right),
\]
which commutes with $c (t)$.
\end{proof}

In the next lemma we show, roughly, that whenever an element
$f \in A_n = C (S^1, M_{4^n})$ is unitarily equivalent
in $C ([0, 1], \, M_{4^n})$, via invariant unitaries,
to a function with small variation over intervals of length $\dt$,
then $\ps_n (f)$ is unitarily equivalent
in $C ([0, 1], \, M_{4^{n + 1}})$, again via invariant unitaries,
to a function with small variation over intervals of length $2 \dt$.

\begin{lem}\label{AlmostConst}
Let the notation be as above.
Let $t \mapsto x (t)$ be a unitary element of $C ([0, 1], \, M_{4^n})$
such that $\mu_n (x (t)) = x (t)$ for every $t \in [0, 1]$.
Then there exists
a unitary element $t \mapsto y (t)$ of $C ([0, 1], \, M_{4^{n + 1}})$
such that $\mu_{n + 1} (y (t)) = y (t)$ for every $t \in [0, 1]$,
with the property that whenever
$\ep > 0$, $\dt > 0$, and $f \in A_n$ satisfy
\[
\| x (s) f (s) x (s)^* - x (t) f (t) x (t)^* \| < \ep
\]
for all $s, \, t \in [0, 1]$ such that $| s - t | < \dt$,
then
\[
\| y (s) \ps_n (f) (s) y (s)^* - y (t) \ps_n (f) (t) y (t)^* \| < \ep
\]
for all $s, \, t \in [0, 1]$ such that $| s - t | < 2 \dt$.
\end{lem}

\begin{proof}
The equation $\mu_n (x (t)) = x (t)$ implies that we can write
$x (t) = x_1 (t) \oplus x_2 (t)$ for unitaries
$x_1, \, x_2 \in C ([0, 1], \, M_{2^{2 n - 1}})$.
For $j = 1, 2$ define
\[
y_j (t) =
  \left( \begin{array}{cc} x_j \left( \frac{1}{2} t \right) & 0
    \\ 0 & x_j \left( \frac{1}{2} (t + 1) \right) \end{array} \right)
   c (t)^*.
\]
Then define
\[
y (t) = \diag ( y_1 (t), \, y_2 (1 - t), \, y_2 (t), \, y_1 (1 - t) ).
\]
Evidently $y$ is a unitary in $C ([0, 1], \, M_{4^{n + 1}})$
and $\mu_{n + 1} (y (t)) = y (t)$ for every $t \in [0, 1]$.

To verify the conclusions of the lemma, it will simplify the
notation to conjugate everything by the permutation matrix
\[
w = \left( \begin{array}{cccc}
  1     &  0     &  0     &  0      \\
  0     &  0     &  1     &  0      \\
  0     &  0     &  0     &  1      \\
  0     &  1     &  0     &  0
\end{array} \right).
\]
(This conjugation is also used in Section~5 of~\cite{Bl0}.)
Thus, let
\[
{\widetilde{\ps}} (f) = w \ps_n (f) w^*
 = \left( \begin{array}{cccc}
\ph^+ (f_{1, 1}) & \ph^+ (f_{1, 2}) &  0     &  0      \\
\ph^+ (f_{2, 1}) & \ph^+ (f_{2, 2}) &  0     &  0      \\
  0     &  0     & \ph^- (f_{1, 1}) & \ph^- (f_{1, 2}) \\
  0     &  0     & \ph^- (f_{2, 1}) & \ph^- (f_{2, 2})
\end{array} \right),
\]
let
\[
{\widetilde{u}} = w u_{n + 1} w^* = \diag (1, \, -1, \, -1, \, 1),
\]
let
\[
{\widetilde{y}} (t) = w y (t) w^* 
 = \diag ( y_1 (t), \, y_2 (t), \, y_1 (1 - t), \, y_2 (1 - t) ),
\]
and similarly define ${\widetilde{\mu}}$, etc.
Note that ${\widetilde{\io}} = \io_{n + 1}$.

Let $\ep > 0$, $\dt > 0$, and $f \in A_n$ be as in the hypotheses.
Define
\[
g (t)
 = \left( \begin{array}{cc}
y_1 (t) \ph^+ (f_{1, 1}) (t) y_1 (t)^*
             & y_1 (t) \ph^+ (f_{1, 2}) (t) y_2 (t)^*  \\
y_2 (t) \ph^+ (f_{2, 1}) (t) y_1 (t)^*
             & y_2 (t) \ph^+ (f_{2, 2}) (t) y_2 (t)^*
\end{array} \right),
\]
and note that
\[
{\widetilde{y}} (t) {\widetilde{\ps}} (f) {\widetilde{y}} (t)^*
  = \diag (g (t), \, g (1 - t)).
\]
Accordingly, it suffices to prove that if
\[
\| x (s) f (s) x (s)^* - x (t) f (t) x (t)^* \| < \ep
\]
for all $s, \, t \in [0, 1]$ such that $| s - t | < \dt$,
then
\[
\| g (s) - g (t) \| < \ep
\]
for all $s, \, t \in [0, 1]$ such that $| s - t | < 2 \dt$.

Let $v$ be the permutation matrix
\[
v = \left( \begin{array}{cccc}
  1     &  0     &  0     &  0      \\
  0     &  0     &  1     &  0      \\
  0     &  1     &  0     &  0      \\
  0     &  0     &  0     &  1
\end{array} \right).
\]
When one calculates $v g (t) v^*$ by
substituting the formulas for $y_j (t)$ and $\ph^+$ in the expression
for $g (t)$,
the factors $c (t)$ and $c (t)^*$ all cancel out,
and the final answer is
\[
v g (t) v^* = \ts{
  \diag \left(
x \left( \frac{1}{2} t \right) f \left( \frac{1}{2} t \right)
             x \left( \frac{1}{2} t \right)^*, \,\,
x \left( \frac{1}{2} (t + 1) \right)
             f \left( \frac{1}{2} (t + 1) \right)
                       x \left( \frac{1}{2} (t + 1) \right)^* \right). }
\]
Since we are assuming
\[
\| x (s) f (s) x (s)^* - x (t) f (t) x (t)^* \| < \ep
\]
for all $s, \, t \in [0, 1]$ such that $| s - t | < \dt$,
it is immediate that
$| s - t | < 2 \dt$ implies
\[
\| v g (s) v^* - v g (t) v^* \| < \ep,
\]
whence also
\[
\| g (s) - g (t) \| < \ep,
\]
as desired.
\end{proof}

\begin{prp}\label{BlAutHasARP}
The automorphism $\af$ of Example~\ref{CAR4}
generates an action of $\Zqt$ with the \aRp.
\end{prp}

\begin{proof}
Let the notation be as before Lemma~\ref{PresConst}.
Let $\ta$ be the unique tracial state on $A = \dirlim A_n$.
Define a tracial state $\ta_n$ on $A_n$ by
\[
\ta_n (f) = \int_0^1 {\mathrm{tr}}_{4^n} (f (t ) ) \, d t,
\]
where ${\mathrm{tr}}_{m}$ is the normalized trace on $M_m$.
Then one checks that $\ta_{n + 1} \circ \ps_n = \ta_n$ for all $n$.
It follows from the uniqueness of $\ta$ that $\ta |_{A_n} = \ta_n$
for all $n$.

We use Theorem~\ref{ARPFromPosElts2} to verify the \aRp.
So let $F \S A$ be finite and let $\ep > 0$.
Choose $m$ and a finite set $S_0 \S A_m$ such that
every element of $F$ is within $\frac{1}{8} \ep$ of an element of $S_0$.

The set $S_0$ is a uniformly equicontinuous set of functions from
$[0, 1]$ to $M_{4^m}$, so there is $\dt > 0$ such that
whenever $s, \, t \in [0, 1]$ satisfy $| s - t | < \dt$, then
\[
\| f (s) - f (t) \| < \ts{ \frac{1}{8} } \ep
\]
for all $t \in [0, 1]$ and all $f \in S_0$.
Choose $n \in \N$ with $n \geq m$ and so large that $2^{n - m} \dt > 1$.
Apply Lemma~\ref{AlmostConst} a total of $n - m$ times,
the first time with $x (t) = 1$ for all $t$, 
obtaining after the last application
a \ct\  unitary path $t \mapsto z (t)$ in $C ([0, 1], \, M_{4^n})$
such that $\mu_{n} (z (t)) = z (t)$ for every $t \in [0, 1]$.
Replacing $z (t)$ by $z (0)^* z (t)$,
we may clearly assume that $z (0) = 1$.
Then, in particular,
\[
\| z (t)^* f (0) z (t) - f (t) \| < \ts{ \frac{1}{8} } \ep
\]
for all $t \in [0, 1]$ and all $f \in S_0$.
Recall that we identify $C (S^1, B)$ with the set of
functions $f \in C ([0, 1], \, B)$ such that $f (0) = f (1)$.
Since the fixed point algebra $A_n^{\af_n} = C (S^1, M_{4^n})^{\af_n}$
is just $C (S^1, \, M_{4^n}^{\mu_n})$,
and since $M_{4^n}^{\mu_n}$ is \fd,
there is an $\af_n$-invariant unitary $y \in A_n$ such that
$y (t) = z (t)$ for $t \in \left[0, \, 1 - \frac{1}{8} \ep \right]$.
Then for each $f \in S_0$, regarded as a subset of $A_n$, there exists
$g \in A_n$ such that $\| y^* g y - f \| < \ts{ \frac{1}{4} } \ep$
and $g (t) = f (0)$
for $t \in \left[0, \, 1 - \frac{1}{8} \ep \right]$.
Let $S$ be the set of all elements $g$ obtained in this way from
elements of $S_0$.
In particular, for every $a \in F$ there is $g \in S$ such that
$\| a - y^* g y \| < \frac{1}{2} \ep$.

We claim that there are orthogonal positive elements
$b_0, \, b_1 \in A_{n + 1} \S A$ such that
$b_j g = g b_j$ for all $g \in S$,
and such that
\[
\af_{n + 1} (b_0) = b_1, \,\,\,\,\,\,
\af_{n + 1} (b_1) = b_0, \,\,\,\,\,\,
0 \leq b_0, \, b_1 \leq 1, \andeqn
0 \leq \ta (1 - b_0 - b_1) < \ep.
\]
For this purpose, it suffices to use in place of $\ps_n$ the
unitarily equivalent \hm\  %
\[
{\widetilde{\ps}} = w \ps_n (-) w^*
  \colon C (S^1, M_{4^n}) \to C (S^1, M_{4^{n + 1}})
\]
in the proof of Lemma~\ref{AlmostConst}
(called $\om_n$ in the proof of Proposition~5.1.1 of~\cite{Bl0}),
and to use in place of $\af_{n + 1}$ the automorphism
\[
{\widetilde{\af}} = \Ad (w) \circ \af_{n + 1} = \Ad ({\widetilde{u}})
 = \Ad (\diag (1, \, -1, \, -1, \, 1) ).
\]
Note that this change does not require any change in the formula
for the trace $\ta_{n + 1}$,
and also does not affect the first part of
the conclusion of Lemma~\ref{PresConst}.
Accordingly, if
\[
g = \left( \begin{array}{cc} g_{1, 1} & g_{1, 2}
    \\ g_{2, 1} & g_{2, 2} \end{array} \right)
  \in S
\]
then ${\widetilde{\ps}} (g) \in  C (S^1, M_{4^{n + 1}})$ satisfies
\[
{\widetilde{\ps}} (g) (t) = \left( \begin{array}{cccc}
g_{1, 1} (0) & g_{1, 2} (0) &  0     &  0      \\
g_{2, 1} (0) & g_{2, 2} (0) &  0     &  0      \\
  0     &  0     & g_{1, 1} (0) & g_{1, 2} (0) \\
  0     &  0     & g_{2, 1} (0) & g_{2, 2} (0)
\end{array} \right)
 = \left( \begin{array}{cc} g (0) & 0 \\ 0 & g (0) \end{array} \right)
\]
for $t \in \left[\frac{1}{4} \ep, \, 1 - \frac{1}{4} \ep \right]$.

Now set
\[
p_0 = \frac{1}{2} \left( \begin{array}{cccc}
  1     &  0     &  1     &  0      \\
  0     &  1     &  0     &  1      \\
  1     &  0     &  1     &  0      \\
  0     &  1     &  0     &  1
\end{array} \right)
\andeqn
p_1 = \frac{1}{2} \left( \begin{array}{cccc}
  1     &  0     & - 1    &  0      \\
  0     &  1     &  0     & - 1     \\
 - 1    &  0     &  1     &  0      \\
  0     & - 1    &  0     &  1
\end{array} \right),
\]
both in $M_4 (M_{4^n})$.
In $2 \times 2$ block form, we can write
\[
{\widetilde{u}}
 = \left( \begin{array}{cc} s & 0 \\ 0 & - s \end{array} \right)
\,\,\,\,\,\, {\mbox{with}} \,\,\,\,\,\,
s = \left( \begin{array}{cc} 1 & 0 \\ 0 & - 1 \end{array} \right),
\]
\[
p_0 = \frac{1}{2}
  \left( \begin{array}{cc} 1 & 1 \\ 1 & 1 \end{array} \right),
\andeqn
p_1 = \frac{1}{2}
  \left( \begin{array}{cc} 1 & - 1 \\ - 1 & 1 \end{array} \right).
\]
With these formulas, it is easy to check that
$p_0$ and $p_1$ are \pj s with $p_0 + p_1 = 1$,
that $p_0$ and $p_1$ commute with ${\widetilde{\ps}} (g) (t)$
for every $g \in S$
and $t \in \left[ \frac{1}{4} \ep, \, 1 - \frac{1}{4} \ep \right]$,
and that $\Ad ({\widetilde{u}})$ exchanges $p_0$ and $p_1$.

Now choose and fix a \cfn\  $h \colon [0, 1] \to [0, 1]$
such that $h (t) = 0$ for
$t \not\in \left[ \frac{1}{4} \ep, \, 1 - \frac{1}{4} \ep \right]$
and $h (t) = 1$ for
$t \in \left[ \frac{1}{2} \ep, \, 1 - \frac{1}{2} \ep \right]$,
and define $b_j (t) = h (t) p_j$ for $j = 0, \, 1$.
Then $b_0$ and $b_1$ are positive elements with $b_0, \, b_1 \leq 1$,
which commute with ${\widetilde{\ps}} (g)$ for every $g \in S$,
which satisfy $b_0 b_1 = 0$,
such that $\Ad ({\widetilde{u}})$ exchanges $b_0$ and $b_1$,
and such that $0 \leq \ta_{n + 1} (1 - b_0 - b_1) < \ep$.
This proves the claim above.

We return to the use of $\ps_{n + 1}$,
and we let $b_0, \, b_1 \in A_{n + 1} \S A$ be as in the claim
(rather than its proof).
In $A_{n + 1} \S A$, define $a_0 = y^* b_0 y$ and $a_1 = y^* b_1 y$.
Since $\af (y) = y$, it follows that
$a_0, \, a_1 \in A_{n + 1} \S A$ satisfy
$a_j y^* g y = y^* g y a_j$ for all $g \in S$,
and
\[
a_0 a_1 = 0, \,\,\,\,\,\,
\af_{n + 1} (a_0) = a_1, \,\,\,\,\,\,
0 \leq a_0, \, a_1 \leq 1, \andeqn
0 \leq \ta (1 - a_0 - a_1) < \ep.
\]
For $a \in F$ choose $g \in S$ such that
$\| a - y^* g y \| < \frac{1}{4} \ep$.
Then
\[
\| [a_j, a ] \| \leq 2 \| a - y^* g y \| + \| [a_j, \, y^* g y] \|
  < 2 \left( \ts{\frac{1}{2} \ep} \right) + 0 = \ep.
\]
This completes the verification of the hypotheses of
Theorem~\ref{ARPFromPosElts2}, so it follows that $\af$ has the \aRp.
\end{proof}

\begin{prp}\label{CAR4_4}
Let $\af \in \Aut (A)$ be as in Example~\ref{CAR4}.
Then:
\begin{itemize}
\item[(1)]
The action of $\Zqt$ generated by $\af$
does not have the strict Rokhlin property.
\item[(2)]
The dual action on $\Cs{2}{A}{\af}$ has
the strict Rokhlin property.
\item[(3)]
The generating automorphism ${\widehat{\af}}$ of the dual action is
tracially approximately inner.
\item[(4)]
The automorphism ${\widehat{\af}}$
acts nontrivially on $K_1 (\Cs{2}{A}{\af})$.
\end{itemize}
\end{prp}

\begin{proof}
We have already observed in
Example~\ref{CAR4} that $\Cs{2}{A}{\af}$ is not AF.
Therefore~(1) follows from Theorem~\ref{SRokhAF}.

It is immediate from the discussion following
Example~\ref{CAR4} that
the action of $\Zqt$ generated by $\af$
is approximately representable
in the sense of Definition~3.6(2) of~\cite{Iz}.
Part~(2) therefore follows from Lemma~3.8(2) of~\cite{Iz}.

We get~(3) from Theorem~\ref{DualHasRokh}.
The algebra $A$ has cancellation of \pj s because it is AF.
It follows from Proposition~\ref{BlAutHasARP} and
Theorem~\ref{RokhTAF}
that $\Cs{2}{A}{\af}$ has tracial rank zero.
So projections in $\Cs{2}{A}{\af}$ have
the weak divisibility property (Definition~\ref{WDivDfn})
by Lemma~\ref{GetWDiv}.
Therefore Theorem~\ref{DualHasRokh} applies.

It remains to prove~(4).
We continue to follow the notation introduced after
Example~\ref{CAR4}.
Let $B_n$ be the fixed point algebra
\[
A_n^{\af_n}
 = C (S^1, \, M_{2^{2 n - 1}}) \oplus C (S^1, \, M_{2^{2 n - 1}})
 \S C (S^1, \, M_{4^{n}}),
\]
with the embedding being as $2 \times 2$ block diagonal matrices.
Let $B = \dirlim B_n$, which is also equal to $A^{\af}$.
Let $\bt_n \in \Aut (B_n)$ be $\bt_n (f, g) = (g, f)$.
By Proposition~5.2.2 of~\cite{Bl0} and the preceding discussion, 
there is a corresponding automorphism $\bt$ of the direct limit,
$A \cong \Cs{2}{B}{\bt}$,
and the isomorphism can be chosen so that $\af$ generates the
dual action.
By Takai duality, it therefore suffices to show that
$\bt$ is nontrivial on $K_1 (B)$.

Following the discussion after Corollary~5.3.2 of~\cite{Bl0},
let $v \in B_1 \S A_1 = C (S^1, M_4)$ be the unitary
\[
v (t) = \diag \left( e^{2 \pi i t}, \, e^{2 \pi i t},
   \, e^{- 2 \pi i t}, \, e^{- 2 \pi i t} \right).
\]
As there, the image of $[v]$ in $K_1 (B)$ is nonzero.
It follows from Proposition~5.3.1 of~\cite{Bl0}
that $K_1 (B)$ is torsion free,
and one checks that $[\bt_1 (v)] = - [v]$,
so $\bt_* ([v]) = - [v] \neq [v]$.
\end{proof}

\section{Questions}\label{Sec:Qst}

\indent
In this section, we state some open questions.
The first one is suggested by attempts to weaken the
hypotheses of Theorem~\ref{DualHasRokh}.
As far as we know, even the version for AF~algebras and
the strict Rokhlin property is open,
so we state that as our second question.

\begin{qst}\label{Question2}
Let $\af$ be an action of $\Zqn$ on a unital \ca\  $A$.
Suppose that
both $A$ and $\CZnAa$ are simple and
have tracial rank zero, and that the automorphisms $\af_g$
and ${\widehat{\af}}_{\ta}$ are all tracially approximately
inner.
Does it follow that $\af$ has the tracial Rokhlin property?
\end{qst}

\begin{qst}\label{Question3}
Let $\af$ be an action of $\Zqn$ on a unital \ca\  $A$.
Suppose that
both $A$ and $\CZnAa$ are simple and~AF,
and that the automorphisms $\af_g$
and ${\widehat{\af}}_{\ta}$ are all approximately inner.
Does it follow that $\af$ has the Rokhlin property?
\end{qst}

The \aRp, as we have defined it,
seems not to be useful much beyond the class of simple \ca s
with tracial rank zero.
Therefore the following question seems important.

\begin{qst}\label{Question5}
What is the correct definition of the \aRp\  for actions of
$\Zqn$ on \ca s which are not stably finite,
which are stably finite but have badly behaved K-theory,
which are not simple, or which have few or no nontrivial \pj s?
\end{qst}

In particular,
it seems of interest to have a suitable version of the \aRp\  %
for actions on purely infinite simple \ca s,
since results in Section~3 of~\cite{Iz2} give strong restrictions
on the K-theory of actions with the strict Rokhlin property.

Similarly, we want to to modify the definition to make it
suitable for actions on \ca s which are not simple.
One possibility is to add to Definition~\ref{ARPDfn} a condition
requiring that for a prespecified nonzero positive element $a$ in the
finite set $F$,
one has $\| e_0 a e_0 \| > \| a \| - \ep$.
Compare with Definition~2.1 of~\cite{LnTAF}.
We do not know if this strengthening is adequate.

Another case of obvious interest is actions on \ca s without
a reasonable supply of \pj s.
A test for any proposed definition is that it should imply our
definition when there are in fact enough \pj s.
This reasoning suggests the following question.

\begin{qst}\label{Question5a}
Is there an analog of Lemma~\ref{ARPFromPosElts} or
Theorem~\ref{ARPFromPosElts2} for actions of $\Zqn$ on
\suca s with tracial rank zero but more than one tracial state?
\end{qst}

One might try simply writing some version of Definition~\ref{ARPDfn}
which uses positive elements in place of \pj s.

We next turn to tracial approximate innerness.
Our first question is motivated by
Lemmas~\ref{TAIAndInf} and~\ref{TAIAndTraces}.

\begin{qst}\label{Question4}
Let $A$ be a simple separable unital nuclear \ca\  with tracial
rank zero and satisfying the \uct.
Let $\af \in \Aut (A)$.
Suppose that $\af_*$ is the identity on the quotient of $K_0 (A)$
by the infinitesimal subgroup ${\mathrm{Inf}} ( K_0 (A))$.
(Notation from~\cite{GPS}.)
Suppose further that for every
$n \in \N$ there exists $\et \in K_0 (A)_+ \SM \{ 0 \}$
such that $\af_* (\et) = \et$ and $n \et \leq [1_A]$.
Does it follow that $\af$ is tracially approximately inner?
\end{qst}

Since the definition of tracial approximate innerness
suffers from the same difficulties as the definition of the \aRp,
we also ask:

\begin{qst}\label{Question6}
What is the correct definition of tracial approximate innerness of
of automorphisms of \ca s which are not stably finite,
or which are stably finite but have badly behaved K-theory?
\end{qst}

Work of H.\  Lin (in preparation) suggests that
Question~\ref{Question4} has a positive answer,
and that our definition is
at least close to the correct one for automorphisms of
\suca s with tracial rank zero and satisfying the \uct.
The condition of approximate invariance of the \pj,
Condition~(1) of Definition~\ref{TAInnDfn},
may only be appropriate for automorphisms of finite order.
Here also, one might try adding for more general
algebras a condition such as
$\| e a e \| > \| a \| - \ep$
for a prespecified nonzero positive element $a$ in the finite set $F$.

One test for a good definition is that the composition of two
tracially approximately inner automorphisms should again be
tracially approximately inner.
We don't know if this is true for our definition even on
\suca s with tracial rank zero.
If we omit approximate invariance of $e$ from the definition,
then this should be true, and not hard to prove, whenever
the order on \pj s over the algebra is determined by traces
(Definition~\ref{OrdDetD}).
Again, a potentially particularly interesting case
is automorphisms of purely infinite simple \ca s.

\end{document}